\documentclass[10pt,reqno]{amsart}
\include{formatAndDefs}
\usepackage{amsfonts,amssymb,amsmath}

\parskip        \baselineskip
\usepackage{enumitem}
\usepackage{aliascnt}
\usepackage[mathscr]{eucal}
\usepackage{graphicx}
\usepackage{latexsym}
\usepackage{psfrag}
\usepackage{epsfig}
\usepackage[utf8]{inputenc}
\usepackage[english]{babel}
\usepackage{tikz}
\usepackage[foot]{amsaddr}
\usepackage{lineno, blindtext}
\usetikzlibrary{decorations.pathreplacing}
\newcommand{\suchthat}{\;\ifnum\currentgrouptype=16 \middle\fi|\;}
\def\eR{\mathbb{R}}
\def\eN{\mathbb{N}}
\renewcommand{\div}{\mbox{div}\,}

\newcommand{\ds}{\mathrm{d}s}
\newcommand{\dt}{\mathrm{d}t}
\newcommand{\dvr}{\operatorname{div}}
\newcommand{\supp}{\operatorname{supp}}
\newcommand{\dist}{\operatorname{dist}}
\newcommand{\tder}{\partial_t}
\newcommand{\Tr}{\operatorname{tr}}
\newcommand{\GG}{\mathbb{G}}
\newcommand{\RR}{\mathbb{R}}
\newcommand{\mA}{\mathcal{A}}
\newcommand{\ba}{\bf{a}}
\newcommand{\mF}{\mathcal{F}}
\newcommand{\weta}{\widetilde{\eta}}
\newcommand{\ww}{\widetilde{w}}

\newcommand{\str}{L^{\infty}(0,T;W^{3,2}(\Gamma))\cap W^{1,\infty}(0,T;L^{2}(\Gamma))}

\numberwithin{equation}{section}
\newtheorem{thm}{Theorem}
\numberwithin{thm}{section}

\newaliascnt{lemma}{thm}
\newtheorem{lem}[lemma]{Lemma}
\aliascntresetthe{lemma}

\newaliascnt{proposition}{thm}
\newtheorem{prop}[proposition]{Proposition}
\aliascntresetthe{proposition}

\newaliascnt{corollary}{thm}

\aliascntresetthe{corollary}

\newaliascnt{definition}{thm}
\newtheorem{mydef}[definition]{Definition}
\aliascntresetthe{definition}

\newaliascnt{remark}{thm}
\newtheorem{remark}[remark]{Remark}
\aliascntresetthe{remark}

\usepackage[colorlinks=true]{hyperref}

\newcommand{\boundellipse}[3]% center, xdim, ydim
{(#1) ellipse (#2 and #3)
}

\begin{document}

	\title[multi component compressible fluid structure interaction problem]{Existence of  weak solution for a compressible multicomponent fluid structure interaction problem }
	
	\date{\today}
	\author{Martin Kalousek$^\dagger$, Sourav Mitra$^\dagger$, \v{S\'arka Ne\v{c}asov\'a}$^\dagger$}
	\address{$^\dagger$Institute of Mathematics, Czech Academy of Sciences, \v{Z}itn\'a 25, 11567 Prague, Czech Republic}
	
	\begin{abstract}
		We analyze a system of PDEs governing the interaction between two compressible mutually noninteracting fluids and a shell of Koiter type encompassing a time dependent 3D domain filled by the fluids. The dynamics of the fluids is modelled by a system resembling compressible Navier-Stokes equations with a physically realistic pressure depending on densities of both the fluids. The shell possesses a non-linear, non-convex Koiter energy.  Considering that the densities are comparable initially we prove the existence of a weak solution until the degeneracy of the energy or the self-intersection of the structure occurs for two cases. In the first case the adiabatic exponents are assumed to solve $\max\{\gamma, \beta\}> 2$, $\min\{\gamma,\beta\}>0,$ and the structure involved is assumed to be non-dissipative. For the second case we assume the critical case $\max\{\gamma,\beta\}\geq 2$ and $\min\{\gamma,\beta\}>0$ and the dissipativity of the structure. The result is achieved in several steps involving, extension of the physical domain, penalization of the interface condition, artificial regularization of the shell energy and the pressure, the almost compactness argument, added structural dissipation and suitable limit passages depending on uniform estimates.
	\end{abstract}
	
	\maketitle
	\section{Introduction}

	Let us first introduce a few notations corresponding to the fluid structure interaction problem. We consider at time $t$ a domain $\Omega_{\eta}(t)\subset\mathbb{R}^{3}$ and a mixture of two compressible fluids confined in it with a nonlinear elastic Koiter shell appearing at the boundary that interacts with the mixture. 
	
	We denote by $\nu_{\eta}$ the unit outward normal to $\Sigma_{\eta}=\partial\Omega_{\eta}.$ \\
	We first consider a reference domain $\Omega\subset \mathbb{R}^{3},$ whose boundary $\partial\Omega$ is parametrized by a $C^{4}$ injective mapping $\varphi:\Gamma\rightarrow \mathbb{R}^{3},$ where $\Gamma\subset \mathbb{R}^{2}.$ More elaborately we first fix $\Gamma,$ a two dimensional surface, which in our case corresponds to the flat middle surface of the shell and is identified with a torus. We make a simplifying assumption that the plate moves in the normal direction to a reference configuration. The time dependent fluid boundary $\Sigma_{\eta}$ can be characterized by an injective mapping $\varphi_{\eta}$ such that for all pairs $x=(x_{1},x_{2})\in \Gamma,$ the pair $\partial_{i}\varphi_\eta(x),$ $i=1,2,$ is linearly independent. More precisely, $\varphi_{\eta}$ is defined as follows
	\begin{equation}\label{varphieta}
		\varphi_{\eta}(t,x)=\varphi(x)+\eta(x,t)\nu(\varphi(x)) \text{ for } x\in \Gamma,
	\end{equation} 
	where $$\nu(y)=\frac{\partial_{1}\varphi(y)\times\partial_{2}\varphi(y)}{|{\partial_{1}\varphi(y)\times\partial_{2}\varphi(y)}|}$$ is the well defined unit normal to $\partial\Omega=\varphi(\Gamma)$ at $y=\varphi(x)$ and the displacement $\eta:\Gamma\rightarrow \mathbb{R}$ solves a nonlinear plate equation. In other words the time dependent surface $\Sigma_{\eta}$ at any instant $t$ can be expressed as
	\begin{equation}\label{Sigmat}
		\Sigma_{\eta}(t)=\{\varphi_\eta(t,x): x\in\Gamma\}.
	\end{equation} 
	
	It is a well known result on the tubular neighborhood, see e.g. \cite[Section 10]{Lee} that there are numbers $a_{\partial\Omega},b_{\partial\Omega}$ such that for $\eta\in(a_{\partial\Omega},b_{\partial\Omega})$ $\varphi_{\eta}(t,\cdot)$ is a bijective parametrization of the surface $\Sigma_\eta(t)$.
	Further we denote by $\nu_{\eta}$ the normal-direction to the deformed middle surface $\varphi_{\eta}(\Gamma)$ at the point $\varphi_{\eta}(x)$ and is given by
	$$\nu_{\eta}(x)=\partial_{1}\varphi_{\eta}(x)\times \partial_{2}\varphi_{\eta}(x).$$
	Now let us introduce the dynamics of a mixture of two compressible fluids contained in the fluid domain $Q^T_\eta=\bigcup_{t\in I}\Omega_\eta(t)\times\{t\}$, where $I=(0,T)$, and its interaction with the shell evolving at the fluid-solid interface. The evolution of the mixture interacting with a Koiter shell is described by the following set of equations. 
	\begin{equation}\label{fluidsysreform}
		\begin{alignedat}{2}
			\partial_{t}(\rho)+\dvr(\rho u)=&0 &&\mbox{ in } Q^T_\eta,\\
			\partial_{t}(Z)+\dvr(Zu)=&0 &&\mbox{ in } Q^T_\eta,\\
			\partial_{t}\left((\rho+Z)u\right)+\dvr\left((\rho+Z)u\otimes u\right)=&\dvr \mathbb{S}(\mathbb D u)-\nabla P(\rho,Z) &&\mbox{ in }Q^T_\eta,\\
			\partial^{2}_{t}\eta-\zeta\partial_{t}\Delta \eta+K'(\eta)=&F\cdot\nu&&\mbox{ on } I\times\Gamma,\\
			u(t,\varphi_\eta(t,x))=&\partial_{t}\eta(x,t)\nu(\varphi(x))&&\mbox{ on }I\times\Gamma,\\
			\rho(\cdot,0)=\rho_{0}(x),\ Z(x,0)=&Z_{0}(x)&& \mbox{ in } \Omega_{\eta_0},\\
			(\alpha+Z)u(\cdot,0)=&M_{0}(x)&& \mbox{ in } \Omega_{\eta_0},\\
			(\eta,\partial_{t}\eta)(\cdot,0)=&(\eta_{0},\eta_{1})&&\mbox{ on }\Gamma,
		\end{alignedat}
	\end{equation}
	where $u$ is the average fluid velocity and $\rho$ and $Z$ are respectively the density of the first and the second species in the mixture. The Lam\'{e} coefficients $\mu$ and $\lambda$ satisfy physically reasonable conditions
	\begin{equation}\label{CoefAss}
		\mu,\lambda>0.
	\end{equation}
	Concerning the structural dynamics \eqref{fluidsysreform}$_{4},$ $K(\eta)$ represents the Koiter energy. Since we need a few more notations to introduce the Koiter energy and the elasticity operator, we present the details only in Section \ref{Sec:KE}. More precisely we refer the readers to \eqref{Koiterenergy} and \eqref{elasticityoperator} for details.\\
	The shell located at the fluid boundary moves due to the force exerted by the fluid mixture and hence $F$ appearing in the R.H.S of \eqref{fluidsysreform}$_{4}$ bears the meaning
	\begin{equation}\label{defF}
		\begin{array}{l}
			F=(-T_{f}\nu_{\eta})\circ\varphi_{\eta}|\nabla\varphi_{\eta(t)}|,
		\end{array}
	\end{equation}
	where the stress tensor of the fluid is denoted by $T_{f}$ and is defined as:
	\begin{equation*}
		\begin{split} T_{f}=T_{f}(\mathbb D u,P(\rho,Z))
			=\mathbb{S}(\mathbb D u)-P(\rho,Z)\mathbb{I}_{3},
		\end{split}
	\end{equation*}
	where 
	\begin{equation}\label{SDef}
		\mathbb{S}(\mathbb D u)=2\mu\left(\mathbb Du-\frac{1}{3}\dvr u\mathbb{I}_{3}\right)+\lambda\dvr u\mathbb{I}_{3}
	\end{equation}
	denotes the stress tensor due to the fluid. Here $\mathbb D u=\frac{1}{2}\left(\nabla u+\nabla^\top u\right)$ is the symmetric part of $\nabla u$, $\nabla^\top$ stands for the gradient transpose and $\mathbb{I}_{3}$ is the $3\times 3$ identity matrix. The term $-\zeta\partial_{t}\Delta\eta$ in \eqref{fluidsysreform}$_{4}$ models the damping of the beam due to friction. In our case the damping parameter $\zeta$ can be both zero or positive.\\
	To make the presentation concise, we start with \eqref{fluidsysreform}, a reformulated version of a more physical model (the derivation of such a model without the structure can be found in \cite{BDGGH}).\footnote{Such model is a version of one velocity  Baer-Nunziato system with dissipation. }\\
	Next, we present some hypotheses on the initial conditions and the structure of the pressure $P(\rho,Z).$
	
	\subsection{Hypotheses}\label{hypothpressure}
	Here we make a list of hypothesis under which we prove the main result. The first two hypothesis are related to initial data.
	\begin{description}
		\item[H1]
		Denoting 
		\begin{equation}\label{COaDef}
			\mathcal O_{\underline a}=\{(\rho,Z)\in\eR^2|\ \rho\in[0,\infty),\ \underline a\rho<Z<\overline a\rho\}
		\end{equation}
		for some $0< \underline a<\overline a<\infty$ we assume
		\begin{equation} \label{eq2.1}
			(\rho_0, Z_0)\in \overline{\mathcal O_{\underline a}}=\{(\rho,Z)\in\eR^2|\ \rho\in[0,\infty),\ \underline a\rho\le Z\le\overline a\rho\}
		\end{equation}
		The following convention for fractions of the form $\frac{Z}{\rho}$ provided $(\rho,Z)\in\overline{\mathcal O_{\underline a}}$  is used systematically: 
		\begin{equation}\label{conv}
			\frac{Z}{\rho}=
			\begin{cases}
				\frac{Z}{\rho}&\text{ if }\rho>0,\\
				0&\text{ if }\rho=0.
			\end{cases}
		\end{equation}
		\item[H2]
		\begin{equation} \label{eq2.6}
			\begin{split}
				&\rho_0,Z_0\geq 0,\ \rho_0,Z_0\not\equiv 0\text{ a.e. in }\Omega_{\eta_0},\ \rho_0 \in L^\gamma(\Omega_{\eta_0}),\ \; Z_0 \in L^\beta(\Omega_{\eta_0}), \\
				& M_0=(\rho_{0}+Z_{0})u_{0}\in L^1(\Omega_{\eta_0}), (\rho_0+Z_0)|u_0|^2\in L^1(\Omega_{\eta_0}),\ \eta_{0}\in W^{2,2}(\Gamma),\ \eta_{1}\in L^{2}(\Gamma).
			\end{split}
		\end{equation}
		\item[H3]
		The pressure function $P$ is supposed  to belong to the class $C(\overline{\mathcal O_{\underline a}})\cap C^1(\mathcal O_{\underline a})$
		and to be such that
		\begin{equation}\label{PowerEst}
			\text{ for all }\rho\in(0,1)\ \sup_{s\in[\underline a,\overline a]}|P(\rho, \rho s)|\le \overline C\rho^\alpha
		\end{equation}
		with some $\alpha >0$,
		\begin{equation}\label{gambetabnd}
			\text{ for all }(\rho,Z)\in\overline{O_{\underline a}}\hspace{1em} \underline{C} (\rho^\gamma+Z^\beta-1)\le P(\rho,Z)\le \overline{C}(\rho^\gamma+Z^\beta+1)
		\end{equation}
		with $\max\{\gamma,\beta\}\ge 2$, $\min\{\gamma,\beta\}>0$ and positive constants $\underline{C}$, $\overline{C}$ and
		\begin{equation}\label{PressDerEst*}
			\text{ for all }(\rho,Z)\in\overline{O_{\underline a}}\ |\partial_Z P(\rho, Z)|\le C(\rho^{-\underline\kappa}+\rho^{\overline\kappa-1}),
		\end{equation}
		with some $0\le\underline{\kappa}<1$ and with some $0<\overline{\kappa}<\max\{\gamma+\gamma_{BOG},\,\beta+\beta_{BOG}\}$ where $\gamma_{BOG}=\min\{\frac{2}{3}\gamma-1,\frac{\gamma}{2}\}$ and $\beta_{BOG}=\min\{\frac{2}{3}\beta-1,\frac{\beta}{2}\}$ are the improvement of the integrability of the densities due to the estimates involving Bogovskii operator. 
		
		\item[H4]
		It is assumed that
		\begin{equation} \label{eq2.4}
			P(\rho,\rho s)=
			{\mathcal P}(\rho,s) - {\mathcal R} (\rho,s),
		\end{equation}
		where $[0,\infty)\ni \rho\mapsto {\mathcal P}(\rho,s)$ is non decreasing  for any $s\in [\underline a,\overline a]$, and $\rho\mapsto {\mathcal R}(\rho,s)$ is for any $s \in [\underline a,\overline a]$ a non-negative $C^2$-function  in $[0,\infty)$ uniformly bounded with respect to $s \in [\underline a, \overline a]$ with compact support uniform with respect to $s \in [\underline a, \overline a]$, i.e., for some $\overline R>0$
		\begin{equation}\label{bndbR}
			\bigcup_{s\in[\underline a,\overline a]}\supp \mathcal R(\cdot,s)\subset[0,\overline R],\ 
			\sup_{s\in[\underline a,\overline a]}\|\mathcal R(\cdot,s)\|_{C^2([0,\overline R])}<\infty.
		\end{equation}
		The constants $\underline a$ and $\overline a$ come from \eqref{COaDef}.
		Moreover, if $\max\{\gamma,\beta\}=2$ it is assumed that
		\begin{equation}\label{?!}
			{\mathcal P}(\rho,s)=f(s)\rho^{\max\{\gamma,\beta\}}+\pi(\rho,s),
		\end{equation}
		where $[0,\infty)\ni \rho\mapsto \pi(\rho,s)$  is  non decreasing  for any $s\in [\underline a,\overline a]$ and $f\in L^\infty(\underline a,\overline a)$, 
		${\rm ess\, inf}_{s\in(\underline a,\overline a)} f(s)\ge \underline f>0$.
		\item[H5]
		It is assumed that the function $\rho\mapsto P(\rho,Z)$, $Z>0$   
		and the function $Z\mapsto \partial_Z P(\rho,Z)$, $\rho>0$ are Lipschitz on $(Z/\overline a, Z/\underline a)\cap(r,\infty)$, $(\underline a\rho,\overline a\rho)\cap(r,\infty)$ respectively, for all $r>0$ with Lipschitz constants
		\begin{equation}\label{eq2.3a-}
			L_P\le C(r)(1+\rho^A),\ L_P\le C(r)(1+Z^A)\text{ respectively,} 
		\end{equation}
		with some non negative number $A$. Number $C(r)$ may diverge to $+\infty$ as $r\to 0_+$.
	\end{description}
	In order to express the influence of pressure $P$ in the energy identity for system \eqref{fluidsysreform},
	we employ the Helmholtz free energy function $H_P$. It is obtained as a solution of the following first order partial differential equation in $\mathcal O_{\underline a}$
	\begin{equation}\label{FOPDE}
		P(\rho,Z)=\rho\partial_\rho H_P(\rho,Z)+Z\partial_ZH_P(\rho,Z)-H_P(\rho,Z).
	\end{equation}
	One of admissible explicit solutions to \eqref{FOPDE}, found by the method of characteristic, is of the form
	\begin{equation}\label{HelmFDef}
		H_P(\rho,Z)=\rho\int_1^\rho\frac{P(s,s\frac{Z}{\rho})}{s^2}\ds\text{ for }\rho>0, H_P(0,0)=0.    
	\end{equation}
	Next we remark on the role of the hypotheses listed in \textbf{H1}-\textbf{H5} and the exact locations in the present article where some of them are used.
	\begin{remark}
		$(1)$ Hypothesis \textbf{H1} is about the comparability of the initial data for the densities $\rho$ and $Z.$ This allows us to prove a comparability of the densities throughout the entire space-time cylinder of existence. Our strategy relies on a time discretization and further solving  structural and fluid sub problems separately. For the fluid part we use the result from \cite{NovoPoko}. In the paper \cite{NovoPoko}, the comparability of densities (from the assumption done initially) is first proved at a viscous approximation layer of the continuity equations by some maximal principle. In the present article we can use this comparability at each discrete layer since the relation is preserved under weak convergences.\\[2.mm]
		$(2)$ Hypothesis \eqref{PowerEst} tells us that $P(0,0)=0$ and further renders the continuity of the Helmotz functional $H_{\rho}$ (introduced in \eqref{HelmFDef}). This is further used while showing the energy inequality, $cf.$ the discussion after \eqref{AuxConv}.\\[2.mm]
		$(3)$ The assumption \eqref{gambetabnd} allows to obtain both estimates on the densities from the one of the pressure (available via energy inequality) and vice-versa. In particular to find an application of the upper bound from \eqref{gambetabnd} to prove the equi-integrability of the pressure we would like to refer the readers to the discussion between \eqref{convdeltapen} and \eqref{afterusingdelavalle}.\\[2.mm]
		$(4)$ The inequality \eqref{PressDerEst*} asserts that the pressure is Lipschitz in its second component and further provides an estimate of the Lipschitz constant depending on the first argument. Such an estimate plays a crucial role for a compactness argument rendering that one of the densities in the expression of pressure can be fixed and thereby providing an access to the Lions-Feireisl theory. We refer the readers to the arguments leading to \eqref{WLPrId} for details.\\[2.mm]
		$(5)$ The decomposition of the pressure \eqref{eq2.4}-\eqref{bndbR} is used to identify the limit of the pressure. For details the readers can have a look into the arguments leading respectively to \eqref{decomposepidelta}, \eqref{PressDec} and \eqref{calafterdiffLKr}.\\[2.mm]
		$(6)$ The particular structural assumption \eqref{?!} for the pressure in the critical case $\max\{\gamma,\beta\}=2$ is needed for controlling the amplitude of density oscillation , more specifically \eqref{osdensity}. Since the proof of \eqref{osdensity} follows the arguments used in \cite[Proposition 14]{NovoPoko}, we do not provide the details for the same in the present article.\\[2.mm]
		$(7)$ The assumption \eqref{eq2.3a-} is directly not used in the present article. In fact this technical assumption is used at a Galerkin level in \cite[Section 4.1]{NovoPoko}. Since our strategy is based on an existence result for a fluid sub-problem (see Theorem \ref{existencediscretefluid}) in a fixed domain proved in \cite{NovoPoko}, the implicit use of \eqref{eq2.3a-} is hidden in Theorem \ref{existencediscretefluid}. 
	\end{remark}
	\begin{remark}
		One notices the appearance of $\gamma_{BOG}$ and $\beta_{BOG}$ in the assumption \eqref{PressDerEst*}. They are precisely the improvement of the integrability exponents of the density due to the argument using Bogovskii operator. Now we do not have this extra integrability of the densities up to the interface $(0,T)\times\Sigma_{\eta},$ since it is not uniformly Lipschitz and hence a Bogovskii type argument can not be used up to the boundary. It is worth noticing that the only places we use \eqref{PressDerEst*} is while freezing one of the densities in the expression of the pressure by using a almost compactness argument (we refer to the proof of \eqref{estdiffpress} and \eqref{estdiffpressdel}). These arguments only uses \eqref{PressDerEst*} applied in parabolic cylinders away from the interface $(0,T)\times\Sigma_{\eta}$ where we still have the Bogovskii type improvement and hence the upper bound of $\overline{\kappa}$ in \eqref{PressDerEst*} is justified.  
	\end{remark}
	Let us give an example of a physical pressure law which solve the Hypotheses \textbf{H1-H5}. We take this example from \cite{NovoPoko} with a minor change in the range of adiabatic exponent (this adaptation is required since in our case the critical adiabatic exponent is $2$ whereas for \cite{NovoPoko} it is $\frac{9}{5}$).\footnote{Let us mention that for multi-component case with nonhomogeneous boundary data, we need $\gamma >2$.}\\
	\begin{equation}\label{examplepress}
		\begin{array}{l}
			\displaystyle P(\rho,Z)=\rho^{\gamma}+Z^{\beta}+\sum^{M}_{i=1}F_{i}(\rho,Z),
		\end{array}
	\end{equation}
	where $F_{i}(\rho,Z)=C_{i}\rho^{r_{i}}Z^{s_{i}},$ $0\leq r_{i}<\gamma,$ $0\leq s_{i}<\beta$ and $r_{i}+s_{i}<\max\{\gamma+\beta\}.$ If $\gamma>2,$ we allow $C_{i}$ to be negative and hence some non-monotone choice of the pressure is possible. For $\gamma= 2,$ we assume $C_{i}\geq 0.$
	\subsection{Definition of weak solution and main result}
	Let us define the notion of bounded energy weak solution to the system \eqref{fluidsysreform}. 
	\begin{mydef}\label{WSDef}
		The quadruple $(\rho,Z,u,\eta)$ is a bounded energy weak solution to the problem \eqref{fluidsysreform} if
		\begin{align}
			&\rho, Z\ge 0\text{ a.e. in }Q^\eta_T,\nonumber\\
			&\rho\in C_w([0,T];L^{\max\{\gamma,\beta\}}(\Omega_{\eta}(t))),\nonumber\\
			&Z\in C_w([0,T];L^{\max\{\gamma,\beta\}}(\Omega_{\eta}(t))),\nonumber\\
			&u\in L^{2}(0,T;W^{1,q}(\Omega_{\eta}(t))),\,\,q<2,\label{listregsol}\\
			&(\rho+Z)u\in C_w([0,T]; L^\frac{2\max\{\gamma,\beta\}}{\max\{\gamma,\beta\}+1}(\Omega_\eta(t))),\nonumber\\
			&(\rho+Z)|u|^{2}\in L^{\infty}(0,T;L^{1}(\Omega_{\eta}(t))),\nonumber\\
			&\eta\in L^{\infty}(0,T;W^{2,2}(\Gamma))\cap L^2(0,T;W^{2+\sigma,2}(\Gamma))\mbox{ for }\sigma>0,\\
			&\tder \eta\in C_w([0,T];L^2(\Gamma))\cap L^2(0,T;W^{\sigma,2}(\Gamma))\mbox{ for }\sigma>0,\nonumber\\
			&P(\rho,Z)\in L^{1}(Q^\eta_T)\nonumber
		\end{align}
		and the following hold.
		\begin{enumerate}[leftmargin=5ex,label=(\roman*), topsep=-1.5ex, itemsep=1ex]
			\item The coupling of $u$ and $\tder\eta$ reads $\Tr_{\Sigma_\eta} u=\tder\eta\nu$, where the operator $\Tr_{\Sigma_\eta}$ is defined in Lemma \ref{Lem:TrOp}.
			\item The momentum equation is satisfied in the sense
			\begin{equation}\label{momentum}
				\begin{split}
					&\int_0^t\int_{\Omega_\eta(s)}(\rho+Z)u\cdot\tder\phi+\int_0^t\int_{\Omega_\eta(s)}\left((\rho+Z)u\otimes u\right)\cdot \nabla\phi-\int_0^t\int_{\Omega_\eta(s)} \mathbb{S}(\mathbb D u)\cdot\nabla\phi\\
					&+\int_0^t\int_{\Omega_\eta(s)}P(\rho,Z)\dvr\phi+\int_0^t\int_{\Gamma}\tder\eta \tder b-\int^{t}_{0}\langle K'(\eta),b\rangle+\zeta\int_{(0,t)\times\Gamma}\partial_{t}\nabla\eta\nabla b\\
					&=\int_{\Omega_\eta(t)}(\rho+Z)u(t,\cdot)\phi(t,\cdot)-\int_{\Omega_{\eta_0}}M_0\phi(0,\cdot)+\int_{\Gamma}\tder\eta(t,\cdot) b(t,\cdot)-\int_{\Gamma}\eta_1 b(0,\cdot)
				\end{split}
			\end{equation}
			for all $t\in[0,T]$, $(\phi,b)\in C^{\infty}([0,T]\times\RR^{3})\times (L^{2}(0,T;W^{2+\sigma,2}(\Gamma))\cap W^{1,\infty}(0,T;L^{2}(\Gamma)),$ for some $\sigma>0$ with $tr_{\Sigma_\eta}\phi=b\nu$. 
			
			\item The continuity equations are solved in the sense
			\begin{equation}\label{contrho}
				\begin{split}
					\int_0^t\int_{\Omega_{\eta}(s)}\rho(\tder\psi+u\cdot\nabla\psi)=&\int_{\Omega_{\eta}(t)}\rho(t,\cdot)\psi(t,\cdot)-\int_{\Omega_{\eta_0}}\rho_0\psi(0,\cdot),\\
					\int_0^t\int_{\Omega_{\eta}(s)}Z(\tder\psi+u\cdot\nabla\psi)=&\int_{\Omega_{\eta}(t)}Z(t,\cdot)\psi(t,\cdot)-\int_{\Omega_{\eta_0}}Z_0\psi(0,\cdot)
				\end{split}
			\end{equation}
			for all $t\in[0,T]$, $\psi\in C^{\infty}([0,T]\times\eR^3).$
			\item The energy inequality
			\begin{equation}\label{energybalance}
				\begin{split}
					&\int_{\Omega_{\eta}(t)}\bigg(\frac{1}{2}(\rho+Z)|u|^{2}+H_{P}(\rho,Z)\bigg)(t,\cdot)+\int_{I}\int_{\Omega_{\eta}(t)}\mathbb{S}(\mathbb D u)\cdot \nabla u+\bigg(\int_{\Gamma}\frac{1}{2}|\partial_{t}\eta|^{2}+\zeta|\tder\nabla\eta|^2+K(\eta)(t,\cdot)\bigg)\\
					&\le\int_{\Omega_{\eta_0}}\bigg(\frac{|M_0|^{2}}{2(\rho_0+Z_0)}+H_{P}(\rho_0,Z_0)\bigg)+\bigg(\frac{1}{2}\int_{\Gamma}|\eta_{1}|^{2}+K(\eta_{0})\bigg)
				\end{split}
			\end{equation}
			holds for a.a. $t\in I$.
		\end{enumerate}
	\end{mydef}
	\begin{remark}
		Notice that the test functions $b$ for the structure in \eqref{momentum} belong to the space $L^{2}(0,T;W^{2+\sigma,2}(\Gamma))$ for some $\sigma>0,$ which by continuous embedding infers $b\in L^{2}(0,T;W^{2,p}(\Gamma))$ for some $p>2.$ This is in coherence with the test functions used in \eqref{elasticityoperator} while introducing the elasticity operator $K'(\eta).$
	\end{remark}
	Having all necessary ingredients introduced we state the main result of the article.
	\begin{thm}\label{Thm:main}
		Assume that $\Omega\subset\eR^3$ is a given bounded domain with the parametrization of its boundary by a $C^4$ injective mapping $\varphi$ via $\partial\Omega=\varphi(\Gamma)$ for a two--dimensional torus $\Gamma$. Suppose that hypotheses \textbf{H1}--\textbf{H5} hold and $\eta_0$ satisfies $\eta_0\in(a_{\partial\Omega}, b_{\partial\Omega})$, $\bar\gamma(\eta_0)>0$ with $\bar \gamma$ defined in \eqref{ovgamma} and one of the following holds. \\
		$\textit{Case I:}$ Let the structural dissipation parameter $\zeta=0$ and suppose $\max\{\gamma,\beta\}> 2$, $0<\min\{\gamma,\beta\}$. 
		$\textit{Case II:}$ Let the structural dissipation parameter $\zeta>0$ and suppose $\max\{\gamma,\beta\}\geq 2,\,\,0<\min\{\gamma,\beta\}$. 
		
		Then there is $T_F\in(0,\infty]$ and a weak solution to the problem \eqref{fluidsysreform} along with the non-linear Koiter energy (the details of the structure of the Koiter energy is presented in \eqref{Koiterenergy}) on the interval $(0,T)$ for any $T<T_F$ in the sense of Definition \ref{WSDef}. Furthermore, for 'Case II' we obtain that $u$ and $\eta$ enjoys the following improved regularity
		\begin{equation}\label{betregu}
			\displaystyle u\in L^{2}(0,T;W^{1,2}(\Omega_{\eta}(t))),\,\,\eta\in W^{1,2}(0,T;W^{1,2}(\Gamma)).
		\end{equation}
		Moreover, in both of the cases ($i.e.$ 'Case I' and 'Case II') above the initial data are attained in the sense
		\begin{equation}\label{InDatAtt}
			\begin{split}
				\lim_{t\to 0_+}\int_{\Omega_\eta(t)}\rho(t)\psi=\int_{\Omega_{\eta_0}}\rho_0\psi,\ \lim_{t\to 0_+}\int_{\Omega_{\eta}(t)}Z(t)\psi=\int_{\Omega_{\eta_0}}Z_0\psi=0,\\
				\lim_{t\to 0_+}\int_{\Omega_\eta(t)}(\rho+Z)u(t)\phi=\int_{\Omega_{\eta_0}}M_0\phi,\ \lim_{t\to 0_+}\int_\Gamma\tder\eta(t)g=\int_\Gamma\eta_1g
			\end{split}
		\end{equation}
		for any $g\in C^\infty(\Gamma)$, $\psi\in C^\infty_c(\eR^3)$, $\phi\in C^\infty_c(\eR^3)$ such that the support of $\psi\circ\tilde\varphi$ and $\phi\circ\tilde\varphi$ as well is compact in $[0,T]\times\Omega$.
		Finally, $T_F$ is finite only if 
		\begin{equation}\label{degenfstkind}
			\text{ either }\lim_{s\to T_F}\eta(s,y)\searrow a_{\partial\Omega}\text{ or }\lim_{s\to T_F}\eta(s,y)\nearrow b_{\partial\Omega}
		\end{equation}
		for some $y\in\Gamma$ or the Koiter energy degenerates, i.e., 
		\begin{equation}\label{degensndkind}
			\lim_{s\to T_F}\bar\gamma(\eta(s,y))=0 
		\end{equation}
		for some $y\in\Gamma.$
	\end{thm}
	
	\begin{remark}
		We point out that \eqref{degenfstkind} excludes the possibility of the self-intersection for the structure. Indeed, for $t<T_F$ we have $\eta(t,\cdot)\in (a_{\partial\Omega},b_{\partial\Omega})$. Hence the flow function $\tilde\varphi_{\eta}$ governing the deformation of $\Omega$ is invertible, see Section~\ref{Sec:GEE}.
	\end{remark}

	\subsection{Ideas, strategy and some further comments} 
	\begin{itemize}
		\item To prove a global existence result (up to a self-intersection of the structure) we extend our problem to larger domain $B$ with the regular boundary such that the moving interface lies in the interior of $B$. At the same time we approximate the viscosity coefficients, initial data and the pressure in a suitable way keeping in mind that we need to recover the weak formulations in the physical domain from the ones in the extended set-up by means of suitable limit passage. For the details about the extension of viscosity coefficients and data (the $\omega,\delta$ level) we refer the readers to Section \ref{extendatacoeff}. 
		
		At this moment we want to point out that we have extended the initial densities by zero outside $\Omega_{\eta_{0}}$ (or more precisely a regularized version of the same, which will be clear from the context). {\it Such an extension guarantees (cf. Lemma \ref{Lem:Fund}) that the densities stay zero outside the physical domain for the entire time horizon.} This Lemma is one of fundamental Lemmas in the proof of existence.  For the proof of Lemma \ref{Lem:Fund} one requires the $W^{1,2}((0,T)\times\Gamma)$ regularity of the interface, which can be achieved due to the structural dissipation $-\zeta\partial_{t}\Delta\eta.$ 
		
		In case the adiabatic exponents solve $\max\{\gamma,\beta\}>2$ and $\min\{\gamma,\beta\}>0$ we will get rid of this structural dissipation $-\zeta\partial_{t}\Delta\eta$ with the aid of suitable uniform estimates later in our analysis.\\[2.mm]
		\item Next in order to solve the new system in the extended set-up we introduce a discretization of the time interval $(0,T)$ in steps of size $\tau.$ Further in each time stepping of length $\tau,$ {\it we split the problem into   two decoupled systems of  equations, one concerning the structure and one for the fluid mixture}. 
		{\bf This splitting does not preserve the interface coupling between the solid and the fluid velocities.} Instead we add penalization terms both to the structural as well as fluid sub-problems which helps later (while passing $\tau\rightarrow 0$) in recovering the kinematic coupling condition on the interface. The penalization we use is of Brinkman type and is inspired from \cite{FeireislNes} and \cite{MaMuNeRoTr}. The fluid sub-problem can be solved by imitating step by step the arguments from \cite{NovoPoko} with very minor modifications. 
		
		Concerning the solid sub-problem the intricate part is to deal with {\it the non-linear, non-convex Koiter energy}. One can notice from the structure of $K'(\eta)$ (we refer to \eqref{elasticityoperator}-\eqref{exaR}) that it consists of a term which is roughly of the form $\int_{\Gamma}\nabla^{2}\eta\cdot\nabla^{2}\eta b$ ($b$ being the test function for the structure).  
		One of difficulty in this point is that the weak$^*$ convergence of the approximates of $\eta$ in the natural energy space $L^{\infty}(W^{2,2}(\Gamma))\cap W^{1,\infty}(L^{2}(\Gamma))$ {\bf is not sufficient for the limit passage in this non-linearity}. Indeed, one can think of the ingenious idea introduced in \cite{MuhaSch} (later used in \cite{Breit2}) of improving the regularity of the structural displacement in the space $L^{2}(W^{2+\sigma,2}(\Gamma))$ and thereby obtaining the strong convergence of the approximating sequence in $L^{2}(W^{2,2}(\Gamma))$ (by using the Aubin-Lions compactness argument). But such an argument can not be applied at this stage since we have lost the fluid-solid interface coupling due to splitting and penalization. To circumvent this difficulty we introduce a further regularization of the Koiter energy $K(\eta)$ by adding a term of the form $\delta^{7}\int_{\Gamma}|\nabla^{3}\eta|^{2}$ with some parameter $\delta$ (we refer to Section \ref{regshell}). This indeed provides us with the required compactness of $\eta.$ 
		
		The structural sub-problem is next solved by a further time discretization with time stepping $\Delta t<< \tau.$ For each $\Delta t$ we solve stationary problems with suitable discretization of the non-linear Koiter energy (such a discretization is inspired from \cite{MuhaSch}). Relying on the estimates uniform in $\Delta t,$ we obtain convergence of interpolants which are sufficient to pass $\Delta t\rightarrow 0$ and recover a solution of the structural sub-problem.\\[2.mm] 
		
		\item The next steps are {\bf limit passages $\tau\rightarrow 0$ and $\delta\rightarrow 0.$} We would like to point out here that to make the presentation concise, we make the limit passages for the regularizing parameter $\delta,$ the dissipation parameter $\zeta$ and the viscosity approximation parameter $\omega$ all at the same level (we refer to \eqref{delomgep0}, Section \ref{reglim0}). As it is typical for compressible fluids, the limit passage in the non-linear pressure term is quite involved. In case of a Lipschitz domain, one can use an argument involving Bogovski\u{i} operator to obtain better integrability of the pressure. But in the present scenario we have uniform (in $\delta$) apriori $L^{\infty}(W^{2,2}(\Gamma))$ regularity of the structural displacement which only allows to obtain that the $\delta-$ approximates of the interface are uniformly $C^{0,\alpha}(\Gamma)$ ($\alpha<1$) regular. So, {\it a Bogovski\u{i} type argument fails}. The way out is to exclude possible concentration of the pressure near the interface and to prove the equi-integrability of the same. The equi-integrability of the pressure furnishes us with a $L^{1}-$ weak sub-sequential limit of the pressure. This is done in the spirit of \cite{Breit}.\\
		\item The next step is to {\bf identify the limit of the pressure}. In order to deal with a compressible bi-fluid model in a time independent smooth domain, the authors in \cite{Vasseur} developed an ingenious idea which amounts in freezing one of the variables in the pressure law and later improved by the authors of \cite{NovoPoko} to incorporate more intricate non-monotone pressure functions. The idea of \cite{NovoPoko} is to prove an almost compactness of the quantity $\frac{Z}{\rho}$ (where $Z$ and $\rho$ are the partial densities of corresponding two fluids). Once such an almost compactness is established, the pressure can be written as a function of a single density $\rho$ and the compactness of $\rho$ can be furnished following the arguments of Feireisl-Lions. We use a similar strategy as that of \cite{NovoPoko} but adapted to the case of a time dependent H\"{o}lder domain. {\it The almost compactness result (in a time varying domain)} is stated in form of Lemma \ref{Lem:AlmComp} and this can be of independent interest to the readers.\\       
		Both for proving the almost compactness result \ref{Lem:AlmComp} and the strong convergence of density (as it is by now classical from Feireisl-Lions approach) we use the renormalized continuity equation.{\it In the present article we prove a result about the existence of a solution to the continuity equation in the renormalized sense in the context of a time varying domain (cf. Lemma \ref{Lem:Renormalization}) in a very general form}. Hence it may be found of independent interest. 
		
		The result concerns two cases;
		
		\noindent (i) $\textit{Case\,I:}$ the function $\eta$ (describing the boundary of the domain $\Omega_{\eta(t)}$) solves a hyperblolic equation (i.e. the dissipation parameter $\zeta=0$), the fluid velocity $u\in L^{2}(0,T;W^{1,q}(\Omega_{\eta(t)})$ for $q<2$ and the fluid density possesses $L^{\infty}(0,T;L^{\widetilde{\gamma}}(\Omega_{\eta(t)}))$ integrability with $\widetilde{\gamma}>2$ 
		
		\noindent (ii) $\textit{Case II:}$ the function $\eta$ solves a parabolic equation (i.e. the dissipation parameter $\zeta>0$ and consequently $\eta\in W^{1,2}(0,T;W^{1,2}(\Gamma)$), the fluid velocity $u\in L^{2}(0,T;W^{1,2}(\Omega_{\eta(t)})$ and the fluid density possesses $L^{\infty}(0,T;L^{2}(\Omega_{\eta(t)}))$ integrability.\\
		Since our fluid boundary is only H\"{o}lder continuous uniformly in time (because $\eta\in L^{\infty}(0,T;W^{2,2}(\Gamma))$) in case there is no structural dissipation, we can only obtain the velocity field $u\in L^{2}(W^{1,q}(\Omega_{\eta})),$ $q<2,$ (this is an application of Lemma \ref{Lem:Korn}) and hence we need the assumption $\widetilde{\gamma}>2$ (where $\widetilde{\gamma}$ is the adiabatic exponent of one of the fluids) in order to use the Friedrichs commutator lemma to furnish a proof of the existence of renormalized continuity equation. In case the structure is of dissipative nature ($i.e.$ $\zeta>0$) we can recover $u\in L^{2}(W^{1,2}(\Omega_{\eta}))$ by a suitable lifting argument and prove that the densities satisfy the renormalized continuity equation even when $\widetilde{\gamma}=2$. 
	\end{itemize}  
	\subsection{Bibliographical remarks} 
	
	In this section we will quote some articles on the theory of existence of compressible Navier-Stokes equations and further comment on works devoted to fluid-structure interaction problems.\\
	\begin{itemize}
		\item {\textit{(i) Mono-fluid compressible Navier-Stokes equations:}} The global existence of strong solutions for a small perturbation of a stable constant state was established in the celebrated work \cite{matnis}. In the article \cite{vallizak} the authors established the local in time existence of strong solutions in the presence of inflow and outflow of the fluid through the boundary. In the same article they also present the proof of global in time existence for small data in the absence of the inflow. P.-L. Lions proved (in \cite{lions}) the global existence of renormalized weak solution with bounded energy for an isentropic fluid (i.e $p(\rho)=\rho^{\gamma}$) with the adiabatic constant $\gamma>3d/(d+2),$ where $d$ is the space dimension. E. Feireisl $\mathit{et\, al.}$ generalized the approach to cover the range $\gamma>3/2$ in dimension $3$ and $\gamma>1,$ in dimension $2$ in \cite{FeireislPet}. Due to the possible concentration of the convection term the global existence theory for the case $1\le\gamma\le \frac{3}{2}$ remains open.\footnote{Let us mention recent result by Abbatiello  et al. \cite{A}, where such case was studied for the so-called dissipative solutions.} Let us mention the celebrated recent work \cite{BreshJabin} where the authors introduce a completely new method to obtain compactness on the density. The well-posedness issues of the compressible Navier-Stokes equations for critical regularity of data can be found in \cite{danchin}. For further references and a very detailed development of the mathematical theory of compressible flow we refer the reader into the book \cite{NovStr04}.\\[2.mm]
		\item {\textit{(ii) References on compressible multi-fluid models in time-independent domains:}} In the past few years the study of compressible bi-fluid models has drawn an immense interest. In the articles \cite{Evkar2, yaoZhu2} the authors deal with one dimensional bi-fluid models with a singular pressure law.\\
		The authors of \cite{Maltese} consider a Navier–Stokes system with variable entropy which share some similarities with a multi-component fluid model since the pressure law they consider is of the form $P(\rho,s)=\rho^{\gamma}\mathcal{T}(s),$ $\gamma\ge \frac{9}{5},$ where $s$ solves an entropy transport. By writing the pressure as $P=(\rho\mathcal{T}^{\frac{1}{\gamma}})^{\gamma}=Z^{\gamma}$ where $Z$ solves a continuity equation, the authors of \cite{Maltese} were able to apply the mono-fluid theory (in the spirit of Lions-Feireisl) to prove a global existence result for the concerned system.\\
		In the seminal work \cite{Vasseur}, the authors establish global existence of weak solutions for a bi-fluid model with a pressure law of the form $P(\rho,Z)=\rho^{\gamma}+Z^{\beta},$ where $\gamma>\frac{9}{5},$ $\beta\ge 1$ and the densities are comparable. The proof of \cite{Vasseur} relies on {\it a new compactness} of the quantity $\frac{Z}{\rho}$ which further allows for a variable reduction in the pressure law.  Improvements on the result of \cite{Vasseur} are obtained in \cite{NovoPoko} where the authors are able to incorporate more intricate non-monotone pressure functions (in the present article we consider a similar structure of the pressure, cf. Hypotheses {\bf H1}-{\bf H5}) and further extending the adiabatic exponents to $\gamma\ge \frac{9}{5},$ $0<\beta<\infty.$ In both the articles \cite{NovoPoko, Vasseur} the densities are comparable. \footnote{ Concerning more general solutions so-called dissipative solution or general boundary conditions we refer to \cite{KKNN,BNN}.} Our strategy relies on  penalization and extension of decoupled equations to a time independent fixed domain, where for the fluid part we apply the result proved in \cite{NovoPoko}. Because of this particular way of constructing solutions, our proof depends on the domination/ comparison of densities. Related to this discussion we quote here a very recent article \cite{Wen}, where the author considers a bi-fluid system in a time-independent domain of class $C^{2+s},$ $s>0$ with a pressure law of the form $P(\rho,Z)=\rho^{\gamma}+Z^{\beta},$ $\gamma,\beta>\frac{9}{5},$ and without any domination/ comparison of the densities involved. The result of \cite{Wen} extends the one proved in \cite{Vasseur} by allowing transition to each single phase flow, meaning that one of the phases can vanish in a point while the other can persist. In yet another recent article \cite{BrZat}, the authors prove the
		global existence theory of weak solutions for a two-fluid Stokes equations on the d-dimensional torus for $d=2,3.$ The proof of \cite{BrZat} relies on the Bresch-Jabin’s new compactness tools for compressible Navier-Stokes equations.\\[2.mm]
		
		\item {\textit{(iii) Fluid-structure interaction problems:}} 
		From the mathematical point of view the incompressible fluid-structure interaction problems are well studied in the literature. For the well posedness and regularity results of incompressible fluid-structure interaction (FSI) models with the structure immersed inside the fluid one can consult the articles \cite{Avalos1,DEES2,cout1,cout2,Galdi2,G2} and for incompressible fluid structure interaction problems with elastic structure at the fluid boundary we refer to \cite{AbelsLiu1, AbelsLiu2,esteban,grand,LeRu14,muhacanic2, Raymond1, veiga}. \\
		Despite of the growing literature on incompressible fluid structure interaction problems the number of articles addressing the compressible fluid structure interaction problems is relatively limited and the literature has been rather recently developed. The strong coupling between the parabolic and hyperbolic dynamics is one of the intricacies in dealing with the compressible Navier-Stokes equations and this results in the regularity incompatibilities between the fluid and the solid structure. However in the past few years there have been works exploring the fluid structure interaction problems comprising the compressible Navier-Stokes equations. For instance we refer to the articles \cite{F4,bougue2, Hieber,KrNePi2,NRRS} (rigid body immersed inside the fluid domain) and \cite{kukavica} (elastic structure inside the fluid). Further we quote the articles \cite{MaityTakahahi, SM2} (strong solutions with damped elastic structure), \cite{Avalos2} (semigroup well posedness with an undamped structure) and \cite{RoyMaity} (strong solution with a wave equation) for the analysis of compressible fluid structure interaction models with the structure/ wave appearing at the fluid boundary.\\
		The first existence result on global weak solutions (until a degeneracy occurs) to a system of compressible Navier-Stokes equations interacting with a hyperbolic elastic structure (appearing at the fluid boundary) appeared in \cite{Breit}. The elastic structure in \cite{Breit} is modeled by a linearized Koiter shell equation (the boundary is described as a graph). Next in \cite{Breit2}, the authors consider a Navier-Stokes-Fourier system and further improve their earlier result by considering non-linear, non-convex Koiter energy. A very interesting part of \cite{Breit2} is that the authors can show that the system under consideration is thermodynamically closed just by using weak regularity of the solution and a further improvement on the regularity of the structural displacement ($\eta\in L^{2}(W^{2+\sigma,2}),$ $\sigma>0$). Such an improvement of the structural regularity (only when the fluid and the shell velocity coincide at the interface) was first observed in \cite{MuhaSch} (for an incompressible FSI problem). The same improved regularity is also a key part of the present article. In yet another recent article \cite{Tri1}, the authors investigate the existence of weak solutions of a system coupling compressible Navier-Stokes and a linear thermoelastic plate equations. Recently the authors of \cite{MaMuNeRoTr} consider a FSI problem with a heat-conducting fluid which is in a thermal equilibrium with a linear thermoelastic structure constituting the fluid-boundary.\\
		Finally we wish to refer the readers to a couple of very interesting articles \cite{BeScKam} and \cite{BreitKMSch} where the authors develop a variational strategy to deal with a fluid-structure interaction model involving a bulk structure. The novel strategy (based on minimizing movement) designed in \cite{BeScKam, BreitKMSch} furnishes a natural way of dealing with non-linear, non-convex elastic energy of a very general form. 
		
	\end{itemize}
	
	\section{Geometry, some key lemmas and properties of non-linear Koiter energy}\label{sec:Geomtry}
	The following subsection is a summary about the description of a moving domain and a few lemmas on the Sobolev embedding, extension operators and Korn's inequality concerning domains with H\"{o}lder continuous boundary.
	\subsection{Geometry, embedding and extension}\label{Sec:GEE}
	This section contains a collection of facts related to domains with a moving boundary.
	Using the notation from the introductory part of this paper, we define the tubular neighbourhood of $\partial\Omega$ as
	\begin{equation*}
		N^b_a=\{\varphi(x)+\nu(\varphi(x))z; x\in\Gamma, z\in(a_{\partial\Omega},b_{\partial\Omega})\},
	\end{equation*}
	the projection $\pi:N^b_a\to\partial\Omega$ as a mapping that assigns to each $x$ a unique $\pi(x)\in\partial \Omega$ such that there is $z\in(a_{\partial\Omega},b_{\partial\Omega})$ and
	\begin{equation*}
		x-\pi(x)=\nu(\pi(x))z
	\end{equation*}
	and the signed distance function $d:N^b_a\to(a_{\partial\Omega},b_{\partial\Omega})$ as
	\begin{equation*}
		d:x\mapsto (x-\pi(x))\cdot\nu(\pi(x)).
	\end{equation*}
	We note that considering the function 
	\begin{equation*}
		\mathfrak d(x,\partial\Omega)=\begin{cases}-\dist(x,\partial\Omega)&\text{ if }x\in\overline\Omega\\
			\dist(x,\partial\Omega)&\text{ if }x\in\eR^3\setminus\Omega\end{cases}
	\end{equation*} 
	$d$ and $\mathfrak d$ coincide in $N^b_a$. Since it is assumed that $\varphi\in C^4(\Gamma)$, it is well known that $\pi$ is well defined and possesses the $C^3$--regularity and $d$ is $C^4$ in a neighbourhood of $\partial\Omega$ containing $N^b_a$, see \cite[Theorem 1 and Lemma 2]{Foote84}.
	Let $\eta:[0,T]\times\Gamma\to\eR$ be a given displacement function with $a_{\partial\Omega}<m\leq \eta\leq M<b_{\partial\Omega}$. We fix such a pair $\{m,M\}$ from the beginning.\\
	Then the flow function $\tilde\varphi_\eta:[0,T]\times \eR^3\to \eR^3$ is defined as
	\begin{equation}\label{FlowMDef}
		\tilde\varphi_\eta(t,x)=x+f_\Gamma(\mathfrak d(x))\eta(t,\varphi^{-1}(\pi(x)))\nu(\pi(x)).
	\end{equation}
	The cut-off function $f_\Gamma\in C^\infty_c(\eR)$, $0\leq f_\Gamma\leq 1$ is defined as 
	\begin{equation*}
		f_\Gamma(x)=(f*\omega_\alpha)(x)
	\end{equation*}
	with a standard mollifying kernel $\omega_\alpha$ possessing the support in $(-\alpha,\alpha)$ for $\alpha<\frac{1}{2}\min\{m'-m'',M''-M'\}$, where $a_{\partial_\Omega}<m''<m'<m< 0< M<M'<M''<b_{\partial\Omega}$. Furthermore the function $f\in W^{1,\infty}(\eR)$ is given by
	\begin{equation*}
		f(x)=\begin{cases}
			1&x\in(m''-m',M''-M'],\\
			1-\frac{x-m''+m'}{m'}&x\in(m'',m''-m'],\\
			1-\frac{x-M''+M'}{M'}&x\in(M''-M',M''],\\
			0&x\in(-\infty,m'']\cup(M'',\infty)
		\end{cases}
	\end{equation*}
	implying 
	\begin{equation}\label{fPrEst}
		f'_\Gamma\in\left[-\frac{1}{M'},-\frac{1}{m'}\right].
	\end{equation}
	We note that for $x\in N^b_a$ we can write
	\begin{equation*}
		\tilde\varphi_\eta(t,x)=(1-f_\Gamma(d(x))x+f_\Gamma(d(x))(\pi(x)+(d(x)+\eta(t,\varphi^{-1}(\pi(x))))\nu(\pi(x)).
	\end{equation*}
	Hence for the inverse $(\tilde\varphi_\eta)^{-1}$ we get
	\begin{equation*}
		(\tilde\varphi_\eta)^{-1}(t,z)=(1-f_\Gamma(d(z))z+f_\Gamma(d(z))(\pi(z)+(d(z)-\eta(t,\varphi^{-1}(\pi(z))))\nu(\pi(z)) \text{ for }z\in N^b_a.
	\end{equation*}
	For $(\tilde\varphi_\eta)^{-1}:[0,T]\times\eR^3\to\eR^3$ we then have 
	\begin{equation}\label{FlowMInv}
		(\tilde\varphi_\eta)^{-1}(t,z)=z-f_\Gamma(\mathfrak d(z))\eta(t,\varphi^{-1}(\pi(z)))\nu(\pi(z)).
	\end{equation}
	Obviously, the mapping $\tilde\varphi_\eta$ and its inverse inherit the regularity of $\eta$.  
	Let us summarize the assumptions on the geometry for the assertions in the rest of this section.
	
	\hypertarget{Assumption}{\textbf{Assumptions (A)}}: Let $\Omega\subset\eR^3$ be a bounded domain of class $C^4$. Let the boundary $\partial\Omega$ of $\Omega$ be parametrized as $\partial\Omega=\varphi(\Gamma)$, where $\varphi$ is a $C^4$ injective--mapping and $\Gamma\subset\eR^2$ is a torus. Let for $t\in[0,T]$ $\Omega_\eta(t)=\tilde\varphi_\eta(t,\Omega)$ with the boundary $\Sigma_\eta(t)=\tilde\varphi_\eta(t,\partial\Omega)$, where $\tilde\varphi_\eta$ is defined in \eqref{FlowMDef} for a displacement $\eta$ satisfying $a_{\partial\Omega}<\eta<b_{\partial\Omega}$. 
	
	\begin{remark}\label{Rem:SimplNot}
		Very often, if there is no threat of confusion we identify for simplicity functions defined on $\Gamma$ with functions on $\partial\Omega$.
	\end{remark}

	Under the validity of Assumption (A) for $\eta\in C([0,T]\times\Gamma)$ we define the underlying function spaces on variable domains in the following way for $p,r\in[1,\infty]$
	\begin{align*}
		L^p(0,T;L^r(\Omega_\eta(t)))=&\{v\in L^1(Q^T_\eta): v(t)\in L^r(\Omega_\eta(t))\text{ for a.e. }t\in (0,T), \|v(t)\|_{L^r(\Omega_\eta(t))}\in L^p((0,T))\},\\
		L^p(0,T;L^r(\Omega_\eta(t)))=&\{v\in L^p(0,T;L^r(\Omega_\eta(t))):\nabla v\in L^p(0,T;L^r(\Omega_\eta(t)))\}.
	\end{align*}
	Moreover, the space $C_w([0,T];L^p(\Omega_{\eta}(t)))$ consists of $v\in L^\infty(0,T;L^p(\Omega_\eta(t))$ such that the mapping $t\mapsto\int_{\Omega_{\eta}(t)}v(t)\theta$ is continuous for any $\theta\in C^\infty_c(\eR^3)$ such that the support of $\theta\circ\tilde\varphi_\eta$ is compact in $[0,T]\times\Omega$, i.e. $\theta$ is compactly supported in $\Omega_\eta(t)$ for each $t\in[0,T]$.
	
	For the purposes of this subsection we define 
	\begin{equation*}
		X=L^\infty(0,T;W^{2,2}(\Gamma))\cap W^{1,\infty}(0,T;L^2(\Gamma)).
	\end{equation*}
	Since 
	\begin{equation*}
		L^\infty(0,T;W^{2,2}(\Gamma))\cap W^{1,\infty}(0,T;L^2(\Gamma))\hookrightarrow C^{0,1-\theta}([0,T];C^{0,2\theta-1}(\Gamma))
	\end{equation*}
	for $\theta\in(\frac{1}{2},1)$, cf. \cite[(2.29)]{LeRu14}, we get $X\hookrightarrow C([0,T]\times\Gamma)$ and the above defined function spaces are meaningful for $\eta\in X$.
	
	In the case when $\tilde\varphi_\eta:\Omega\to\Omega_\eta$ induced by the mapping $\eta$ is not bi--Lipschitz, we do not have an isomorphism between corresponding Lebesgue and Sobolev spaces. The next lemma quantifies the loss in the regularity for transformations.
	\begin{lem}\label{Lem:TrCont}
		Let \hyperlink{Assumption}{Assumptions (A)} hold true with $\eta\in X$ and $p\in[1,\infty], q\in(1,\infty]$. The mapping $v\mapsto v\circ\tilde\varphi_\eta$ is continuous from $L^p(0,T;L^q(\Omega_\eta))$ to $L^p(0,T;L^r(\Omega))$ and from  $L^p(0,T;W^{1,q}(\Omega_\eta))$ to $L^p(0,T;W^{1,r}(\Omega))$ and $v\mapsto v\circ (\tilde\varphi_\eta)^{-1}$ is continuous from $L^p(0,T;L^q(\Omega))$ to $L^p(0,T;L^r(\Omega_\eta(t)))$ and from  $L^p(0,T;W^{1,q}(\Omega))$ to $L^p(0,T;W^{1,r}(\Omega_\eta))$ for any $1\leq r<q$.
		\begin{proof}
			The assertion for fixed $t\in(0,T)$, i.e. for $v(t)\mapsto (v\circ\tilde\varphi_\eta)(t)$, $v(t)\mapsto (v\circ(\tilde\varphi_\eta)^{-1})(t)$ was shown in \cite[Lemma 2.6.]{LeRu14} with the continuity constant depending also on $\|\eta(t)\|_{W^{2,2}(\Gamma)}$, which is now uniformly bounded in $t$. Hence the assertion of this lemma follows.  
		\end{proof}
	\end{lem}
	The following lemma concerns the continuity of the trace operator on a domain with the moving boundary. It is obtained similarly as the previous lemma by a combination of already proven time independent result in \cite[Corollary 2.9]{LeRu14} and the Sobolev embedding theorem.
	\begin{lem}\label{Lem:TrOp}
		Let \hyperlink{Assumption}{Assumptions (A)} hold true with $\eta\in X$ and $p\in[1,\infty]$, $q\in(1,\infty)$. Then the linear mapping $\Tr_{\Sigma_\eta}:v\mapsto v \circ \tilde\varphi_\eta|_{\partial\Omega}$ is well defined and continuous from $L^p(0,T;W^{1,q}(\Omega_\eta(t)))$ to $L^p(L^r(\partial\Omega))$ for all $r\in (1,\frac{2q}{3-q})$, respectively from $L^p(0,T;W^{1,q}(\Omega_\eta(t)))$ to $L^p(0,T;W^{1-\frac{1}{r},r}(\Sigma_\eta(t))$ for any $1\leq r<q$.
	\end{lem}
	
	The next lemma is devoted to the extension of a function defined on a moving domain to the whole space. It is obtained similarly as Lemma \ref{Lem:TrCont} by a combination of already proven time independent results \cite[Lemma 2.5 and Remark 2.6]{Breit}. 
	\begin{lem}\label{Lem:Extension}
		Let \hyperlink{Assumption}{Assumptions (A)} hold true with $\eta\in X$ ($i.e.$ $\eta\in L^{\infty}(0,T;C^{0,\kappa}(\Gamma)$ for $0<\kappa<1$) and  $p\in[1,\infty]$, $q\in(1,\infty]$. Then there is a continuous linear operator\\ $\mathcal E_\eta:L^p(0,T;W^{1,q}(\Omega_\eta(t)))\to L^p(0,T;W^{1,r}(\eR^3))$ for any $r\in[1,q)$ such that $\mathcal E_\eta|_{Q^T_\eta}$ is the identity. 
	\end{lem}
	Next, we state a variant of the Korn inequality on domains with varying boundaries. It is obtained similarly as Lemma \ref{Lem:TrCont} by the application of already proven time independent result on H\"older domains \cite[Lemma 3.8]{MaMuNeRoTr}.
	\begin{lem}\label{Lem:Korn}
		\textit{[Korn type inequality]} Let \hyperlink{Assumption}{Assumptions (A)} hold true with $\eta\in X$ ($i.e.$ $\eta\in L^{\infty}(0,T;C^{0,\kappa}(\Gamma)),$ $\kappa<1$) and $p\in[1,\infty]$. Moreover, let $M,L>0$, $\gamma\in \left(\frac{3}{2},\infty\right)$ and $q\in[1,2)$ be given. Then there exists a positive constant $C=C(q,M,L,\|\eta\|_{L^\infty(0,T;W^{2,2}(\Gamma)}))$ such that
		\begin{equation}\label{TiSpKorn}
			\|u\|^2_{L^2(0,T;W^{1,q}(\Omega_\eta(t)))}\leq C\left(\|\mathbb Du\|^2_{L^2(0,T;L^2(\Omega_\eta(t)))}+\int_{\Omega_\eta(t)}\rho|u|^2\right)
		\end{equation}
		for any pair $\rho,u$ such that the right hand side is finite and $\rho\geq 0$ a.e. in $Q^T_\eta$, $\|\rho\|_{L^\infty(0,T;L^\gamma(\Omega_\eta(t)))}\leq L$, $\int_{\Omega_\eta(t)}\rho\geq M$. 
	\end{lem}
	
	Next we state a result on the solenoidal extension operator which is taken from \cite[Prop. 3.3]{MuhaSch} (we also refer to \cite[Proposition 2.9]{Breit2}).
	\begin{prop}\label{smestdivfrex}
		Let \hyperlink{Assumption}{Assumptions (A)} hold true for a given $\eta \in X$ with $a_{\partial\Omega}<m\leqslant \eta\leqslant M< b_{\partial\Omega},$ there exists a tubular neighborhood $S_{m,M}$ of $\partial\Omega$ such that
		\begin{equation}\label{tbnbd}
			\begin{array}{l}
				\{\varphi(x)+z\nu(\varphi(x))\suchthat m\le z\le M\}\Subset S_{m,M}
			\end{array}
		\end{equation}
		and there are linear operators
		$$\mathcal{K}_{\eta}:L^{1}(\Gamma)\rightarrow \mathbb{R},\,\,\mathcal{F}^{\div}_{\eta}:\{\xi\in L^{1}(0,T;W^{1,1}(\Gamma))\suchthat \mathcal{K}_{\eta}(\xi)=0\}\rightarrow L^{1}(0,T;W^{1,1}_{\div}(B)),$$ 
		such that the couple $(\mathcal{F}^{\div}(\xi-\mathcal{K}_{\eta}(\xi)),\xi-\mathcal{K}_{\eta}(\xi))$ solves
		\begin{equation}\nonumber
			\begin{array}{ll}
				&\displaystyle\mathcal{F}^{\div}_{\eta}(\xi-\mathcal{K}_{\eta}(\xi))\in L^{\infty}(0,T;L^{2}(\Omega_{\eta}))\cap L^{2}(0,T;W^{1,2}_{\div}(\Omega_{\eta})),\\
				& \displaystyle\xi-\mathcal{K}_{\eta}(\xi)\in L^{\infty}(0,T;W^{2,2}(\Gamma))\cap W^{1,\infty}(0,T;L^{2}(\Gamma)),\\
				&\displaystyle tr_{\Sigma_\eta}(\mathcal{F}^{\div}_{\eta}(\xi-\mathcal{K}(\xi)))=\xi-\mathcal{K}_{\eta}(\xi),\\
				&\displaystyle \mathcal{F}^{\div}_{\eta}(\xi-\mathcal{K}_{\eta}(\xi))(t,x)=0\,\,\mbox{for}\,\,(t,x)\in (0,T)\times (\Omega\setminus S_{m,M}),
			\end{array}
		\end{equation}
		where
		\begin{equation}\label{defB}
			\begin{array}{l}
				B=B_{m,M}=\Omega\cup S_{m,M}.
			\end{array}
		\end{equation}
		Provided that $\eta,\xi\in L^{\infty}(0,T;W^{2,2}(\Gamma))\cap W^{1,\infty}(0,T;L^{2}(\Gamma)),$ one has the following estimates
		\begin{equation}\label{Fdivetaest1}
			\begin{array}{ll}
				&\displaystyle\|\mathcal{F}^{\div}_{\eta}(\xi-\mathcal{K}_{\eta}(\xi))\|_{L^{q}(0,T;W^{1,p}(B))}\lesssim \|\xi\|_{L^{q}(0,T;W^{1,p}(\Gamma))}+\|\xi\nabla\eta\|_{L^{q}(0,T;L^{p}(\Gamma))},\\
				&\displaystyle \|\partial_{t}\mathcal{F}^{\div}_{\eta}(\xi-\mathcal{K}_{\eta}(\xi))\|_{L^{q}(0,T;L^{p}(B))}\lesssim\|\partial_{t}\xi\|_{L^{q}(0,T;L^{p}(\Gamma))}+\|\xi\partial_{t}\eta\|_{L^{q}(0,T;L^{p}(\Gamma))},
			\end{array}
		\end{equation}
		for any $p\in (1,\infty),$ $q\in (1,\infty].$\\
		Further in the same spirit of the proof of \cite[Proposition 3.3]{MuhaSch} and with the assumption $\eta,\xi\in L^{\infty}(0,T;W^{3,2}(\Gamma))\cap W^{1,\infty}(0,T;L^{2}(\Gamma))$ one in particular proves that
		\begin{equation}\label{higherderest1}
			\begin{array}{ll}
				&\displaystyle \|\nabla^{3}\mathcal{F}^{\div}_{\eta}(\xi-\mathcal{K}_{\eta}(\xi))\|_{L^{\infty}(0,T;L^{2}(B))}\\
				&\displaystyle\lesssim \bigg(\||\nabla\eta||\nabla^{2}\xi|\|_{L^{\infty}(L^{2}(\Gamma))}+\||\nabla^{2}\eta||\nabla\xi|\|_{L^{\infty}(L^{2}(\Gamma))}+\|\nabla^{3}\xi\|_{L^{\infty}(L^{2}(\Gamma))}+\||\nabla\eta|^{2}|\nabla\xi|\|_{L^{\infty}(L^{2}(\Gamma))}\\
				&\displaystyle \qquad+ \||\nabla\eta|^{3}|\xi|\|_{L^{\infty}(L^{2}(\Gamma))}+\||\nabla\eta||\nabla^{2}\eta||\xi|\|_{L^{\infty}(L^{2}(\Gamma))}+\||\nabla^{3}\eta||\xi|\|_{L^{\infty}(L^{2}(\Gamma))}\bigg)\\
				&\displaystyle \lesssim \|\nabla\eta\|_{L^{\infty}(L^{\infty}(\Gamma))}\|\nabla^{2}\xi\|_{L^{\infty}(W^{1,2}(\Gamma))}+\|\nabla^{2}\eta\|_{L^{\infty}(W^{1,2}(\Gamma))}\|\nabla\xi\|_{L^{\infty}(W^{2,2}(\Gamma))}+\|\nabla^{3}\xi\|_{L^{\infty}(L^{2}(\Gamma))}\\
				&\displaystyle\qquad +\|\nabla\eta\|^{2}_{L^{\infty}(L^{\infty}(\Gamma))}\|\nabla\xi\|_{L^{\infty}(L^{\infty}(\Gamma))}+\|\nabla\eta\|^{3}_{L^{\infty}(L^{\infty}(\Gamma))}\|\xi\|_{L^{\infty}(L^{\infty}(\Gamma))}\\
				&\displaystyle\qquad+\|\nabla\eta\|_{L^{\infty}(L^{\infty}(\Gamma))}\|\nabla^{2}\eta\|_{L^{\infty}(W^{1,2}(\Gamma))}\|\xi\|_{L^{\infty}(L^{\infty}(\Gamma))}+\|\nabla^{3}\eta\|_{L^{\infty}(L^{2}(\Gamma))}\|\xi\|_{L^{\infty}(L^{\infty}(\Gamma))}
			\end{array}
		\end{equation}
	\end{prop}
	\subsection{Renormalized weak solution of continuity equation in time dependent H\"{o}lder domains}
	The next lemma is one of the most important observations of the present article and it concerns the extension of the renormalized weak solutions of continuity equation considered in varying domains. 
	\begin{lem}\label{Lem:Renormalization}
		Let \hyperlink{Assumption}{Assumptions (A)} hold true with $\eta\in X$. Let  the functions $r^{(i)}\in L^\infty(0,T;L^{\gamma_i}(\Omega_\eta(t)))$, $i=1,\ldots,M$, $\mathfrak m=\min_{i=1,\ldots,M}\{\gamma_i\}\geq 2$ with the velocity $u\in L^2(0,T;W^{1,q}(B))$, where $B$ is defined in Proposition \ref{smestdivfrex}, $q\in [1,2)$ if $\mathfrak m>2$ and $q=2$ for $\mathfrak m=2$, satisfy the continuity equation in the sense
		\begin{equation}\label{ContEq}
			\int_0^T \int_{\Omega_{\eta}(s)}r^{(i)}(\tder \phi+u\cdot\nabla\phi)=0
		\end{equation}
		for any $\phi\in C^\infty_c((0,T)\times\eR^3)$. Then for any $\mathcal B:C^1([0,\infty)^M)\to\eR$, $\nabla \mathcal B\in L^\infty((0,\infty)^M)$, $\mathcal B(0)=0$ the function $\mathcal B(r)$, where $r=(r^{(1)},\ldots, r^{(M)})$ satisfies the renormalized continuity equation 
		\begin{equation}\label{TIREq}
			\int_0^t\int_{\Omega_{\eta}(s)} \mathcal B(r)(\tder\phi+u\cdot\nabla\phi)-\left(\nabla_r \mathcal B(r)r-\mathcal B(r) \right)\dvr u\phi=\int_{\Omega_{\eta}(s)}\mathcal B(r)\phi|^{s=t}_{s=0}
		\end{equation}
		for any $t\in[0,T]$ and any $\phi\in C^\infty([0,T]\times\eR^3)$. 	
		\begin{proof} 
			As the first step we extend each density $r^{(i)}$ by zero in $[0,T]\times B\setminus Q^T_\eta$.  We then deduce from \eqref{ContEq}
			\begin{equation}\label{REq}
				\int_{(0,T)\times B}r^{(i)}(\tder \phi+u\cdot\nabla\phi)=0
			\end{equation}
			for any $\phi\in C_c^\infty((0,T)\times B)$. Let us consider a standard mollifying operator $S_\varepsilon$. Regularizing equation \eqref{REq} we get for $r^{(i)}_\varepsilon= S_\varepsilon(r^{(i)})$
			\begin{equation}\label{REqPoint}
				\tder r^{(i)}_\varepsilon+\dvr(r^{(i)}_\varepsilon u)=\dvr(r^{(i)}_\varepsilon  u)-\dvr(S_\varepsilon(r^{(i)} u))=:R_\varepsilon(r^{(i)})\text{ a.e. in }(0,T)\times B.
			\end{equation}        
			The properties of mollifiers imply 
			\begin{equation}\label{MollProp}
				r^{(i)}_\varepsilon \to r^{(i)}\text{ in }L^\infty(0,T;L^{\gamma_i}_{loc}(B))\text{ and a.e. in }(0,T)\times B
			\end{equation}
			for $i\in\{1,\ldots,M\}$. By the Friedrichs lemma on commutators we conclude
			\begin{equation}\label{CommLim}
				R_\varepsilon(r^{(i)})\to 0\text{ in }L^1(0,T;L^1_{loc}(B))
			\end{equation}
			since for $\gamma_i\geq 2$ we always find $q_{\gamma_i}\in [1,2]$ such that $\gamma^{-1}_{i}+q^{-1}_{\gamma_i}=1$. Multiplying \eqref{REqPoint} on $\partial_i \mathcal B(r_\varepsilon)$ and summing the resulting identity over $i\in\{1,\ldots, M\}$ we conclude denoting $r_\varepsilon=(r^{(1)}_\varepsilon,\ldots, r^{(M)}_\varepsilon)$
			\begin{equation}\label{BEqPoint}
				\begin{split}
					&\tder \mathcal B(r_\varepsilon)+\dvr(\mathcal B(r_\varepsilon) u)+(r_\varepsilon\cdot\nabla \mathcal B(r_\varepsilon)-\mathcal B(r_\varepsilon))\dvr u\\
					&=\sum_{i=1}^MR_\varepsilon(r^{(i)})\partial_i\mathcal B(r_\varepsilon)\text{ a.e. in }(0,T)\times B.
				\end{split}
			\end{equation} 
			From now on we fix an arbirary $\phi\in C^\infty([0,T]\times\eR^3)$, $B'$ such that $\Omega_{\eta}(t)\subset B'\subset\overline{ B'}\subset B$ for a.a. $t$ with the Lipschitz boundary and $\varepsilon_0$ such that for any $\varepsilon<\varepsilon_0$ we have $\supp r_\varepsilon(t)\subset B'$ for a.a. $t$.
			Multiplying \eqref{REqPoint} on $\phi$ and integrating over $(0,t)\times B'$ yields
			\begin{equation}\label{AuxIdH}
				\begin{split}
					\int_{B'}(\mathcal B(r_\varepsilon)\phi)(s,\cdot)|_{s=0}^{s=t}=&\int_0^t\int_{B'}\mathcal B(r_\varepsilon)\tder\phi+\int_0^t\int_{B'}\mathcal B(r_\varepsilon) u\cdot\nabla\phi \\
					&-\int_0^t\int_{B'}(r_\varepsilon\cdot\nabla \mathcal B(r_\varepsilon)-\mathcal B(r_\varepsilon))\dvr u\phi+\int_0^t\int_{B'}\sum_{i=1}^MR_\varepsilon(r^{(i)})\partial_i\mathcal B(r_\varepsilon)\phi
				\end{split} 
			\end{equation} 
			Using the assumed boundedness of $\nabla \mathcal B$ and convergence \eqref{MollProp} we infer
			\begin{equation}\label{BEpsCnv}
				\mathcal B(r_\varepsilon)\to \mathcal B(r)\text{ in }L^\infty(0,T;L^{\mathfrak m}(B'))
			\end{equation} and using additionaly the Vitali convergence theorem we get
			\begin{equation}\label{RHSCnv}
				r_\varepsilon\cdot\nabla \mathcal B(r_\varepsilon)\to r\cdot\nabla \mathcal B(r)\text{ in }L^1(0,T;L^1(B')).
			\end{equation}
			Employing \eqref{CommLim}, \eqref{BEpsCnv}, \eqref{RHSCnv} and the assumption on the boundednes of $\nabla B$ in \eqref{AuxIdH} we obtain
			\begin{equation}\label{RenId}
				\begin{split}
					\int_{B'}(\mathcal B(r)\phi)(s,\cdot)|_{s=0}^{s=t}=&\int_0^t\int_{B'}\mathcal B(r)\tder\phi+\int_0^t\int_{B'}\mathcal B(r)u\cdot\nabla\phi \\
					&-\int_0^t\int_{B'}(r\cdot\nabla \mathcal B(r)-\mathcal B(r))\dvr u\phi
				\end{split}
			\end{equation}
			for almost all $t\in [0,T]$. An immediate consequence of \eqref{RenId} is that $\tder \int_{B'}\mathcal B(r)\psi\in L^1((0,T))$ for any $\psi\in C^\infty_c(B')$. This implies that after changing $\mathcal B(r)$ on a zero measure subset of $[0,T]$ we have $\mathcal B(r)\in C_w([0,T];L^{\mathfrak m}(B'))$ and \eqref{RenId} holds for all $t\in [0,T]$. Finally, taking into account that $\mathcal B(r)=0$ in $\left((0,T)\times B'\right)\setminus Q^T_\eta$ we conclude \eqref{TIREq} from \eqref{RenId}.
		\end{proof}
	\end{lem}
	Inspired from \cite{NovoPoko}, in the following subsection we prove a compactness criterion which will help later in identifying the limit of the pressure by considering it as a function of single density (which helps to adapt some arguments from the mono-fluid theory). Compared to \cite{NovoPoko}, here we prove a version of almost compactness  which is suitable to adapt for a time varying domain.  
	\subsection{Almost compactness in the context of moving domain}
	The ensuing lemma deals with almost compactness property of sequences of solutions to transport equations on varying domains.
	\begin{lem}\label{Lem:AlmComp}
		Let a sequence $\{(\eta^n,\rho^n,Z^n,u^n)\}$ be such that
		\begin{enumerate}
			\item \hyperlink{Assumption}{Assumptions (A)} hold true for each $\eta^n$,
			\item $(\rho^n, Z^n)$ is a pair of solutions to the continuity equation in $Q^T_{\eta^n}$ prolonged by zero on $(0,T)\times B\setminus Q^T_{\eta^n}$, where $B$ comes from Proposition \ref{smestdivfrex} with the corresponding velocity $u^n$,
			\item the following estimate holds
			\begin{equation}\label{BddAss}
				\begin{split}
					\sup\left(\|\eta^n\|_{L^\infty(0,T;W^{2,2}(\Gamma))}\right.&\left.+\|\tder \eta^n\|_{L^\infty(0,T;L^2(\Gamma))}+\|\rho^n\|_{L^\infty(0,T;L^\gamma(\Omega_\eta(t)))}+\|Z^n\|_{L^\infty(0,T;L^\beta(\Omega_\eta(t)))}\right.\\
					&\left.+\|u^n\|_{L^2(0,T;W^{1,q}(B))}\right)<\infty
				\end{split}
			\end{equation}
			
			\begin{equation}\label{valqgb}
				\text{ where } q\in[1,2) \text{ if } \min\{\gamma,\beta\}>2\text{ and }q=2\text{ if }\min\{\gamma,\beta\}=2.
			\end{equation}
		\end{enumerate}
		Furthermore, let
		\begin{equation}\label{InVCnv}
			\lim_{n\to\infty}\int_{\Omega_{\eta^n_0}}\frac{(b^n_0)^2}{d^n_0}=\int_{\Omega_{\eta_0}}\frac{(b_0)^2}{d_0},
		\end{equation}
		where $b^n_0=\rho^n_0$ or $b^n_0=Z^n_0$, $d^n_0=\rho^n_0+Z^n_0$ and $b_0$ be the limit of either $\{\rho^n_0\}$ in $L^\gamma(\Omega_{\eta^0})$ or $\{Z^n_0\}$ in $L^\beta(\Omega_{\eta^0})$. Then, up to a subsequence, 
		\begin{equation}\label{ConvergN}
			\begin{alignedat}{2}
				\eta^n&\to\eta&&\text{ in }C([0,T];C^{0,\kappa}\Gamma),\ \kappa\in(0,1),\\
				\rho^n&\rightharpoonup\rho&&\text{ in }C_w([0,T];L^\gamma(B)),\\
				Z^n&\rightharpoonup Z&&\text{ in }C_w([0,T]; L^\beta(B)),\\
				u^n&\rightharpoonup u&&\text{ in }L^2(0,T;W^{1,q}(B)),\ q\in[1,2)\text{ if }\min\{\gamma,\beta\}>2,\ q=2\text{ if }\min\{\gamma,\beta\}=2
			\end{alignedat}
		\end{equation} 
		the pairs $(\rho,u)$ and $(Z,u)$ solve continuity equations in $(0,T)\times B$ and 
		\begin{equation}\label{ZeroLim}
			\lim_{n\to\infty}\int_{B}d^n|a^n-a|^p(t,\cdot)=0
		\end{equation}
		for any $p\in [1,\infty)$ and $t\in[0,T]$, where $a^n=\frac{b^n}{d^n}$ and $a=\frac{b}{d}$ keeping in mind convention \eqref{conv}.\\
		Moreover, if
		\begin{equation}\label{CompAssum}
			\text{for any } n\in\eN\  (\rho^n,Z^n)\in \overline{\mathcal O_{\underline a}}\text{ a.e. in }(0,T)\times B,\ (\rho,Z)\in \overline{\mathcal O_{\underline a}}\text{ a.e. in }(0,T)\times B,
		\end{equation}
		then
		\begin{equation}\label{ZeroLimFr}
			\lim_{n\to\infty}\int_{B}\rho^n|s^n-s|^p(t,\cdot)=0
		\end{equation}
		for any $t\in[0,T]$ and $p\geq 1$, where we define $s^n(t,x)=\frac{Z^n(t,x)}{\rho^n(t,x)}$, $s(t,x)=\frac{Z(t,x)}{\rho(t,x)}$.\\[2.mm]
		Further notice that when the comparability \eqref{CompAssum} is satisfied the requirements of \eqref{valqgb} and \eqref{ConvergN}$_{4}$ can be replaced by
		\begin{equation}\label{valqgb2}
			q\in[1,2)\text{ if } \max\{\gamma,\beta\}>2\text{ and } q=2\text{ if } \max\{\gamma,\beta\}=2.
		\end{equation}

		\begin{proof}
			Using the fact that $W^{2,2}(\Gamma)$ is compactly embedded in $C^{0,\kappa}(\Gamma)$ for any $\kappa\in(0,1)$ and the continuous embedding of $C^{0,\kappa}(\Gamma)$ into $L^2(\Gamma)$ we conclude \eqref{ConvergN}$_1$ by the Aubin--Lions lemma.
			The existence of a nonrelabeled sequence $\{(\rho^n,Z^n)\}$ with a limit $(\rho, Z)$ satisfying \eqref{ConvergN}$_{2,3}$ follows from \eqref{BddAss}. We note that details of the proof of \eqref{ConvergN}$_{2,3}$ can be found in \cite[Section 7.10.1]{NovStr04}. Convergence \eqref{ConvergN}$_4$ follows immediately from \eqref{BddAss}.\\
			Further in view of \eqref{ConvergN}$_{2,3}$ and the bound of $\{u^{n}\}$ from \eqref{BddAss} one can conclude that
			$$(\rho^{n}u^{n},Z^{n}u^{n})\rightharpoonup^{*} (\rho u,Z u)\,\,\mbox{in}\,\, L^{\infty}(0,T;L^{\frac{2\gamma}{\gamma+1}}(B))\times L^{\frac{2\beta}{\beta+1}}(B))$$
			(the proof follows the same line of arguments used later while showing \eqref{RUZUTWeakSt}). This convergence is sufficient for the passage $n\rightarrow \infty$ in the continuity equations solved by $(\rho^{n},u^{n})$ and $(Z^{n},u^{n})$ to conclude that the pairs $(\rho,u)$ and $(Z,u)$ as well solve continuity equations in $(0,T)\times B.$\\
			For the proof of \eqref{ZeroLim}
			it is necessary to show that $\frac{(b^n)^2}{d^n}$, where $b^n=\rho^n$ or $b^n=Z^n$ and $d^n=\rho^n+Z^n$, for any $n\in\eN$ as well as the limits $\frac{b^2}{d}$, where $b=\rho$ or $b=Z$ and $d=\rho+Z$, satisfy the time integrated renormalized continuity equation up to the boundary. Obviously, for $r=(b,d)$ the function $\mathcal B(r)=\frac{b^2}{d}$ does not fulfill the assumptions in Lemma~\ref{Lem:Renormalization}. Therefore one has to employ the latter lemma with the function $\mathcal B_\sigma(r)=\frac{b^2}{d+\sigma}$ for $\sigma>0$ and then use the Lebesgue dominated convergence theorem for the limit passage $\sigma\to 0_+$ to conclude that \eqref{TIREq} holds with $\eta=\eta^n$ and $\mathcal B((b^n,d^n))=\frac{(b^n)^2}{d^n}$ for any $n$ and $\mathcal B((b, d))=\frac{b^2}{d}$ as well.
			
			Fixing $t\in[0,T]$ we have
			\begin{equation}\label{LimS}
				\begin{split}
					\lim_{n\to\infty}\int_{B} d^n(t)(a^n-a)^2(t)=&\lim_{n\to\infty}\int_{B}d^n(t)(a^n)^2(t)-2\lim_{n\to\infty}\int_{B}d^n(t)a^n(t)a(t)+\lim_{n\to\infty}\int_{B}d^n(t)(a(t))^2\\
					=&\sum_{j=1}^3 I_j.
				\end{split}
			\end{equation}
			By the definition of $I_1$ and $a^n$, we get employing \eqref{TIREq} for $\eta=\eta^n$, $B((b^n,d^n))=\frac{(\rho^n)^2}{d^n}$ with $\phi=1$ and assumption \eqref{InVCnv}
			\begin{equation*}
				I_1=\lim_{n\to\infty}\int_{\Omega_{\eta^n}(t)}\frac{(b^n)^2(t)}{d^n(t)}=\lim_{n\to\infty}\int_{\Omega_{\eta^n}(0)}\frac{(b^n)^2(0)}{d^n(0)}=\int_{\Omega_{\eta_0}}\frac{b_0^2}{d_0}.
			\end{equation*}
			Next, thanks to \eqref{ConvergN}$_{2,3}$, the definition of $a^n, a$ and equation \eqref{TIREq} for $\mathcal B(b,d)=\frac{b^2}{d}$ with $\phi=1$, we deduce
			\begin{equation*}
				\begin{split}
					I_2=&-2\lim_{n\to\infty}\int_{B}b^n(t)a(t)=-2\int_{B}b(t)a(t)=-2\int_{\Omega_{\eta}(t)}\frac{b^2(t)}{d(t)}=-2\int_{\Omega_{\eta_0}}\frac{b^2_0}{d_0},\\
					I_3=&\lim_{n\to\infty}\int_{B}d^n(t)a^2(t)=\int_{B}d(t)a^2(t)=\int_{\Omega_{\eta}(t)}\frac{b^2(t)}{d(t)}=\int_{\Omega_{\eta_0}}\frac{b^2_0}{d_0}.
				\end{split}
			\end{equation*}
			Hence going back to \eqref{LimS} we conclude 
			\begin{equation}\label{ZeroLimSqr}
				\lim_{n\to\infty}\int_{B} d^n(t)(a^n-a)^2(t)=0.
			\end{equation}
			Using \eqref{ZeroLimSqr} and the fact that $a^n-a$ is bounded by definition, \eqref{ZeroLim} for $p>2$ immediately follows. Moreover, using the H\"older inequality along with \eqref{ZeroLimSqr} and the bound on $\{d^n\}$ in $L^\infty(0,T;L^1(B))$ following from assumption \eqref{BddAss} we deduce \eqref{ZeroLim} also for $p<2$
			\begin{equation*}
				\begin{split}
					\lim_{n\to\infty}\int_{B} d^n|a^n-a|^p(t,\cdot)=&\lim_{n\to\infty}\int_{B} (d^n)^\frac{p}{2}|a^n-a|^p(t,\cdot)(d^n)^{1-\frac{p}{2}}\\
					\leq&\lim_{n\to\infty}\left(\int_{B} d^n(a^n-a)^2(t,\cdot)\right)^\frac{p}{2}\left(\int_{B}d^n\right)^\frac{2-p}{2}=0,
				\end{split}
			\end{equation*}
			which concludes \eqref{ZeroLim}.
			
			We now focus on the proof of \eqref{ZeroLimFr}. Let us observe that \eqref{ZeroLim} with $p=1$ implies for any $t\in[0,T]$ that
			\begin{equation}\label{PointCnv}
				\begin{split}
					\left(\rho^n-(\rho^n+Z^n)\frac{\rho}{\rho+Z}\right)(t,\cdot)\to 0\text{ in }L^1(B),\\
					\left(Z^n-(\rho^n+Z^n)\frac{Z}{\rho+Z}\right)(t,\cdot)\to 0\text{ in }L^1(B).
				\end{split}
			\end{equation}
			Next, we show
			\begin{equation}\label{AuxPConv}
				\left(Z^n-\rho^n\frac{Z}{\rho}\right)(t,\cdot)\to 0\text{ in }L^1(B)
			\end{equation}
			for any $t\in[0,T]$.
			We rewrite
			\begin{equation*}
				Z^n-\rho^n\frac{Z}{\rho}=Z^n-(\rho^n+Z^n)\frac{Z}{\rho+Z}-\left(\rho^n-(\rho^n+Z^n)\frac{\rho}{\rho+Z}\right)\frac{Z}{\rho}
			\end{equation*}
			and deduce \eqref{AuxPConv} employing \eqref{PointCnv}. Taking into account the assumed bound on $\{s^n-s\}$ in $L^\infty((0,T)\times B)$ one concludes \eqref{ZeroLimFr} from \eqref{AuxPConv}. 
		\end{proof}
	\end{lem}
	\subsection{Non-linear Koiter energy and estimates for the structure}\label{Sec:KE}
	The following discussion on the non-linear Koiter shell energy and its properties is a summary of \cite[Section 4]{MuhaSch} and \cite{Breit2} (some of them are inspired from the reference literature \cite{CiarletIII}).\\
	The non-linear Koiter model is given in terms of the difference of the first and the second fundamental forms of $\Sigma_{\eta}$ and $\Gamma.$ We recall that $\nu_{\eta}$ denotes the normal-direction to the deformed middle surface $\varphi_{\eta}(\Gamma)$ at the point $\varphi_{\eta}(x)$ and is given by
	$$\nu_{\eta}(x)=\partial_{1}\varphi_{\eta}(x)\times \partial_{2}\varphi_{\eta}(x)={\bf{a}}_{1}(\eta)\times{\bf{a}_{2}}(\eta).$$
	%%%%%%%%%%%%%%%%%%%%%%%%%%%%%%%%%%%%%%%%%%%%%%%%%%%%%%%%%%%%%%%%
	In view of \eqref{varphieta}, these tangential derivatives $\bf{a}_{i}(\eta)$ can be computed as follows
	\begin{equation}\label{compdertan}
		\begin{array}{ll}
			\displaystyle {\bf{a}}_{i}(\eta)=\partial_{i}\varphi_{\eta}={{a}}_{i}+\partial_{i}\eta\nu+\eta\partial_{i}\nu,\quad\mbox{in}\quad i\in\{1,2\},
		\end{array}
	\end{equation}
	where ${{a}}_{i}=\partial_{i}\varphi(x).$\\
	Hence the components of the first fundamental form of the deformed configuration are given by
	\begin{equation}\label{fff}
		\begin{array}{ll}
			\displaystyle a_{ij}(\eta)={\ba}_{i}(\eta)\cdot{\ba}_{j}(\eta)
			=a_{ij}+\partial_{i}\eta\partial_{j}\eta+\eta({a}_{i}\cdot\partial_{j}\nu+{a}_{j}\cdot\partial_{i}\nu)+\eta^{2}\partial_{i}\nu\cdot\partial_{j}\nu,
		\end{array}
	\end{equation}
	where $a_{ij}=\partial_{i}\varphi(x)\cdot\partial_{j}\varphi(x).$\\
	%%%%%%%%%%%%%%%%%%%%%%%%%%%%%%%%%%%%%%%%%%%%%%%%%%%%%%%%%%%%%%%%
	Now in order to introduce the elastic energy $K=K(\eta)$ associated with the non-linear Koiter model we will use the description presented in \cite{MuhaSch} (which is inspired from \cite{Ciarlet}).\\
	In order to introduce the Koiter shell energy we first define two quantities $\mathbb{G}(\eta)$ and $\mathbb{R}(\eta).$ The change of metric tensor $\mathbb{G}(\eta)=(G_{ij}(\eta))_{i,j}$ is defined as follows
	\begin{equation}\label{Geta}
		\begin{array}{ll}
			\displaystyle G_{ij}(\eta)&\displaystyle=\partial_{i}\varphi_{\eta}\cdot\partial_{j}\varphi_{\eta}-\partial_{i}\varphi\cdot\partial_{j}\varphi=a_{ij}(\eta)-a_{ij}\\
			&\displaystyle =\partial_{i}\eta\partial_{j}\eta+\eta({a}_{i}\cdot\partial_{j}\nu+{a}_{j}\cdot\partial_{i}\nu)+\eta^{2}\partial_{i}\nu\cdot\partial_{j}\nu.
		\end{array}
	\end{equation}
	Further we define the tensor $\mathbb{R}(\eta)=(R_{ij}(\eta))_{i,j}$ which is a variant of the second fundamental form to measure the change of curvature
	\begin{equation}\label{Rsharp}
		\begin{array}{ll}
			\displaystyle R_{ij}(\eta)=\frac{\partial_{ij}\varphi_{\eta}\cdot \nu_{\eta}}{|\partial_{1}\varphi\times\partial_{2}\varphi|}-\partial_{ij}\varphi\cdot\nu=\frac{1}{|{a}_{1}\times {a}_{2}|}\partial_{i}{\ba}_{j}(\eta)\cdot\nu_{\eta}-\partial_{i}a_{j}\cdot\nu,\quad i,j=1,2.
		\end{array}
	\end{equation}
	Next, the elasticity tensor is defined as
	\begin{equation}\label{elastictytensor}
		\begin{array}{ll}
			&\displaystyle \mathcal{A}\mathbb{E}=\frac{4\lambda_{s}\mu_{s}}{\lambda_{s}+2\mu_{s}}(\mathbb{A}:\mathbb{E})\mathbb{A}+4\mu_{s}\mathbb{A}\mathbb{E}\mathbb{A},\qquad \mathbb{E}\in \mbox{Sym}(\mathbb{R}^{2\times 2}),
		\end{array}
	\end{equation}
	where $\mathbb{A}$ is the contravariant metric tensor associated with $\partial\Omega$ and $\lambda_{s},\mu_{s}>0$ are the Lam\'{e} coefficients. The Koiter energy of the shell is given by:
	\begin{equation}\label{Koiterenergy}
		\displaystyle K(\eta)=\frac{h}{4}\int_{\Gamma}\mathcal{A}\mathbb{G}(\eta(\cdot,t)):\mathbb{G}(\eta(\cdot,t))+\frac{h^{3}}{48}\int_{\Gamma}\mathcal{A}\mathbb{R}(\eta(\cdot,t))\otimes \mathbb{R}(\eta(\cdot,t)),
	\end{equation} 
	where $h>0$ is the thickness of the shell.\\ 
	Now in view of the Koiter energy \eqref{Koiterenergy} we write (following \cite{MuhaSch}) the elasticity operator $K'(\eta)$ as follows
	\begin{equation}\label{elasticityoperator}
		\begin{array}{ll}
			&\displaystyle \langle K'(\eta),b\rangle = a_{G}(t,\eta,b)+a_{R}(t,\eta,b),\qquad\forall b\in W^{2,p}(\Gamma)\,\,\mbox{where}\,\, p>2.
		\end{array}
	\end{equation}
	In the previous expression $a_{G}(t,\eta,b)$ and $a_{R}(t,\eta,b)$ are defined respectively as
	\begin{equation}\label{amab}
		\begin{split}
			&a_{G}(t,\eta,b)=\frac{h}{2}\int_{\Gamma}\mathcal{A}\mathbb{G}(\eta(\cdot,t)):\mathbb{G}'(\eta(\cdot,t))b,\\
			&a_{R}(t,\eta,b)=\frac{h^{3}}{24}\int_{\Gamma}\mathcal{A}\mathbb{R}(\eta,\cdot,t):\mathbb{R}'(\eta(\cdot,t))b,
		\end{split}
	\end{equation}
	where $\mathbb{G}'$ and $(\mathbb{R})'$ represent respectively the Fr\'{e}chet derivative of $\mathbb{G}$ and $\mathbb{R}.$ It is important to know the structure of $a_{G}(t,\eta,b)$ and $a_{R}(t,\eta,b),$ which will play a key role during the limit passages in suitably constructed approximate equations. Since $G_{ij}(\eta)$ is given by \eqref{Geta}, $G'_{ij}(\eta)b$ can simply be calculated as
	\begin{equation}\label{Gij'}
		G_{ij}'(\eta)b=\partial_{i}b\partial_{j}\eta+\partial_{i}\eta\partial_{j}b+b(a_{i}\cdot\partial_{j}\nu+a_{j}\cdot\partial_{i}\nu)+2\eta b\partial_{i}\nu\cdot\partial_{j}\nu.
	\end{equation}
	Hence one checks that $a_G(t,\eta,b)$ is a polynomial in $\eta$ and $\nabla\eta$ of order three and further the coefficients are in $L^{\infty}(\Gamma).$\\
	As in \cite[Section 4.1.]{MuhaSch}, $R_{ij}(\eta)$ (introduced in \eqref{Rsharp}) can be written in the following form which is easier to handle
	\begin{equation}\label{Rijcomp}
		R_{ij}(\eta)=\overline{\gamma}(\eta)\partial^{2}_{ij}\eta+P_{0}(\eta,\nabla\eta),
	\end{equation}
	where $P_{0}$ is a polynomial of order three in $\eta$ and $\nabla\eta$ such that all terms are at most quadratic in $\nabla\eta$ and the coefficients of $P_{0}$ depend on $\varphi$ and the geometric quantity $\overline{\gamma}(\eta)$ (depending on $\partial\Omega$ and $\eta$) is defined as follows
	\begin{equation}\label{ovgamma}
		\overline{\gamma}(\eta)=\frac{1}{|a_{1}\times a_{2}|}\bigg(|a_{1}\times a_{2}|+\eta(\nu\cdot(a_{1}\times\partial_{2}\nu+\partial_{1}\nu\times a_{2}))+\eta^{2}\nu\cdot(\partial_{1}\nu\times\partial_{2}\nu)\bigg).
	\end{equation}
	Hence $R_{ij}(\eta)$ can be written as follows
	\begin{equation}\label{rewriteRij}
		R_{ij}'(\eta)b=\overline{\gamma}(\eta)\partial^{2}_{ij}b+(\overline{\gamma}'(\eta)b)\partial^{2}_{ij}\eta+P'_{0}(\eta,\nabla\eta)b,
	\end{equation}
	$i.e.$ we have
	\begin{equation}\label{exaR}
		\begin{split}
			a_{R}(t,\eta,b)=&\frac{h^{3}}{24}\int_{\Gamma}\bigg[\mA(\overline{\gamma}(\eta)\nabla^{2}\eta):(\overline{\gamma}(\eta)\nabla^{2}b)+\mA(\overline{\gamma}(\eta)\nabla^{2}\eta):(\overline{\gamma}'(\eta)b\nabla^{2}\eta)\\
			&+\mA(\overline{\gamma}(\eta)\nabla^{2}\eta):P'_{0}(\eta,\nabla\eta)b+\mA(P_{0}(\eta,\nabla\eta)):(\overline{\gamma}(\eta)\nabla^{2}b)\\
			&+\mA(P_{0}(\eta,\nabla\eta)):(\overline{\gamma}'(\eta)b\nabla^{2}\eta)+\mA(P_{0}(\eta,\nabla\eta)):(P'_{0}(\eta,\nabla\eta))b\bigg].
		\end{split}
	\end{equation}
	One notices that $\overline{\gamma}'(\eta)$ is a linear in $\eta$ and $\overline\gamma({\eta})$ is quadratic in $\eta.$
	
	%%%%%%%%%%%%%%%%%%%%%%%%%%%%%%%%%%%%%%%%%%%%%%%%%%%%%%%%%%%%%%
	%%%%%%%%%%%%%%%%%%%%%%%%%%%%%%%%%%%%%%%%%%%%%%%%%%%%%%%%%%%%%%%%%%%%%%%%%%%%%%%%%%%%%%%%%%%%%%%%%%%%%%%%%%%%%%%%%%%%%%%%%%%%%%%%%%%%%%%%%%%%%%%%%%%%%%%%%%%%%%%%%%%%%%%%%%%%%%%%%%%%%%%%%%%%%%%%%%%%%%%%%%%%%%%%%%%%%%%%%%%%%%%%%%%%%%%%%%%%%%%%%%%%%%%%%%%%%%%%%%%%%%%%%%%%%%%%%%%%%%%%%%%%%%%%%%%%%%%%%%%%%%%%%%%%%%%%
	\section{Artificial regularization, extension of the problem in a larger domain, further approximations of the pressure and data}\label{Sec:ExtProb}
	In this section we first introduce a regularization of the shell energy and next suitable approximations of the viscosity coefficients, pressure and the initial data.\\
	The regularization of the shell energy is needed to solve a structural sub-problem (cf. Section \ref{Sec:StrSub}), more precisely to obtain suitable compactness properties of the structural displacement.\\
	On the other hand suitable approximations of the Lam\'{e} coefficients and the initial data play a crucial role to first solve a dummy problem in a smooth larger domain and then to return back to the physical domain by means of suitable limit passages.
	%%%%%%%%%%%%%%%%%%%%%%%%%%%%%%%%%%%%%%%%%%%%%%%%%%%%%%%%%%%%%%%%%%%%%%%%%%%%%%%%%%%%%%%%%%%%%%%%%%%%%%%%%%%%%%%%%%%%%%%%%%%%%%%%%%%%%%%%%%%%%%%%%%%%%%%%%%%%%%%%%%%%%%%%%%%%%%%%%%%%%%%%%%%%%%%%%%%%%%%%%%%%%%%%%%%%%%%%%%%%%%%%%%%%%%%%%%%%%%%%%%%%%%%%%%%%%%%%%%%%%%%%%%%%%%%%%%%%%%%%%%%%%%%%%%%%%%%%%%%%%%%%%%%%%%%%%%%%%%%%%%%%%%%%%%%%%%%%%%%%%%%%%%%%%%%%%%%%%%%%%%%%%%%%%%%%%%%%%%%%%%%%%%%%%%%%%%%%%%%%%%%%%%%%
	\subsection{Regularization of the shell energy and artificial dissipation}\label{regshell}
	Let us introduce a regularization of the shell energy as follows
	\begin{equation}\label{regshellenergy}
		K_{\delta}(\eta)=K(\eta)+\delta^{7}\int_{\Gamma}|\nabla^{3}\eta|^{2}.
	\end{equation}
	In connection with the regularization above, we further regularize the initial condition for the structural displacement, $i.e.$ we consider a sequence $\{\eta_{0}^{\delta}\}_{\delta}\subset W^{3,2}(\Gamma)$ such that
	\begin{equation}\label{regintdens}
		\begin{split}
			&\eta_{0}\in W^{2,2}(\Gamma);\ \eta_{0}^{\delta}\rightarrow \eta_{0}\mbox{ in } W^{2,2}(\Gamma),\,\, \mbox{and  }\delta^{7}\int_\Gamma|\nabla^{3}\eta_0^{\delta}|^2\rightarrow 0\,\,\mbox{as}\,\,\delta\rightarrow 0.
		\end{split}
	\end{equation}
	Considering an arbitrary function $\omega\in C^\infty(\Gamma)$ we define an approximate identity $\omega_{\delta}(\cdot)=\delta^{-2}\omega\left(\frac{\cdot}{\delta}\right)$ and construct the sequence $\{\eta_0^{\delta}\}$ via $\eta^{\delta}_0=\eta_0*\omega_{\delta}$. Then \eqref{regintdens} follows by arguments similar to the proof of \cite[Theorem 7.38]{Adams75}. Taking into account obvious inequalities $\|\nabla^{3}\eta_0^{\delta}\|_{L^2(\Gamma)}\leq\|\eta_0\|_{L^2(\Gamma)}\|\nabla^{3}\omega_{\delta}\|_{L^1(\Gamma)}$ and $\|\nabla^{3}\omega_{\delta}\|_{L^1(\Gamma)}\leq\delta^{-3}\|\nabla^{3}\omega\|_{L^1(\Gamma)}$ and \eqref{regintdens}$_{3}$ follows.\\
	Further such a regularization of the shell energy makes the boundary Lipschitz in space for a.e. time which justifies some arguments involving integration by parts. Those justifications are comparatively intricate to perform in H\"{o}lder domains.\\ 
	In a strong form the evolution of the structure \eqref{fluidsysreform}$_{4},$ now takes the following form 
	\begin{equation}
		\begin{array}{ll}
			&\displaystyle \partial^{2}_{t}\eta+K'_{\delta}(\eta)-\zeta\partial_{t}\Delta \eta=F\cdot\nu\,\,\mbox{on}\,\,\Gamma\times I,
		\end{array}
	\end{equation}
	for some positive parameter $\zeta$ and in view of \eqref{regshellenergy}, $K'_{\delta}(\eta)$ can be defined as
	\begin{equation}\label{K'ep0}
		\begin{array}{ll}
			\langle K'_{\delta},b\rangle =\langle K'(\eta), b\rangle+\delta^{7}\langle \nabla^{3}\eta,\nabla^{3}b\rangle,\,\,\mbox{for all}\,\,b\in \str.
		\end{array}
	\end{equation}
	The dissipation (more specifically the $W^{1,2}((0,T)\times\Gamma)$ regularity of the structure when $\zeta>0$) plays an important role in Lemma \ref{Lem:Fund}. Indeed we will pass with the dissipation parameter $\zeta$ to zero (using some uniform estimates) in Section \ref{reglim0} for the case $\max\{\gamma,\beta\}>2$ and $\min\{\gamma,\beta\}>0$. Whereas to handle the case $\max\{\gamma,\beta\}=2$ and $\min\{\gamma,\beta\}>0$ we will need the dissipation $\zeta$ to be positive.\\
	As the first part of constructing an approximate solution we extend the weak formulation of our problem \eqref{fluidsysreform} in a larger domain. In other words we first embed our physical domain $\Omega_{\eta}(t)$ in $B=B_{m,M},$ where $B$ is the neighborhood of $\Omega$ introduced in \eqref{defB}. By construction \eqref{tbnbd}-\eqref{defB}, $\Omega_{\eta}(t)\Subset B=B_{m,M}
	$ for all $t\in[0,T]$ since $m\le\eta\le M.$ Indeed such a pair $(m,M)$ (and consequently $B=B_{m,M}$) can be fixed from the beginning since $\eta$ is bounded as a consequence of the energy estimates.
	
	To begin with we extend the viscosity coefficients and the data of the problem as explained in the following section. 
	\subsection{Extension of coefficients and data}\label{extendatacoeff}
	For a fixed $0<\omega\ll 1$ we approximate the viscosity coefficients $\mu$ and $\lambda$ from \eqref{SDef} by
	\begin{equation}\label{ViscAppr}
		\mu^\eta_\omega:=f^\eta_\omega\mu,\ \lambda^\eta_\omega:=f^\eta_\omega\lambda,    
	\end{equation}
	where the function $f^\eta_\omega\in C^\infty_c([0,T]\times\eR^3)$ satisfies
	\begin{equation}\label{ExtendingFProp}
		\begin{split}
			0<\omega\leq f^\eta_\omega\leq 1\text{ in }[0,T]\times B,\\
			f^\eta_\omega(t,\cdot)|_{\Omega_\eta}=1\text{ for all }t\in[0,T],\\
			\|f^\eta_\omega\|_{L^p(((0,T)\times B)\setminus Q^T_{\eta}}\leq c\omega\text{ for some }p\geq 1,\\
			\text{ the mapping }\eta\mapsto f^\eta_\omega\text{ is Lipschitz.}
		\end{split}
	\end{equation}
	The function $f^\eta_\omega$ is defined via 
	\begin{equation}
		f^\eta_\omega(t,X):=f_\omega((\tilde {\varphi}_\eta)^{-1}(t,X)),
	\end{equation}
	where we set
	\begin{equation}
		f_\omega=\chi_{\Omega}+\chi_{B\setminus \Omega}g_\omega
	\end{equation}
	with a suitable cut-off function $g_\omega\in C^\infty_c(\eR^3)$ satisfying
	\begin{equation}
		g_\omega(x)\begin{cases}
			=1&\text{ if }x\in\partial\Omega,\\
			\in(2\omega, 1]&\text{ if }0<\dist(x,\partial\Omega)<\omega,\\
			\in (\omega, 2\omega)&\text{ if }\dist(x,\partial\Omega)\geq\omega.
		\end{cases}
	\end{equation}
	We note that the first three properties in \eqref{ExtendingFProp} immediately follow by the definition of $f_\omega^\eta$, whereas the fourth property is a consequence of the regularity of $f_\omega$ and the Lipschitz continuity of the mapping $\eta\mapsto (\tilde{\varphi}_\eta)^{-1}$ defined in \eqref{FlowMInv}.
	
	We introduce an approximate the pressure $P(\rho,Z).$ In that direction we consider $\delta>0$ and a sufficiently large $\kappa\gg \max\{4,\gamma,\beta\}$.
	The approximation of the pressure is defined as 
	\begin{equation}\label{approxpressure}
		P_{\delta}(\rho,Z)=P(\rho,  Z)+\delta\left(\rho^{\kappa}+Z^{\kappa}+\frac{1}{2}\rho^{2}Z^{\kappa-2}+\frac{1}{2}Z^{2}\rho^{\kappa-2}\right).
	\end{equation}
	
	The initial data $\rho_{0},$ $Z_{0}$ and $M_{0}$ are extended and approximated in $B$ such that the approximating functions $\rho_{0,\delta},$ $Z_{0,\delta}$ and $M_{0,\delta}$ satisfy
	\begin{equation}\label{extensionindata}
		\begin{split}
			&\rho_{0,\delta},\ Z_{0,\delta}\geqslant 0,\ \rho_{0,\delta}|_{\eR^3\setminus\Omega_{\eta^\delta_0}}= Z_{0,\delta}|_{\eR^3\setminus\Omega_{\eta^\delta_0}}=0,\ \rho_{0,\delta},\ Z_{0,\delta}\not\equiv 0,\\
			&(\rho_{0,\delta}(x), Z_{0,\delta}(x))\in\overline{\mathcal O_{\underline a}}\text{ for a.a. }x\in B,\ \rho_{0,\delta}, Z_{0,\delta}\in L^\kappa(B),\ M_{0,\delta}\in L^{\frac{2\kappa}{\kappa+1}}(B),\\
			&\rho_{0,\delta}\rightarrow \rho_{0}\text{ in }L^\gamma(\Omega_{\eta_0}),\ Z_{0,\delta}\rightarrow Z_{0}\text{ in } L^\beta(\Omega_{\eta_0}), M_{0,\delta}\to M_0\text{ in }L^1(\Omega_{\eta_0}),\\& \int_B \frac{|M_{0,\delta}|^2}{\rho_{0,\delta}+Z_{0,\delta}}\to \int_{\Omega_0}\frac{|M_0|^2}{(\rho_0+Z_0)},\ \delta\int_B\left(|\rho_{0,\delta}|^\kappa+|Z_{0,\delta}|^\kappa\right)\to 0\text{ as }\delta\to 0.
		\end{split}
	\end{equation}
	The interested reader can consult \cite[Section 7.10.7]{NovStr04}, where the analogous regularization of initial data is performed for the single--fluid case.\\
	Next we define the notion of weak solution for a bi-fluid system considered on the fixed domain $B$ containing the structure.
	
	\subsection{Definition of weak solution in the extended set up}\label{Def:EPSol}
	The weak solution for the extended problem ($i.e.$ for a system defined in $B$) is defined as follows.
	\begin{mydef}\label{WSExtProb}
		The quadruple $(\rho,Z,u,\eta)$ is a bounded energy weak solution to the extended problem in $B$  if 
		\begin{align*}
			&\rho,\ Z\ge 0\text{ a.e. in }B,\\ &\rho\in L^{\infty}(0,T;L^{\kappa}(B)),\\
			& Z\in L^{\infty}(0,T;L^{\kappa}(B)),\\
			& u\in L^{2}(0,T;W^{1,2}_0(B)),\\
			&(\rho+Z)|u|^{2}\in L^{\infty}(0,T;L^{1}(B)),\\
			& P_\delta(\rho,Z)\in L^{1}((0,T)\times B)\\
			& \eta\in L^{\infty}(0,T;W^{2,2}(\Gamma)\cap W^{1,\infty}(0,T;L^{2}(\Gamma)),\\
			&\eta\in L^{\infty}(0,T;\delta^\frac{7}{2}W^{3,2}(\Gamma)),\\
			& \eta\in W^{1,2}(0,T;\sqrt{\zeta}W^{1,2}(\Gamma))
		\end{align*}
		and the following hold.
		\begin{enumerate}[leftmargin=5ex, label=(\roman*), topsep=-1.5ex, itemsep=1ex]
			\item The coupling of $u$ and $\tder\eta$ reads $\Tr_{\Sigma_\eta} u=\tder\eta\nu$, where the operator $\Tr_{\Sigma_\eta}$ is defined in Lemma \ref{Lem:TrOp}.		
			\item The momentum equation is satisfied in the sense
			\begin{equation}\label{momentumex}
				\begin{split}
					&\int_{(0,t)\times B}(\rho+Z)u\cdot\tder\phi+\int_{(0,t)\times B}\left((\rho+Z)u\otimes u\right)\cdot \nabla\phi-\int_{(0,t)\times B} \mathbb{S}^{\eta}_{\omega}(\mathbb D u)\cdot\nabla\phi\\
					&+\int_{(0,t)\times B}P_\delta(\rho,Z)\mathrm{div}\,\phi+\int_{(0,t)\times \Gamma}\tder\eta \tder b-\int^{t}_{0}\langle K'_{\delta}(\eta),b\rangle+\zeta\int_{(0,t)\times\Gamma}\partial_{t}\nabla\eta\nabla b
					\\&=\int_{B}(\rho+Z)u(t,\cdot)\phi(t,\cdot)+\int_{\Gamma}\partial_{t}\eta(t,\cdot)b(t,\cdot)-\int_{B}M_{0,\delta}\phi(0,\cdot)-\int_{\Gamma}\eta_1 b(0,\cdot)
				\end{split}
			\end{equation}
			for a.a. $t\in(0,T)$ and all $(b,\phi)\in \str\times C^{\infty}([0,T]\times \mathbb{R}^{3})$ with $tr_{\eta}\phi=b\nu$ and $\mathbb{S}^{\eta}_{\omega}(\mathbb D u)$ is defined as follows
			\begin{equation}\label{SDefetaomega}
				\mathbb{S}^{\eta}_{\omega}(\mathbb D u)=2\mu^{\eta}_{\omega}\left(\mathbb D u-\frac{1}{3}\dvr u\mathbb{I}_{3}\right)+\lambda^{\eta}_{\omega}\dvr u\mathbb{I}_{3}.
			\end{equation}
			We recall that in \eqref{momentumex}, the regularized Koiter energy $K_{\delta}$ is given by \eqref{Koiterenergy}-\eqref{regshellenergy}, the approximate pressure $P_{\delta}(\cdot,\cdot)$ is as introduced in \eqref{approxpressure} and $M_{0,\delta}$ is defined in \eqref{extensionindata}. 
			\item The continuity equations are satisfied in the sense
			\begin{equation}\label{contrhoex}
				\begin{split}
					\int_B\left(\rho(t,\cdot)\psi(t,\cdot)-\rho_{0,\delta}\psi(0,\cdot)\right)=&\int_{(0,t)\times B}\rho(\tder\psi+u\cdot\nabla\psi),\\
					\int_B\left(Z(t,\cdot)\psi(t,\cdot)- Z_{0,\delta}\psi(0,\cdot)\right)=&\int_{(0,t)\times B}Z(\tder\psi+u\cdot\nabla\psi)
				\end{split}
			\end{equation}
			for all $t\in[0,T]$ and all $\psi\in C^{\infty}([0,T]\times \mathbb{R}^{3})$. Where one recalls that the approximated initial densities $\rho_{0,\delta}$ and $Z_{0,\delta}$ were introduced in \eqref{extensionindata}.
			\item The energy inequality
			\begin{equation}\label{energybalanceex}
				\begin{split}
					&\int_B\bigg(\frac{1}{2}(\rho+Z)|u|^{2}+{\mathcal{H}}_{P,\delta}(\rho,Z)\bigg)(\cdot,t)+\int_0^t\int_B\mathbb{S}^{\eta}_{\omega}(\mathbb D u)\cdot \nabla u+\bigg(\int_{\Gamma}\frac{1}{2}|\partial_{t}\eta|^{2}+K_{\delta}(\eta)\bigg)(\cdot,t)\\
					&+\zeta\int_{\Gamma\times I}|\nabla\partial_{t}\eta|^{2}\le\int_B\bigg(\frac{|M_{0,\delta}|^{2}}{2(\rho_{0,\delta}+Z_{0,\delta})}+{\mathcal{H}}_{P,\delta}(\rho_{0,\delta},Z_{0,\delta})\bigg)+\bigg(\frac{1}{2}\int_{\Gamma}|\eta_{1}|^{2}+K_{\delta}(\eta_{0}^{\delta})\bigg)
				\end{split}
			\end{equation}
			holds for a.a. $t\in I,$ where
			\begin{equation}\label{HPdelta}
				\displaystyle {\mathcal{H}}_{P,\delta}(\rho,Z)=H_{P_{\delta}}(\rho,Z)+h_{\delta}(\rho,Z).
			\end{equation}
			In \eqref{HPdelta}, $h_{\delta}(\cdot,\cdot)$ is defined as
			\begin{equation}\label{hdelta}
				h_{\delta}(\rho,Z)=\frac{\delta}{\kappa-1}\left(\rho^{\kappa}+Z^{\kappa}+\frac{1}{2}\rho^{2}Z^{\kappa-2}+\frac{1}{2}Z^{2}\rho^{\kappa-2}\right)
			\end{equation}
			and $H_{P_{\delta}}(\cdot,\cdot)$ is as defined in \eqref{HelmFDef}.
		\end{enumerate}
	\end{mydef}
	
	\begin{thm}\label{resultstaulayer}
		Assume that the Hypotheses \textbf{H1}--\textbf{H5} hold. Further we recall the artificial regularization of the shell energy, added structural dissipation from Section \ref{regshellenergy} and the extension of the Lam\'{e} coefficients, initial data and pressure regularization from Section \ref{extendatacoeff}. Then there is $T\in(0,\infty]$ and a weak solution to the extended problem in the sense of Definition~\ref{WSExtProb} on $(0,T)$. The time $T$ is finite only if 
		\begin{equation*}
			\text{either }\lim_{s\to T}\eta(s,y)\searrow a_{\partial\Omega}\text{ or }\lim_{s\to T}\eta(s,y)\nearrow b_{\partial\Omega}
		\end{equation*}		
		for some $y\in\Gamma$.
	\end{thm}
	The proof of Theorem \ref{resultstaulayer} involves two crucial steps:\\
	\begin{itemize}
		\item Splitting of \eqref{momentumex} into two sub-problems (a decoupling penalization technique), namely the fluid and structural problems and  to prove their existence separately. All of these splitting and penalization is done at an approximate level where the time interval $(0,T)$ is divided into sub-intervals of length $\tau$ (the approximation parameter). Considering the length  we have decided to devote an entire section (Section \ref{sppen}) for this part of our analysis.
		\item Next we pass $\tau$ to zero and prove Theorem \ref{resultstaulayer}. This is done in Section \ref{sec:TauCnv}, more precisely in Section \ref{proofthmtaulevel}.
	\end{itemize}

	%%%%%%%%%%%%%%%%%%%%%%%%%%%%%%%%%%%%%%%%%%%%%%%%%%%%%%%%%%%%%%%%%%%%%%%%%%%%%%%%%%%%%%%%%%%%%%%%%%%%%%%%%%%%%%%%%%%%%%%%%%%%%%%%%%%%%%%%%%%%%%%%%%%%%%%%%%%%%%%%%%%%%%%%%%%%%%%%%%%%%%%%%%%%%%%%%%%%%%%%%%%%%%%%%%%%%%%%%%%%%%%%%%%%%%%%%%%%%%%%%%%%%%%%%%%%%%%%%%%%%%%%
	\section{The splitting and penalized problem}\label{sppen}		
	Let us divide the time interval into $N\in\mathbb{N}$ sub-intervals of length $\tau=\frac{T}{N}$.\\
	As a first step to solve the extended problem in $B$ (as introduced in \ref{momentumex}), we solve two sub-problems corresponding to the fluid part and the elastic part. The splitting does not preserve the kinematic coupling condition. Instead we will include penalization terms in the weak formulations of the decoupled equations and this will ensure the recovery of the interface couplings as $\tau\rightarrow 0.$ We further introduce an auxiliary unknown $v$ representing the trace of the fluid velocity on the interface, more precisely
	\begin{equation}\label{vinterface}
		v=\Tr_{\Sigma_\eta}u.
	\end{equation}
	\subsection{The structural sub-problem}\label{Sec:StrSub}
	The notion of weak solution for the structural sub-problem will be defined inductively. More precisely for $n\geqslant 0,$ it is defined as follows
	\begin{equation}\label{n0}
		\begin{split}
			{\text{For }}\ n=0:&\ \eta^{0}(0,\cdot)= \eta_{0}^{\delta},\ \partial_{t}\eta^{0}(0,\cdot)=\eta_{1}\text{ such that}\\
			&v^{0}(x+\eta^{0}\nu(x),t)=\eta_{1}\nu(x)\text{ on }\Gamma\text{ for }t\in[-\tau,0];
		\end{split}
	\end{equation}
	and for $n\geqslant 1,$ we assume the existence of $\eta^{n}$ and also the existence of a solution $(\rho^n, Z^n, u^n)$ of the fluid sub-problem and solve for $\eta^{n+1}$ such that:
	\begin{enumerate}[label=\arabic*., topsep=-1.5ex, itemsep=1ex]
		\item $\eta^{n+1}\in W^{1,\infty}(n\tau,(n+1)\tau;L^{2}(\Gamma))\cap L^{\infty}(n\tau,(n+1)\tau;W^{2,2}(\Gamma))\cap L^{\infty}(n\tau,(n+1)\tau;\sqrt{\delta}W^{3,2}(\Gamma))$\\
		$\,\,\,\,\cap W^{1,2}(n\tau,(n+1)\tau;\sqrt{\zeta}W^{1,2}(\Gamma)),$\\
		\item $\eta^{n+1}(\cdot,n\tau)=\eta^{n}(\cdot,n\tau)$, $\partial_{t}\eta^{n+1}(\cdot,n\tau)=\partial_{t}\eta^{n}(\cdot,n\tau)$ in the weakly continuous sense in time.
		\item The following structural equation
		\begin{equation}\label{structuralpen}
			\begin{split}
				&(1-\delta)\int_{n\tau}^{(n+1)\tau}\int_{\Gamma}\partial_{t}\eta^{n+1}\partial_{t}b-\delta\int_{n\tau}^{(n+1)\tau}\int_{\Gamma}\frac{\partial_{t}\eta^{n+1}-v^{n}\cdot\nu}{\tau}b+\zeta\int^{(n+1)\tau}_{n\tau}\int_{\Gamma}\partial_{t}\nabla\eta^{n+1}\nabla b\\
				&-\int_{n\tau}^{(n+1)\tau}\langle K'_{\delta}(\eta^{n+1}),b\rangle=(1-\delta)\int_{n\tau}^{(n+1)\tau}\frac{d}{dt}\int_{\Gamma}\partial_{t}\eta^{n+1} b
			\end{split}
		\end{equation}
		holds for all $b\in L^{\infty}(n\tau,(n+1)\tau;W^{3,2}(\Gamma))\cap W^{1,\infty}(n\tau,(n+1)\tau;L^{2}(\Gamma))$, where $v^n=\Tr_{\Sigma_\eta^n}u^n$.
		\item The following energy like inequality
		\begin{equation}\label{energytypestr}
			\begin{split}
				&\frac{\delta}{2\tau}\int_{n\tau}^t\left(\|\partial_{t}\eta^{n+1}-v^{n}\cdot\nu\|^{2}_{L^{2}(\Gamma)}+\|\partial_{t}\eta^{n+1}\|^{2}_{L^{2}(\Gamma)}\right)+\zeta\int^{t}_{n\tau}\|\partial_{t}\nabla\eta^{n+1}\|_{L^{2}(\Gamma)}^{2}\\
				&+\frac{1-\delta}{2}\|\partial_{t}\eta^{n+1}(t)\|^{2}_{L^{2}(\Gamma)}+ K_{\delta}(\eta^{n+1})(t)\\
				& \leqslant \frac{1-\delta}{2}\|\partial_{t}\eta^{n+1}(n\tau)\|^{2}_{L^{2}(\Gamma)}+K(\eta^{n+1}(n\tau))+\frac{\delta}{2\tau}\int_{n\tau}^{t}\|v^n\|^{2}_{L^{2}(\Gamma)}
			\end{split}
		\end{equation}
		holds for all $t\in(n\tau,(n+1)\tau].$
	\end{enumerate}
	\subsection{The fluid sub-problem}\label{Sec:FlSub} Similarly to the structural sub-problem, the notion of solution to the fluid sub-problem is also defined using induction as follows
	\begin{equation}\label{n0rhou}
		\mbox{For }n=0:\,\,\rho^{0}(0,\cdot)= \rho_{0,\delta},\ Z^0(0,\cdot)=Z_{0,\delta},\,((\rho+Z)u)^{0}(\cdot,0)=M_{0,\delta};
	\end{equation}
	and for $n\geqslant 1,$ we assume the existence of $(\rho^{n},u^{n})$ and solve for $(\rho^{n+1},u^{n+1})$ such that:
	\begin{enumerate}[label=\arabic*., topsep=-1.5ex, itemsep=1ex]
		\item $\rho^{n+1}, Z^{n+1}\geqslant 0$, $\rho^{n+1}, Z^{n+1}\in L^{\infty}(n\tau,(n+1)\tau;L^{\kappa}(B))$,
		$u^{n+1}\in L^{2}(n\tau,(n+1)\tau;W^{1,2}_0(B))$,\\ $(\rho^{n+1}+Z^{n+1})|u^{n+1}|^{2}\in L^{\infty}(n\tau,(n+1)\tau;L^{1}(B))$,
		\item $\rho^{n+1}(n\tau)=\rho^{n}(n\tau)$, $Z^{n+1}(n\tau)=Z^{n}(n\tau)$,  $((\rho+Z)u)^{n+1}(n\tau)=(\rho u)^{n}(n\tau)$ in weakly continuous sense in time.
		\item
		The continuity equations of the form
		\begin{equation}\label{contdicouplevel}
			\begin{split}
				\int_{n\tau}^{(n+1)\tau}\frac{d}{dt}\int_{B}\rho^{n+1}\psi-\int_{n\tau}^{(n+1)\tau}\int_{B}\bigg(\rho^{n+1}\partial_{t}\psi+\rho^{n+1} u^{n+1}\cdot\nabla\psi\bigg)&=0,\\
				\int_{n\tau}^{(n+1)\tau}\frac{d}{dt}\int_{B}Z^{n+1}\psi-\int_{n\tau}^{(n+1)\tau}\int_{B}\bigg(Z^{n+1}\partial_{t}\psi+Z^{n+1} u^{n+1}\cdot\nabla\psi\bigg)&=0
			\end{split}
		\end{equation}
		hold for all $\psi\in C^{\infty}([n\tau,(n+1)\tau]\times \RR^{3})$.
		\item The following momentum equation
		\begin{equation}\label{momentumdecouple}
			\begin{split}
				&\int_{n\tau}^{(n+1)\tau}\int_{B}(\rho^{n+1}+Z^{n+1})\bigg(u^{n+1}\cdot\partial_{t}\phi+(u^{n+1}\otimes u^{n+1})\cdot \nabla\phi\bigg)\\
				&-\int_{n\tau}^{(n+1)\tau}\int_{B} \mathbb{S}^{n+1}_\omega(\mathbb D u^{n+1})\cdot\nabla\phi
				+\int_{n\tau}^{(n+1)\tau}\int_{B}P_{\delta}(\rho^{n+1},Z^{n+1})\mathrm{div}\,\phi
				\\&-\delta\int_{n\tau}^{(n+1)\tau}\int_{\Gamma}\frac{v^{n+1}-\partial_{t}\eta^{n+1}\nu}{\tau}\cdot b=\int_{n\tau}^{(n+1)\tau}\frac{d}{dt}\int_{B}(\rho^{n+1}+Z^{n+1})u^{n+1}\cdot\phi
			\end{split}
		\end{equation}
		holds for all $(b,\phi)\in L^{\infty}(n\tau,(n+1)\tau;W^{3,2}(\Gamma))\cap W^{1,\infty}(n\tau,(n+1)\tau;L^{2}(\Gamma))\times C^{\infty}([n\tau,(n+1)\tau]\times \RR^{3})$ with $\Tr_{\Sigma_{\eta^{n+1}}}\phi=b\nu$ where 
		\begin{align}
			v^{n+1}&=\Tr_{\Sigma_{\eta^{n+1}}} u^{n+1},\label{vn+1un+1}\\
			\mathbb{S}^{n+1}_\omega(\mathbb D u^{n+1})&=2\mu_\omega^{\eta^{n+1}}\left(\mathbb D u^{n+1}-\frac{1}{3}\dvr u^{n+1}\mathbb I_3\right)+\lambda^{\eta^{n+1}}_\omega\dvr u^{n+1}\mathbb I_3\nonumber\\
		\end{align}
		and the viscosity coefficients $\mu^\eta_\omega$, $\lambda^\eta_\omega$ are defined in \eqref{ViscAppr}.
		\item The following energy inequality
		\begin{equation}\label{energybalancedecouple}
			\begin{split}
				&\int_{B}\bigg(\frac{1}{2}(\rho^{n+1}+Z^{n+1})|u^{n+1}|^{2}+{\mathcal{H}}_{P,\delta}(\rho^{n+1},Z^{n+1})\bigg)(t)\\
				&+\int_{n\tau}^{t}\int_{B}\mathbb{S}^{n+1}_\omega(\mathbb D u^{n+1})\cdot\nabla u^{n+1}+\frac{\delta}{2\tau}\int_{n\tau}^{t}\int_{\Gamma}\bigg(|v^{n+1}-\partial_{t}\eta^{n+1}\cdot\nu|^{2}+|v^{n+1}|^{2}\bigg)\\
				&  \leqslant \int_{B}\bigg(\frac{1}{2}(\rho^{n}+Z^{n})|u^{n}|^{2}+\mathcal{H}_{P,\delta}(\rho^{n},Z^{n})\bigg)(n\tau)+\frac{\delta}{2\tau}\int_{n\tau}^{t}\int_{\Gamma}|\partial_{t}\eta^{n+1}|^{2}
			\end{split}
		\end{equation}
		holds for a.a. $t\in [n\tau,(n+1)\tau],$ where
		\begin{equation}\label{HPdelta*}
			{\mathcal{H}}_{P,\delta}(\rho,Z)=H_{P_{\delta}}(\rho,Z)+\delta(\rho^\kappa +Z^\kappa+\frac{1}{2}\rho^{\kappa-2}Z^2+\frac{1}{2}\rho^2Z^{\kappa-2}).
		\end{equation}
		and $H_{P_{\delta}}$ is as defined in \eqref{HelmFDef}.
	\end{enumerate}
	%%%%%%%%%%%%%%%%%%%%%%%%%%%%%%%%%%%%%%%%%%%%%%%%%%%%%%%%%%%%%%%%%%%%%%%%%%%%%%%%%%%%%%%%%%%%%%%%%%%%%%%%%%%%%%%%%%%%%%%%%%%%%%%%%%%%%%%%%%%%%%%%%%%%%%%%%%%%%%%%%%%%%%%%%%%%%%%%%%%%%%%%%%%%%%%%%%%%%%%%%%%%%%%%%%%%%%%%%%%%%%%%%%%%%%%%%%%%%%%%%%%%%%%%%%%%%%%
	\subsection{Existence of solution for the sub-problems}
	In this section we will present results on the existence
	of $\eta^{n+1}$ solving the structural sub problem \eqref{structuralpen} along with the estimate \eqref{energytypestr} and $(\rho^{n+1},Z^{n+1},u^{n+1})$
	solving the fluid sub problem \eqref{contdicouplevel}--\eqref{momentumdecouple} along with the estimate \eqref{energybalancedecouple}. The first theorem concerns the existence of solution to the structural sub-problem
	\begin{thm}\label{labelstructuralsp}
		Let $(\eta^{n},\partial_{t}\eta^{n})(\cdot,0)=(\eta^{n\tau},\eta^{n\tau}_{1})\in W^{3,2}(\Gamma)\times L^{2}(\Gamma).$ Further let $v^{n}\in L^{2}(\Gamma).$ Then for $n\in \mathbb{N}\cup \{0\}$ and a positive $\tau<1,$ the problem \eqref{structuralpen} admits of a solution $\eta^{n+1}$ such that $(\eta^{n+1},\partial_{t}\eta^{n+1})\in W^{3,2}(\Gamma)\times L^{2}(\Gamma).$
	\end{thm}
	The proof of Theorem \ref{labelstructuralsp} borrows ideas from \cite[Section 6]{MuhaSch} and \cite[Section 5.2]{CanicMuha}. Since we are using a different scheme to decouple the fluid and the structural sub-problems (one recalls the operator splitting scheme used in \cite{CanicMuha} and \cite{MuhaSch}) we have a penalization term appearing in the weak formulation of the structural sub-problem. Further we have an extra visco-elastic term which appears with a parameter $\zeta.$ Appearance of these terms needs some modified adaptations of the arguments used in \cite{CanicMuha} and \cite{MuhaSch}.  Further for an application of fixed point argument we need to use different functional spaces compared to the ones used in \cite{CanicMuha}. Hence we prefer to provide an independent proof of Theorem \ref{labelstructuralsp} in Section \ref{proofdecstr}.\\[2.mm]
	%%%%%%%%%%%%%%%%%%%%%%%%%%%%%%%%%%%%%%%%%%%%%%%%%%%%%%%%%%%%%%%%%%%%%%%%%%%%%%%%%%%%%%%%%%%%%%%%%%%%%%%%%%%%%%%%%%%%%%%%%%%%%%%%%%%%%%%%%%%%%%%%%%%%%%%%%%%%%%%%%%%%%%%%%%%%%%%%%%%%%%%%%%%%%%%%%%%%%%%%%%%%%%%%%%%%%%%%%%%%%%%%%%%%%%%%%%%%%%%%%%%%%%%%%%%%%%%%%%%%%%%%%%%%%%%%%%%%%%%%%%%%%%%%%%%%%%%%%%%%%%%%%%%%%%%%%%%%%%%%%%%%%%%%%%%%%%%%%%%%%%%%%%%%%%%%%%%%%%%%%%%%%%%%%%%%%%%%%%%%%%%%%%%%%%%%%%%%%%%%%%%%%%%%%%%%%%%
	The next theorem corresponds to solving the fluid sub problem in a fixed domain of class $C^{2}.$ 
	\begin{thm}\label{existencediscretefluid}
		Let $(\rho_{0,\delta}, Z_{0,\delta}, M_{0,\delta})$ satisfy \eqref{extensionindata}, $\eta^{n+1}$ solves the items $1$ and $2$ of Section \ref{Sec:StrSub}. Further let hypotheses $(H1-H5)$ hold. Then for $\tau>0$ there exists at least one weak solution $(\rho^{n+1},Z^{n+1},u^{n+1})$ solving \eqref{contdicouplevel}--\eqref{momentumdecouple} in an iterative manner. Moreover, inequality \eqref{energybalancedecouple} holds and for all $t\in [n\tau,(n+1)\tau]$ and almost all $x\in B$ $(\rho^{n+1}(t,x),Z^{n+1}(t,x))\in\overline{\mathcal O_{\underline a}}$.
	\end{thm}
	For the proof of Theorem \ref{existencediscretefluid}, we refer to \cite[Theorem 1, p. 365]{NovoPoko}. Note that for a fixed $\tau>0$ the weak formulation of the momentum equation \eqref{momentumdecouple} differs slightly from that of \cite[p. 364, (26)]{NovoPoko}. First we have an extra term $\displaystyle\delta\int_{n\tau}^{(n+1)\tau}\int_{\Gamma}\frac{v^{n+1}-\partial_{t}\eta^{n+1}\nu}{\tau}\cdot b$ which is of lower order and hence can be handled with minor modifications. Secondly, the viscosity coefficients $\mu$ and $\lambda$ in \cite{NovoPoko} are assumed to be constants whereas in our case $\mu^{\eta^\tau}_\omega$ and $\lambda^{\eta^\tau}_\omega$ are functions. Even this does not cause any problem to adapt arguments from the proof of \cite{NovoPoko} because of the non-degeneracy construction \eqref{ExtendingFProp}$_{1}.$ 
	
	\subsection{Uniform bounds on approximate solutions and weak formulations}
	Let us assume that structure and fluid sub--problems have been solved during the iteration process described in subsections \ref{Sec:StrSub} and \ref{Sec:FlSub} for fixed $N\in\eN$ and $\{(\rho^n,Z^n,u^n,\eta^n)\}$ be a sequence of corresponding solutions.
	
	For the purposes of this subsection we use the notation
	\begin{equation}\label{piecewdef}
		f^{\tau}(t):=f^{n+1}(t)\text{ for }t\in(n\tau,(n+1)\tau],
	\end{equation}
	where $\tau=\frac{T}{N}$. $f^{n+1}$ stands for one of the functions $\rho^{n+1}$, $Z^{n+1}$, $\eta^{n+1}$, $u^{n+1}$. Accordingly, $f^{\tau}$ stands for one of the functions $\rho^{\tau}$, $Z^{\tau}$, $\eta^{\tau}$, $u^{\tau}$. Moreover, we set $\mathbb S^{\eta^{\tau}}_\omega:=\mathbb S^{\eta^{n+1}}_\omega$.  Using the energy inequalities for the decoupled sub--problems we derive the total energy inequality. To this end we fix $m\in\{1,\ldots,N-1\}$. Setting $t=(n+1)\tau$ in \eqref{energytypestr}, \eqref{energybalancedecouple} respectively, summing over $n\in\{0,\ldots, m-2\}$ and adding \eqref{energytypestr}, \eqref{energybalancedecouple} with $t\in [(m-1)\tau,m\tau]$, we obtain
	\begin{equation}\label{UBTau}
		\begin{split}
			&\int_B\left(\frac{1}{2}({\rho^{\tau}}+{Z^{\tau}})|u^{\tau}|^2+\mathcal{H}_{P,\delta}({\rho^{\tau}},{Z^{\tau}})\right)(t)+\frac{1-\delta}{2}\|\tder{\eta^{\tau}}(t)\|^2_{L^2(\Gamma)}+K_{\delta}({\eta^{\tau}})(t)+\zeta\int_0^t\int_\Gamma|\tder\nabla{\eta^{\tau}}|^2\\&+\int_0^t\int_B \mathbb S^{\eta^{\tau}}_\omega(\mathbb D u^{\tau})\cdot\nabla u^{\tau}+\frac{\delta}{2\tau}\int_0^t\left(\|\tder{\eta^{\tau}}- v^{\tau}(\cdot-\tau)\cdot\nu\|^2_{L^2(\Gamma)}+\|v^{\tau}-\tder{\eta^{\tau}}\nu\|^2_{L^2(\Gamma)}\right)+\frac{\delta}{2\tau}\int^t_{t-\tau}\|v^{\tau}\|^2_{L^2(\Gamma)}\\
			&\leq \int_B\left(\frac{|M_{0,\delta}|^2}{2({\rho^{\tau}}_{0,\delta}+{Z^{\tau}}_{0,\delta})}+\mathcal H_{P,\delta}({\rho^{\tau}}_{0,\delta}, {Z^{\tau}}_{0,\delta})\right)+\frac{1-\delta}{2}\|{\eta_{1}}\|^2_{L^2(\Gamma)}+ K_{\delta}({\eta}_0^\delta)+\frac{\delta}{2}\|v^0\|^2_{L^2(\Gamma)}.
		\end{split}
	\end{equation}
	We point out that the relations $v^{\tau}(s-\tau)=v^n(s-\tau)$ for $s\in (n\tau,(n+1)\tau]$ if $n\geqslant 1$ and $v^{\tau}(s-\tau)=v^0$ if $s\in [0,\tau]$ being in accordance with \eqref{piecewdef} and the notation
	\begin{equation}\label{SEODef}
		\mathbb{S}^{\eta^{\tau}}_\omega(\mathbb Du^{\tau})=2\mu^{\eta^{\tau}}_\omega\left(\mathbb D u^{\tau}-\frac{1}{3}\dvr u^{\tau}\mathbb I_3\right)+\lambda_\omega^{\eta^{\tau}}\dvr u^{\tau}\mathbb I_3  
	\end{equation}
	were also used. Based on inequality \eqref{UBTau} and the Korn inequality the functions $({\eta^{\tau}}, {\rho^{\tau}}, {Z^{\tau}}, u^{\tau})$ satisfy the following estimates 
	\begin{equation}\label{TauEst}
		\begin{split}
			\delta^\frac{1}{2} \|v^{\tau}-\tder{\eta^{\tau}}\nu\|_{L^2(0,T;L^2(\Gamma))}\leq& c\tau^\frac{1}{2},\\
			\zeta^\frac{1}{2}\|\tder{\eta^{\tau}}\|_{L^2(0,T;W^{1,2}(\Gamma))}+\delta^{\frac{1}{2}}\|{\eta^{\tau}}\|_{L^\infty(0,T;W^{3,2}(\Gamma))}\leq &c,\\
			\delta^\frac{1}{\kappa}\left(\|{\rho^{\tau}}\|_{L^\infty(0,T;L^\kappa(B))}+\|{Z^{\tau}}\|_{L^\infty(0,T;L^\kappa(B))}\right)\leq& c,\\
			\|\sqrt{{\rho^{\tau}}+{Z^{\tau}}}u^{\tau}\|_{L^\infty(0,T;L^2(B))}\leq &c,\\
			\|u^{\tau}\|_{L^2(0,T;W^{1,2}(B))}\leq &c\omega^{-\frac{1}{2}}.
		\end{split}
	\end{equation}
	Further, it directly follows from Theorem \ref{existencediscretefluid} that
	\begin{equation}\label{DensComparT}
		(\rho^{\tau}(t,x),Z^{\tau}(t,x))\in\overline{\mathcal O_{\underline a}} \text{ for all }t\in(0,T)\text{ and almost all }x\in B.
	\end{equation}
	Moreover, by interpolation of estimates \eqref{TauEst}$_{3,4,5}$ and the Sobolev embedding of $W^{1,2}(B)$ into $L^6(B)$ we get
	\begin{equation}\label{AddTauEst}
		\begin{split}
			\|({\rho^{\tau}}+{Z^{\tau}}) u^{\tau}\|_{L^\infty(0,T;L^\frac{2\kappa}{\kappa+1}(B))}&\leq c,\\
			\|({\rho^{\tau}}+{Z^{\tau}}) u^{\tau}\otimes u^{\tau}\|_{L^2(0,T;L^\frac{6\kappa}{4\kappa+3}(B))}&\leq c.
		\end{split}
	\end{equation}
	Since the interface is uniform in time Lipschitz (recall the bound from the second summand of \eqref{TauEst}), one uses an argument involving the classical Bogovskii operator (for instance we refer to \cite[Section 4.4.]{NovoPoko}) to furnish the following
	\begin{equation}\label{improvedpressurestawayint}
		\int^{T}_{0}\int_{{B\setminus {\Sigma_{{\eta^{\tau}}}}}}(({\rho^{\tau}})^{\gamma+1}+({Z^{\tau}})^{\beta+1}+\delta(({\rho^{\tau}})^{\kappa+1}+({Z^{\tau}})^{\kappa+1}))\leqslant C
	\end{equation} 
	where $C>0$ is independent of $\tau$. To be precise for the proof of \eqref{improvedpressurestawayint} we first use test functions of the form $(\phi,b)=(\psi\mathfrak B({\rho^{\tau}}-[{\rho^{\tau}}]_{\Omega_{{\eta^{\tau}}}}),0)$ and $(\phi,b)=(\psi\mathfrak B({\rho^{\tau}}-[{\rho^{\tau}}]_{B\setminus\Omega_{{\eta^{\tau}}}}),0)$ in \eqref{momentumdecouple} 
	(where $\mathfrak B$ is the Bogovskii operator for the domains $\Omega_{{\eta^{\tau}}},$ $B\setminus\Omega_{{\eta^{\tau}}}$ respectively and $\psi\in C^1_c((0,T))$) and next repeat the arguments with the test functions $(\psi\mathfrak B({Z^{\tau}}-[{Z^{\tau}}]_{\Omega_{{\eta^{\tau}}}}),0)$ and $(\psi\mathfrak B({Z^{\tau}}-[{Z^{\tau}}]_{B\setminus\Omega_{{\eta^{\tau}}}}),0).$ For the definition of Bogovskii operator $\mathfrak B$ we refer the readers to \cite[Chapter 3.]{NovStr04} and further for the obtainment of \eqref{improvedpressurestawayint} we use similar calculations as that of \cite[Section 4.3.]{NovoPoko}. It is possible to obtain a constant $C>0$ independent of $\tau$ in \eqref{improvedpressurestawayint} since $\Sigma_{{\eta^{\tau}}}$ is uniform in time Lipschitz and hence norm of the linear operator $\mathfrak B=\mathfrak B_{\tau}$ (where $\mathfrak B_{\tau}$ corresponds to $\Omega_{{\eta^{\tau}}}$) is independent of $\tau$ (see for instance \cite[Lemma 4.1]{BuFe09} and the remark that follows).

	By the above iterative procedure we obtain functions $({\rho^{\tau}}, {Z^{\tau}}, u^{\tau})$ satisfying the bounds from \eqref{UBTau} on the interval $[0,T]$. For further analysis we also need the weak formulations of the continuity equations and coupled momentum equation, which are satisfied by $({\rho^{\tau}}, {Z^{\tau}}, u^{\tau})$. We obtain directly from \eqref{contdicouplevel} that
	\begin{equation}\label{ConEqAppr}
		\begin{split}
			\int_{B}\left({\rho^{\tau}}(t,\cdot)\psi(t,\cdot)-{\rho^{\tau}}_{0,\delta}\psi(0,\cdot)\right)&=\int_0^t\int_{B}\bigg({\rho^{\tau}}\partial_{t}\psi+{\rho^{\tau}} u^{\tau}\cdot\nabla\psi\bigg),\\
			\int_B\left( {Z^{\tau}}(t,\cdot)\psi(t,\cdot)- {Z^{\tau}}_{0,\delta}\psi(0,\cdot)\right)&=\int_0^t\int_{B}\bigg({Z^{\tau}}\partial_{t}\psi+{Z^{\tau}} u^{\tau}\cdot\nabla\psi\bigg)
		\end{split}
	\end{equation}
	hold for $t\in[0,T]$ and all $\psi\in C^{\infty}([0,T]\times \RR^{3})$. To obtain the weak formulation of the coupled momentum equation we fix 
	\begin{equation}\label{AdmTest}
		\phi\in C^\infty([0,T]\times \RR^{3}),\ b\in\str, b\nu=\Tr_{\Sigma_{\eta^\tau}} \phi.
	\end{equation}
	Let $t\in[0,T]$ be fixed. We first find $m\in\eN$ such that $t\in[m\tau,(m+1)\tau).$
	Next we add \eqref{structuralpen} tested by $b|_{[n\tau,(n+1)\tau]},$ \eqref{momentumdecouple} tested by $(\phi|_{[n\tau,(n+1)\tau]},b|_{[n\tau,(n+1)\tau]}),$ sum the resulting identity over $n=1,\ldots, m-1.$ Once again adding the resulting expression with \eqref{structuralpen} tested by $b|_{[m\tau,t]}$ and \eqref{momentumdecouple} tested by $(\phi|_{m\tau,t]},b|_{[m\tau,t]})$ to conclude
	\begin{equation}\label{MomEqAppr}
		\begin{split}
			&\int_0^t\int_B \left(({\rho^{\tau}}+{Z^{\tau}})\left(u^{\tau}\tder \phi + (u^{\tau}\otimes u^{\tau})\cdot\nabla \phi\right)+P_\delta({\rho^{\tau}},{Z^{\tau}})\dvr\phi -\mathbb S^{\eta^{\tau}}_\omega(\mathbb D u^{\tau})\cdot\nabla\phi\right)\\
			&\quad-\delta\int_0^t\int_\Gamma \frac{(v^{\tau}-v^{\tau}(\cdot-\tau))\cdot\nu}{\tau}b+(1-\delta)\int_0^t\int_\Gamma \tder{\eta^{\tau}} \tder b -\int_0^t\langle K'_{\delta}({\eta^{\tau}}),b\rangle-\zeta\int^{t}_{0}\int_{\Gamma}\partial_{t}\nabla{\eta^{\tau}}\nabla b\\
			&=\int_B({\rho^{\tau}}+{Z^{\tau}}) u^{\tau}(t,\cdot)\phi(t,\cdot)-\int_BM_{0,\delta}\cdot\phi(0,\cdot)+(1-\delta)\left(\int_\Gamma\tder{\eta^{\tau}}(t,\cdot)b(t,\cdot)-\int_\Gamma{\eta}_1b(0,\cdot)\right)
		\end{split}
	\end{equation}
	for any $t\in[0,T]$ and any pair $(\phi, b)$ satisfying \eqref{AdmTest}.\\
	The next subsection contains a lemma on the extension of the densities in a larger time independent domain. 
	\subsection{On the extension of density in a time independent domain}
	The following lemma deals with the solution of the continuity equation and states that the fluid densities vanish outside the physical domain $\Omega_{\eta}$ if they vanishes initially outside $\Omega_{\eta_0}$.
	\begin{lem}\label{Lem:Fund}
		Let $\rho,Z\in L^\infty(0,T;L^3(B))$, $u\in L^2(0,T;W^{1,2}_0(B))$ satisfy the continuity equation \eqref{contrhoex} with the initial condition $\rho_{0,\delta}$, $Z_{0,\delta}$ respectively, given in \eqref{extensionindata}. Let Assumptions (A) hold with $\eta\in W^{1,\infty}(0,T;L^2(\Gamma))\cap W^{1,2}((0,T)\times\Gamma)\cap L^\infty(0,T;C^{0,1}(\Gamma))$ and $u(x+\eta(t,\varphi^{-1}(x))\nu(x))=\tder\eta(t,\varphi^{-1}(x)))\nu(x)$ hold on $(0,T)\times\partial\Omega$ in the sense of traces. Then it follows that
		\begin{equation*}
			\rho|_{B\setminus \Omega_\eta(t)}=Z|_{B\setminus \Omega_\eta(t)}\equiv 0\text{ for a.a. }t\in(0,T).
		\end{equation*}
	\end{lem}
	We refer to the proof of Lemma \ref{Lem:Fund} presented in the appendix, Section \ref{extension0proof}.\\
	The Lemma \ref{Lem:Fund} will be used in the upcoming sections and this will also play a crucial role to come back to the physical domain $\Omega_{\eta}$ from the extended domain $B.$\\
	%%%%%%%%%%%%%%%%%%%%%%%%%%%%%%%%%%%%%%%%%%%%%%%%%%%%%%%%%%%%%%%%%%%%%%%%%%%%%%%%%%%%%%%%%%%%%%%%%%%%%%%%%%%%%%%%%%%%%%%%%%%%%%%%%%%%%%%%%%%%%%%%%%%%%%%%%%%%%%%%%%%%%%%%%%%%%%%%%%%%%%%%%%%%%%%%%%%%%%%%%%%%%%%%%%%%%%%%%%%%%%%%%%%%%%%%%%%%%%%%%%%%%%%%%%%%%%%%%%%%%%%%%%%%%%%%%%%%%%%%%%%%%%%%%%%%%%%%%%%%%%%%%%%%%%%%%%%%%%%%%%%%%%%%%%%%%%%%%%%%%%%%%%%%%%%%%%%%%%%%%%%%%%%%%%%%%%%%%%%%%%%%%%%%%%%%%%%%%%%%%%%%%%%%%%%%%%%%%%%%%%%%
	\section{Limit passage \texorpdfstring{$\tau\to 0_+$}{} layer}\label{sec:TauCnv}
	The goal of this section is the limit passage in equations \eqref{ConEqAppr} and \eqref{MomEqAppr} and to prove Theorem \ref{resultstaulayer}.\\
	We continue by stating the convergences of the approximates $(\eta^{\tau},\rho^{\tau},Z^{\tau},u^{\tau})$ in the following subsection.
	\subsection{Convergence of the approximates and some consequences  }\label{Sec:TauLim} To this end we use the following convergences that are direct consequence of \eqref{TauEst} (note that at this stage we are only interested on bounds independent on $\tau$ and hence in the following convergences we do not specify explicitly the dependence of the spaces on other parameters $\delta$ and $\zeta$)
	\begin{equation}\label{TauConv}
		\begin{alignedat}{2}
			\eta^\tau&\rightharpoonup^*\eta&&\text{ in }L^\infty(0,T;W^{3,2}(\Gamma)),\\
			\tder\eta^\tau&\rightharpoonup^*\tder\eta&&\text{ in }L^\infty(0,T;L^2(\Gamma))\cap L^{2}(0,T;W^{1,2}(\Gamma)),\\
			\rho^\tau&\rightharpoonup^*\rho&&\text{ in }L^\infty(0,T;L^\kappa(B)),\\
			Z^\tau&\rightharpoonup^*Z &&\text{ in }L^\infty(0,T;L^\kappa(B)),\\
			u^\tau&\rightharpoonup u &&\text{ in }L^2(0,T;W^{1,2}(B)).\\
		\end{alignedat}
	\end{equation}
	Since 
	\begin{equation}\label{EpsDepEmb}
		\begin{split}
			\|\eta(t)-\eta(s)\|_{(L^2(\Gamma),W^{3,2}(\Gamma))_{\theta,2}}&\leq \|\eta(t)-\eta(s)\|^\theta_{W^{3,2}(\Gamma)}\|\eta(t)-\eta(s)\|^{1-\theta}_{L^2(\Gamma)}\\
			&\leq c|t-s|^{1-\theta}\|\eta\|^\theta_{L^\infty(0,T;W^{3,2}(\Gamma))}\|\eta\|^{1-\theta}_{W^{1,\infty}(0,T;L^2(\Gamma))}
		\end{split}
	\end{equation}
	we get
	\begin{equation*}
		L^\infty(0,T;W^{3,2}(\Gamma))\cap W^{1,\infty}(0,T;L^2(\Gamma))\hookrightarrow C^{0,1-\theta}([0,T]; W^{3\theta ,2}(\Gamma))\hookrightarrow C^{0,1-\theta}([0,T]; C^{0,1}(\Gamma)),
	\end{equation*}
	for $\theta\in \left(\frac{2}{3},1\right)$. Hence we obtain due to \eqref{TauConv}$_{1,2}$ that
	\begin{equation}\label{EtaTCnv}
		\eta^\tau\to\eta\text{ in }C^{0,\frac{1}{4}}([0,T]; C^{0,1}(\Gamma)).
	\end{equation}
	Since at this level $\delta$ is fixed we obtain the following strong convergence (up to a non-relabeled subsequence) of $\eta^{\tau}$ as a consequence of \eqref{TauConv}$_{1}$ and \eqref{TauConv}$_{2}$ and the classical Aubin-Lions theorem
	\begin{equation}\label{strngconvetatau}
		\begin{array}{l}
			\eta^{\tau}\rightarrow \eta\,\,\mbox{in}\,\, L^{\infty}(0,T;W^{2,4}(\Gamma)).
		\end{array}
	\end{equation}
	Convergence \eqref{strngconvetatau} suffices to conclude that
	\begin{equation}\label{limpassageenenergyeta}
		\begin{array}{l}
			\displaystyle \int_0^t\langle K'_{\delta}(\eta^{\tau}),b\rangle\rightarrow \int_0^t\langle K'_{\delta}(\eta),b\rangle,
		\end{array}
	\end{equation}
	for any $b\in L^{\infty}(0,T;W^{3,2}(\Gamma))\cap W^{1,\infty}(0,T;L^{2}(\Gamma)),$ where we have used the structure of $K'_{\delta}(\eta)$ (cf. \eqref{elasticityoperator},\eqref{amab},\eqref{Geta},\eqref{Gij'},\eqref{exaR} and \eqref{regshellenergy}).
	
	Furthermore, it follows that the limit densities $\rho, Z$ satisfy 
	\begin{equation}\label{VanDens}
		\rho|_{B\setminus\Omega_\eta(t)}=Z|_{B\setminus\Omega_\eta(t)}\equiv 0\text{ for a.a. }t\in (0,T)
	\end{equation}
	by Lemma \ref{Lem:VanSeq}.
	
	%As a consequence of Lemma \ref{Lem:Coer} we get the existence of a minimal time where the following holds (due to the $W^{2,2}$ coercivity of the Koiter energy)
	%\begin{equation*}
	%\|\eta^\tau\|_{L^\infty(0,T;W^{2,2}(\Gamma))}\leq c
	%\end{equation*}
	%(also uniformly with respect to $\delta$),
	%which similarly to \eqref{EpsDepEmb} implies
	%\begin{equation}\label{EtaHConv}
	%\eta^\tau\to\eta\text{ in }C^{0,1-\theta}([0,T];C^{0,2\theta-1}(\Gamma)).
	%\end{equation}
	As an immediate consequence of \eqref{EtaTCnv} we get
	\begin{equation}\label{TauSetCnv}
		|(\Omega^{\eta^\tau}\setminus\Omega^\eta)\cup(\Omega^{\eta}\setminus\Omega^{\eta^\tau})|\to 0\text{ in }C^{0,\frac{1}{4}}[0,T]
	\end{equation}
	and 
	\begin{equation}\label{AuxApprCnv}
		f^{\eta^\tau}_\omega\to f^\eta_\omega\text{ in }C^{0,\frac{1}{4}}([0,T]\times\overline B)
	\end{equation}
	by \eqref{ExtendingFProp}$_4$.\\ 
	%%%%%%%%%%%%%%%%%%%%%%%%%%%%%%%%%%%%%%%%%%%%%%%%%%%%%%%%%%%%%%%%%%%%%%%%%%%%%%%%%%%%%%%%%%%%%%%%%%%%%%%%%%%%%%%%%%%%%%%%%%%%%%%%%%%%%%%%%%%%%%%%%%%%%%%%%%%%%%%%%%%%%%%%%%%%%%%%%%%%%%%%%%%%%%%%%%%%%%%%%%%%%%%%%%%%%%%%%%%%%%%%%%%%%%%%%%%%%%%%%%%%%%%%%%%%%%%%%%%%%%%%%%%%%%%%%%%%%%%%%%%%%%%%%%%%%%%
	Using the so far obtained convergences, we will pass to the limit \eqref{MomEqAppr} and further obtain a energy analogue solved by $(\rho,Z,u,\eta)$ (where $(\rho,Z,u,\eta)$ is as introduced in \eqref{TauConv}).

	The next section if devoted for the proof of Theorem \ref{resultstaulayer}.

	\subsection{Proof of Theorem \ref{resultstaulayer}}\label{proofthmtaulevel}
	
	\subsubsection{Passage $\tau\rightarrow 0_{+}$ in the non-linear terms, construction of test functions and obtaining \eqref{momentumex}:}\label{tauto0eq}
	In this section we focus on the convergences of the non-linear terms appearing in \eqref{MomEqAppr} and further conclude the proof of \eqref{momentumex}.\\[2.mm]
	$(1)$ \textit{Convergence of a term linked to convection:} 
	In a way that is now standard for the mono--fluid case we will show that
	\begin{equation}\label{CnvTPas}
		(\rho^\tau+Z^\tau) u^\tau\otimes u^\tau\rightharpoonup(\rho+Z) u\otimes u\text{ in }L^1((0,T)\times B).
	\end{equation}
	Indeed, as $\rho^\tau$ and $Z^\tau$ satisfy the equations in \eqref{ConEqAppr} and estimate \eqref{TauEst}$_{3,4}$, we deduce by the arguments based on the abstract Arzel\`a--Ascoli theorem, cf. \cite[Section 7.10.1]{NovStr04}, that
	\begin{equation}\label{RZTStrongly}
		(\rho^\tau,Z^\tau)\to(\rho,Z)\text{ in }C_w([0,T];L^\kappa(B)).
	\end{equation}
	As a consequence of the latter convergence and \eqref{DensComparT} we get
	\begin{equation}\label{DensComparD}
		(\rho(t,x),Z(t,x))\in\overline{\mathcal O_{\underline a}} \text{ for all }t\in(0,T)\text{ and almost all }x\in B.
	\end{equation}
	Using \eqref{RZTStrongly}, the compact embedding $L^\kappa(B)$ into $W^{-1,2}(B)$ and convergence \eqref{TauConv}$_{5}$ we infer
	\begin{equation}\label{RUZUTWeakSt}
		(\rho^\tau u^\tau,Z^\tau u^\tau)\rightharpoonup^*(\rho u,Zu)\text{ in }L^\infty(0,T;L^\frac{2\kappa}{\kappa+1}(B)).
	\end{equation}
	Employing momentum equation \eqref{MomEqAppr} with test functions $(\phi,b)\in C^\infty([0,T];\RR^{3})\times \str$ satisfying $b=\Tr_{\Sigma_\eta}\phi\cdot \nu$, the H\"older inequality, the estimates in \eqref{TauEst} and \eqref{AddTauEst}, the Sobolev embedding theorem, Lemma \ref{Lem:TrOp} and the trivial embedding of $L^p(B)$ into $L^p(\Omega_\eta(t))$ we conclude the uniform continuity of the sequence $\{(\rho^\tau+Z^\tau) u^\tau\}$ in $C([0,T]; W^{-3,2}(B))$. Using the abstract Arzel\`a--Ascoli theorem we infer
	\begin{equation}
		(\rho^\tau+Z^\tau) u^\tau\to(\rho +Z)u\text{ in }C([0,T];W^{-3,2}(B)).
	\end{equation}
	provided that $L^\frac{2\kappa}{\kappa-1}(B)$ is compactly embedded in $W^{-3,2}(B)$. Following the lines of the proof of \cite[Lemma 6.2]{NovStr04} after (6.1.5) we get 
	\begin{equation}\label{RZUTCnv}
		(\rho^\tau+Z^\tau) u^\tau\to(\rho +Z)u\text{ in }C_w([0,T];L^\frac{2\kappa}{\kappa+1}(B)).
	\end{equation}
	As $\frac{2\kappa}{\kappa+1}>\frac{6}{5}$, we have the compact embedding $L^\frac{2\kappa}{\kappa+1}(B)$ in $W^{-1,2}(B)$ and
	\begin{equation}\label{RZUTCnv2}
		(\rho^\tau+Z^\tau) u^\tau\to(\rho +Z)u\text{ in }L^2(0,T;W^{-1,2}(B))
	\end{equation}
	accordingly. Combining the latter convergence, \eqref{TauConv}$_5$ and the boundedness of $(\rho^{\tau}+Z^{\tau})|u^{\tau}|^{2}$ in $L^{2}(0,T;L^{\frac{6\kappa}{4\kappa+3}}(B))$ (which follows by interpolation from $(\rho^{\tau}+Z^\tau)|u^{\tau}|^{2}\in L^{\infty}(L^{1})\cap L^{1}(L^{\frac{3\kappa}{\kappa+3}})$) we conclude \eqref{CnvTPas}. We note that we can obtain that 
	\begin{equation}\label{rhou2conv}
		(\rho^\tau+Z^\tau)|u^\tau|^2\rightharpoonup (\rho+Z)|u|^2\text{ in }L^1((0,T)\times B)
	\end{equation}
	in the exactly same way as \eqref{CnvTPas}.\\[2.mm] 
	$(2)$ \textit{A convergence related to the penalization term possessing the factor $-\delta:$}
	In connection with the penalization term in \eqref{MomEqAppr}, containing a factor $-\delta$ we compute the following
	\begin{equation}\label{PenTermDecomp}
		\begin{split}
			&\int_0^T\int_0^t\int_\Gamma \frac{(v^{\tau}(s)-v^\tau(s-\tau))\cdot\nu}{\tau}b(s)\psi(t)\ds\dt=-\int_0^T\int_\tau^{t-\tau}\int_\Gamma v^\tau(s)\cdot\nu\frac{b(s+\tau)-b(s)}{\tau}\ds\psi(t)\dt\\
			&+\int_0^T\frac{1}{\tau}\int_{t-\tau}^t\int_\Gamma v^\tau(s)\cdot\nu b(s)\ds\psi(t)\dt-\frac{1}{\tau}\int_0^T\int_0^\tau\int_\Gamma \eta_1b\psi=\sum_{i=1}^3I^\tau_i.
		\end{split}
	\end{equation} 
	for any $b\in L^{\infty}(0,T;W^{3,2}(\Gamma))\cap W^{1,\infty}(0,T;L^{2}(\Gamma))$ and $\psi\in C^\infty_c((0,T))$. We immediately obtain 
	\begin{equation}\label{IT3Pass}
		I^\tau_3\to -\int_0^T\int_\Gamma \eta_1b(0)\psi(t)\dt.
	\end{equation}
	Using the convergence 
	\begin{equation}\label{PenConv}
		v^\tau-\tder\eta^\tau\nu\to 0\text{ in }L^2(0,T;L^2(\Gamma))
	\end{equation}
	coming from \eqref{TauEst}$_1$ and convergence \eqref{TauConv}$_2$ and further by using  
	$$\frac{b(\cdot+\tau)-b(\cdot)}{\tau}\rightarrow \partial_{t}b(\cdot)\,\,\mbox{in}\,\, L^{2}(0,T;L^{2}(\Gamma))$$
	(which follows from the fact that $b$ is Lipschitz continuous and hence a.e. differentiable in time with values in $L^{2}(\Gamma)$)
	we get
	\begin{equation}\label{IT1Pass}
		I^\tau_1\to-\int_0^T\int_0^t\int_\Gamma \tder\eta(s)\tder b(s)\ds\psi(t)\dt. 
	\end{equation}
	Concerning the term $I^\tau_2$, we have
	\begin{equation}
		\begin{split}
			I^\tau_2=&\int_0^T\tau^{-1}\int_{t-\tau}^t\int_\Gamma (v^\tau(s)-\tder\eta^\tau(s)\nu)\cdot\nu b(s)\ds\psi(t)\dt+\int_0^T\tau^{-1}\int_{t-\tau}^t\int_\Gamma(\tder\eta^\tau-\tder\eta)b\psi\\
			&+\int_0^T\tau^{-1}\int_{t-\tau}^t\int_\Gamma\tder\eta b\psi= \sum_{j=1}^3J^\tau_j
		\end{split}
	\end{equation}
	provided we consider the extension $\eta^{\tau}=\eta_0^\delta$ on $[-\tau,0]$. In order to pass to the limit in the term $J^\tau_1$ we define $w^\tau\in L^2(\eR;L^2(\Gamma))$ as
	\begin{equation*}
		w^\tau(s)=\begin{cases}
			v^\tau(s)-\tder\eta^\tau(s)\nu&\text{ if }s\in(0,T),\\
			0&\text{ if }s\in\eR\setminus(0,T)
		\end{cases}
	\end{equation*}
	and denote $(f)_\tau=\tau^{-1}\int_{t-\tau}^tf(s)\ds$. We infer by the Jensen inequality 
	\begin{equation*}
		\begin{split}
			\int_0^T\|(w^\tau)_\tau\|_{L^2(\Gamma)}^2\leq& \tau^{-1}\int_0^T\int_{t-\tau}^t\|w^\tau(s)\|^2_{L^2(\Gamma)}\ds \dt\leq \tau^{-1}\int_0^T\left(\int_0^t\|w^\tau(s)\|^2_{L^2(\Gamma)}\ds-\int_0^{t-\tau}\|w^\tau(s)\|^2_{L^2(\Gamma)}\ds\right)\dt\\
			=&\tau^{-1}\left(\int_0^T\int_0^t\|w^\tau(s)\|^2_{L^2(\Gamma)}\ds \dt-\int_{-\tau}^{T-\tau}\int_0^t\|w^{\tau}(s)\|^2_{L^2(\Gamma)}\ds \dt\right)\\
			\leq&\tau^{-1}\int_{T-\tau}^T\int_0^t\|w^\tau(s)\|^2_{L^2(\Gamma)}\ds \dt\leq \|w^\tau\|^2_{L^2(0,T;L^2(\Gamma))}.
		\end{split}
	\end{equation*}
	The latter inequality and \eqref{PenConv} imply
	\begin{equation*}
		(v^\tau-\tder\eta^\tau\nu)_\tau\to 0\text{ in }L^2(0,T;L^2(\Gamma))
	\end{equation*}
	by which we conclude 
	\begin{equation}\label{JT1Pass}
		J^\tau_1\to 0.
	\end{equation}
	We obtain for $J^\tau_2$ by the integration by parts
	\begin{equation}
		\begin{split}
			J^\tau_2=&\tau^{-1}\int_0^T\int_\Gamma [(\eta^\tau-\eta)b]^t_{t-\tau}\psi(t)\dt -\tau^{-1}\int_0^T\int_{t-\tau}^t\int_\Gamma(\eta^\tau-\eta)(s)\tder b(s)\ds\psi(t)\dt\\
			=& \int_0^{T-\tau}\int_\Gamma(\eta^\tau-\eta)b\frac{\psi(t)-\psi(t+\tau)}{\tau}\dt+\tau^{-1}\int_{T-\tau}^T\int_\Gamma (\eta^\tau-\eta)b\psi\dt\\
			&-\tau^{-1}\int_0^T\int_{t-\tau}^t\int_\Gamma(\eta^\tau-\eta)(s)\tder b(s)\ds\psi(t)\dt
		\end{split}
	\end{equation}
	as $\eta^\tau=\eta^\delta_0$ on $[-\tau,0]$. Employing convergence \eqref{EtaTCnv} and the fact that $\psi$ possesses a compact support in $(0,T)$ we conclude
	\begin{equation}\label{JT2Pass}
		\lim_{\tau\to 0_+}J^\tau_2=0.
	\end{equation}
	Eventually, by the Lebesgue differentiation theorem we deduce
	\begin{equation}\label{JT3Pass}
		\lim_{\tau\to 0_+}J^\tau_3=\int_0^T\int_\Gamma\tder\eta b\psi.
	\end{equation}
	Taking into consideration \eqref{IT1Pass}, \eqref{IT3Pass}, \eqref{JT1Pass}, \eqref{JT2Pass} and \eqref{JT3Pass} we deduce from \eqref{PenTermDecomp} that
	\begin{equation}\label{convdeltapen}
		\begin{split}
			&-\delta  \int_0^T\int_0^t\int_\Gamma \frac{(v^{\tau}(s)-v^\tau(s-\tau))\cdot\nu}{\tau}b(s)\ds\psi(t)\dt+(1-\delta)\int_0^T\int_0^t\int_\Gamma \tder\eta(s) \tder b(s)\ds\psi(t)\dt\\&-(1-\delta)\int^{T}_{0}\left(\int_{\Gamma}\partial_{t}\eta(t) b(t)-\int_\Gamma \eta_1b(0)\right)\psi(t)\dt\rightarrow\int^{T}_{0}\left(\int^{t}_{0}\int_{\Gamma}\partial_{t}\eta\partial_{t}b-\int_{\Gamma}\tder \eta b+\int_{\Gamma}\eta
			_{1}b(\cdot,0)\right)\psi\dt\text{ as }\tau\rightarrow 0.
		\end{split}
	\end{equation}
	(3) \textit{Convergence of the pressure:}
	Estimate \eqref{improvedpressurestawayint} (especially the fact that the constant $C>0$ is independent of $\tau$) and the assumption \eqref{gambetabnd} on the structure of $P(\cdot,\cdot)$ at once furnishes the equi-integrability of $\{P_{\delta}(\rho^{\tau},Z^{\tau})\}_{\tau}$ in $(0,T)\times B$. 
	Hence by using de la Vall\'{e}e-Poussin criterion we have the following
	\begin{equation}\label{afterusingdelavalle}
		P_{\delta}(\rho^{\tau},Z^{\tau})\rightharpoonup \overline{P_{\delta}(\rho,Z)}\quad\mbox{in}\quad L^{1}((0,T)\times B).
	\end{equation}
	The next goal is to identify the weak limit $\overline{P_{\delta}(\rho,Z)}$ with $P_{\delta}(\rho,Z)$. In that direction we first define in accordance with the convention from \eqref{conv} for $(t,x)\in[0,T]\times B$ and 
	\begin{equation}
		s^\tau(t,x)=\frac{Z^\tau(t,x)}{\rho^\tau(t,x)},\ s(t,x)=\frac{Z(t,x)}{\rho(t,x)}
	\end{equation}
	The application of Lemma \ref{Lem:AlmComp} on the sequence $\{\rho^\tau, Z^\tau, u^\tau\}$ (and the corresponding constant sequence of displacements) yields
	\begin{equation}\label{AlmostCompLim}
		\lim_{\tau\to 0_+}\int_B\rho^\tau(t,\cdot)|s^\tau(t,\cdot)-s(t,\cdot)|^p=0\text{ for all }t\in[0,T]\text{ and any }p\in[1,\infty).
	\end{equation}
	Next, we write
	$$\displaystyle P_{\delta}(\rho^{\tau},Z^{\tau})=P_{\delta}(\rho^{\tau},\rho^{\tau}s^{\tau})=P_{\delta}(\rho^{\tau},\rho^{\tau}s^{\tau})-P_{\delta}(\rho^{\tau},\rho^{\tau}s)+P_{\delta}(\rho^{\tau},\rho^{\tau}s).$$
	We claim that
	\begin{equation}\label{limpassfixs}
		\begin{split}
			\lim_{\tau\rightarrow 0_+} \int_{\mathcal Q} \left|(P_{\delta}(\rho^{\tau},\rho^{\tau}s^{\tau})-P_{\delta}(\rho^{\tau},\rho^{\tau}s))\right|=0\text{ for any }\mathcal Q\Subset[0,T]\times (\overline B\setminus\Sigma_\eta(t)).
		\end{split}
	\end{equation}
	We notice that due to \eqref{EtaTCnv} for fixed $\mathcal Q$ there is $\tau_0$ such that $\mathcal Q\Subset[0,T]\times (\overline B\setminus\Sigma_{\eta^\tau})$ for any $\tau<\tau_0$. Then applying \eqref{PressDerEst*}, the identity 
	\begin{equation*}
		a^r-b^r=(a-b)(a^{r-1}+a^{r-2}b+\ldots+ab^{r-2}+b^{r-1})\text{ for }a,b\geq 0,r\in\eN
	\end{equation*}
	it follows that
	\begin{equation}\label{estdiffpress}
		\int_\mathcal{Q}\left|(P_{\delta}(\rho^{\tau},\rho^{\tau}s^{\tau})-P_{\delta}(\rho^{\tau},\rho^{\tau}s))\right|\leq c(\delta)\left(\int_\mathcal{Q} ((\rho^\tau)^{-\underline{\kappa}+1}+(\rho^\tau)^{\overline{\kappa}})|s^\tau-s|+\int_{\mathcal Q}(\rho^\tau)^\kappa|s^\tau-s|\right)
	\end{equation}
	for any $\tau<\tau_0,$ where the first summand on the right hand side of the above estimate is obtained by using mean-value theorem and the assumption \eqref{PressDerEst*}.  
	Applying the H\"older inequality, uniform estimate \eqref{improvedpressurestawayint} and \eqref{AlmostCompLim} we conclude \eqref{limpassfixs}. Let us note that, the interface $(0,T)\times\Sigma_{\eta}$ is H\"{o}lder continuous ($cf.$ \eqref{EtaTCnv}) and hence for each member of $[0,T]\times(\overline{B}\setminus\Sigma_{\eta}(t))$ it is always possible to choose a parabolic neighborhood $\mathcal{Q}$ of the same such that $\mathcal{Q}\Subset [0,T]\times(\overline{B}\setminus\Sigma_{\eta}(t)).$ Hence we immediately infer from \eqref{limpassfixs}
	\begin{equation}\label{WLPrId}
		\overline{P_\delta(\rho,Z)}=\overline{p_\delta(\rho)}\text{ a.e. in }(0,T)\times B,
	\end{equation}
	where $p_\delta(r)=P_\delta(r,rs)$.
	
	$(4)$ \textit{Continuity of the fluid and structural velocities on the interface, choice of test functions in \eqref{momentumex}}
	We want to verify item $(i)$ of Definition~\ref{Def:EPSol} for $\eta$ and $u$ obtained in \eqref{TauConv}. In particular, we have $\Tr_{\Sigma_{\eta^\tau}}u^\tau-\tder\eta^\tau\nu\to 0$ as $\tau\to 0_+$ by \eqref{TauEst}$_1$. With regard to \eqref{TauConv}$_2$ it remains to show that $u^\tau\circ\tilde\varphi_{\eta^\tau}\rightharpoonup u\circ\tilde\varphi_\eta$ in $L^1((0,T)\times \eR^3)$ after extending $u^\tau$ by zero in $\eR^3\setminus B$. Since the same is proven in Section~\ref{contvelo} but for less regular flow maps, we refer for details therein.\\
	Further we notice that the test functions used at the approximate layer solve the compatibility $b\nu=tr_{\Sigma_{\eta^{\tau}}}\phi$ at the interface $(0,T)\times\Sigma_{\eta^{\tau}}$ (we refer to \eqref{AdmTest}). Using same pair of test functions both at the approximate level and at the limit (as $\tau\rightarrow 0$) might not guarantee the interface compatibility in limit. The way is to construct a test function for the limiting equation by suitable approximation. Rather than giving a details here we would like to refer the readers to Section \ref{constestfn} for such a construction (where it is done even with restricted regularities of the unknowns). We remark that the strong convergence of the sequence $\{\partial_{t}\eta^{\delta}\}$ is a crucial part of such a construction and the former can be proved by following the arguments used to show \eqref{TDerEtaStrong}.
	
	$(5)$ \textit{The limit passage $\tau\rightarrow 0_{+}$ in the equations:} Having necessary convergences we can perform the limit passage $\tau\to 0_+$ in \eqref{ConEqAppr} and \eqref{MomEqAppr} for $(\rho, Z, u,\eta, \phi)=(\rho^\tau, Z^\tau, u^\tau,\eta^\tau, \phi^\tau)$. Indeed, using \eqref{RZTStrongly} and \eqref{RUZUTWeakSt} we conclude \eqref{contrhoex}. In order to perform the limit passage in the momentum equation, we fix an arbitrary pair $(\phi,b)$ of admissible test functions in \eqref{momentumex}. Next, fixing an arbitrary $\psi\in C^\infty((0,T))$, multiplying \eqref{MomEqAppr} by $\psi(t)$, integrating the identity over $(0,T)$, employing \eqref{limpassageenenergyeta}, \eqref{RZUTCnv}, \eqref{CnvTPas}, \eqref{TauConv}$_{1,2,5}$, \eqref{convdeltapen}, \eqref{afterusingdelavalle} and \eqref{WLPrId} we conclude that
	\begin{equation}\label{weakformaftertau0}
		\begin{split}
			&\int_0^t\int_B \left((\rho+Z)\left(u\cdot\tder \phi + (u\otimes u)\cdot\nabla \phi\right)+\overline{p_\delta(\rho)}\dvr\phi -\mathbb S^\eta_\omega(\mathbb D u)\cdot\nabla\phi\right)+\int_0^t\int_\Gamma \tder\eta \tder b -\int_0^t\langle K'_\delta(\eta),b\rangle\\
			&-\zeta\int_{(0,T)\times\Gamma}\partial_{t}\nabla\eta\nabla b=\int_B(\rho+Z) u(t,\cdot)\phi(t,\cdot)-\int_BM_{0,\delta}\cdot\phi(0,\cdot)+\int_\Gamma\tder\eta(t,\cdot)b(t,\cdot)-\int_\Gamma\eta_1b(0,\cdot)
		\end{split}
	\end{equation}
	for a.a. $t\in(0,T)$ and all $(\phi,b)\in C^\infty_c([0,T]\times B)\times L^\infty(0,T;W^{3,2}(\Gamma))\cap W^{1,\infty}(0,T;L^2(\Gamma))$ such that $b\nu=\Tr_{\Sigma_\eta}\phi$. 
	
	The next task is to identify $\overline{{p}_{\delta}(\rho)}$ or equivalently $\overline{{P}_{\delta}(\rho,Z)}$ with $P_{\delta}(\rho,Z).$ Thanks to \eqref{WLPrId}, we can apply a strategy similar to the theory of mono-fluid with non-monotone pressure law developed in \cite{FeireislPet}. This will be done in the spirit of \cite{NovoPoko} adapted to our case in order to suitably handle the presence of a moving interface $\Sigma_{\eta}$ inside the fixed domain $B.$ We first state a local version of the effective viscous flux equality, which can be proved by using the arguments presented in \cite[Section 3.6.5]{FeiNovo}. The following equality in the context of FSI problems can also be found in \cite{Tri1} and \cite{Breit}. 
	\begin{lem}\label{effvisflux}
		Upto a non-explicitly relabeled subsequence of $\tau\rightarrow 0_{+}$ the following identity holds
		\begin{equation}\label{effvisfleq}
			\lim\limits_{\tau\rightarrow 0} \int_{(0,T)\times B} \phi\bigg(p_{\delta}(\rho^{\tau})\rho^{\tau}-(\lambda+2\mu)\rho^{\tau}\mathrm{div}\,u^{\tau}\bigg)=\int_{(0,T)\times B}\phi\bigg(\overline{p_{\delta}(\rho)}-(\lambda+2\mu)\mathrm{div}\,u\bigg)\rho
		\end{equation}
		for all $\phi\in C^{\infty}_{c}(((0,T)\times B)\setminus((0,T)\times\Sigma_{\eta} )).$
	\end{lem}
	We note that for any $\phi\in C^{\infty}_{c}(((0,T)\times B)\setminus((0,T)\times\Sigma_{\eta} ))$ we have $\phi\in C^{\infty}_{c}(((0,T)\times B)\setminus((0,T)\times\Sigma_{\eta^\tau} ))$ for any $\tau<\tau_0$ for $\tau_0$ small enough by \eqref{EtaTCnv}.
	From the arbitrariness of the test function $\phi,$ the convergence \eqref{effvisfleq} leads to the following equality which holds a.e in $(0,T)\times B:$\\
	\begin{equation}\label{effvisfluxae}
		\overline{p_{\delta}(\rho)\rho}-(\lambda+2\mu)\overline{\rho\dvr{u}}=\overline{p_{\delta}(\rho)}\rho-(\lambda+2\mu)\rho\dvr{u}\text{ a.e. in }(0,T)\times B,
	\end{equation}
	where $\overline{p_{\delta}(\rho)\rho}$ and $\overline{\rho\dvr\,u}$ denote respectively the $L^{1}$ weak limits of $p_{\delta}(\rho^{\tau})\rho^{\tau}$ and $\rho^{\tau}\dvr\,u^{\tau}$ respectively.\\[2.mm] 
	$(6)$ \textit{Strong convergence of $\{\rho^{\tau}\}$:}
	It is by now classical in the literature the effective viscous flux identity \eqref{effvisfluxae} relates the quantity $\overline{\rho\dvr{u}}-\rho\dvr{u}$ with the defect measure of the density oscillations described via the renormalized continuity equations. This will be clear from the following discussion. In order to prove the strong convergence of $\{\rho^{\tau}\}$ we will adapt the arguments used in \cite{NovoPoko} with suitable modifications. 
	
	Let us now show that 
	\begin{equation}\label{almstconvrhotau}
		\rho^{\tau}\rightarrow \rho\mbox{ a.e. in } (0,T)\times B.
	\end{equation}
	Since both the pairs $(\rho^{\tau},u^{\tau})$ and $(\rho,u)$ are solutions to continuity equations with $u^{\tau},u\in L^{2}(0,T;W^{1,2}(B))$ and $\rho^{\tau},\rho\in L^{\infty}(0,T;L^{\kappa}(B))$ they both solve renormalized continuity equations, cf. Lemma \ref{Lem:Renormalization} with a function $\mathcal B(r)=L_k(r)$. The truncation $L_k$ of the function $\rho\ln(\rho)$, is defined as
	\begin{equation}\label{Lkrho}
		L_{k}(\rho)=\rho\int^{\rho}_{1}\frac{T_{k}(z)}{z^{2}}\mbox{ for }k>1,
	\end{equation}
	and $T_k$ stands for an $L^\infty$--truncation defined for any $k>1$ via
	$$T_{k}(z)=kT\left(\frac{z}{k}\right)$$
	where 
	\begin{equation}\label{truction}
		T(z)=\begin{cases}
			z&\mbox{ for } z\in[0,1),\\
			\mbox{concave} & \mbox{ for } z\in[1,3),\\
			2&\mbox{ for } z\geqslant 3.
		\end{cases}
	\end{equation}
	By using the renormalized continuity equation we obtain (we have stated a version of renormalized continuity equation in Lemma \ref{Lem:Renormalization} and the one we are using now is a simpler version of that in a fixed domain) for $r\in\{\rho^\tau,\rho\}$ and the corresponding $v\in\{u^\tau,u\}$
	\begin{equation*}
		\int_B (L_k(r)\phi)(\cdot,t)-\int_B (L_k(r)\phi)(\cdot,0)=\int_0^t\int_B T_k(r)\dvr v\phi+L_k(r)(\tder\phi+v\cdot\nabla \phi)
	\end{equation*}
	for $t\in[0,T]$ and $\phi\in C^\infty([0,T]\times\eR^3)$.
	We perform the passage $k\to \infty$ employing the obvious convergences $T_k(r)\to r$ in $L^2((0,T)\times B)$ and $L_k(r)\to r\log r$ in $C_w([0,T];L^{\frac{\kappa}{2}}(B))$ and $L_k(r)(0)\to \rho_{0,\delta}\log(\rho_{0,\delta})$ in $L^1(B)$ and arrive at 
	\begin{equation}\label{RenEqRLR}
		\int_B (r\log r\phi)(\cdot,t)-\int_B (\rho_{0,\delta}\log(\rho_{0,\delta})\phi)(\cdot,0)=\int_0^t\int_B r\dvr v\phi+r\log r(\tder\phi+v\cdot\nabla \phi)
	\end{equation}
	for $t\in[0,T]$ and $\phi\in C^\infty([0,T]\times\eR^3)$. 
	Setting $(r,v)=(\rho^{\tau},u^{\tau})$ and $\phi=1$ in \eqref{RenEqRLR} and further passing to the limit $\tau\rightarrow 0_+$ we obtain
	\begin{equation}\label{transaftau0}
		\int_{B}(\overline{\rho\log(\rho)})(\cdot,t)-\int_{B}\rho_{0,\delta}\log(\rho_{0,\delta}))(0,\cdot)=\int^{t}_{0}\int_{B}\overline{\rho\dvr{u}}
	\end{equation}
	for $t\in[0,T].$ Now \eqref{transaftau0} and \eqref{RenEqRLR} with $(r,v,\phi)=(\rho, u,1)$ together furnish
	\begin{equation}\label{inequalityrhologrho}
		\int_{B}(\overline{\rho\log(\rho)}-\rho\log(\rho))(\cdot,t)=\int^{t}_{0}\int_{B}(\overline{\rho\dvr{u}}-{\rho\dvr{u}})
	\end{equation}
	for $t\in[0,T]$.\\
	Next we recall from \eqref{eq2.4} and \eqref{approxpressure} that we have the following decomposition of $p_{\delta}(\rho)$ at our disposal
	\begin{equation}\label{decomposepidelta}
		p_{\delta}(\rho)=P_{\delta}(\rho,\rho s)=\mathcal{P}(\rho,s)+\mathcal{M}_{\delta}(\rho,s)-\mathcal{R}(\rho,s)
	\end{equation}
	where $\rho\mapsto \mathcal{M}_{\delta}(\rho,s)$ is a monotone non-decreasing function. Now the monotonicity of the maps $\rho\mapsto \mathcal{P}(\rho,s)$ and $\rho\mapsto\mathcal{M}_{\delta}(\rho,s)$ render that
	\begin{equation}\label{inqfrmmono}
		\bigg(\overline{\mathcal{P}(\rho,s)\rho}-\overline{\mathcal{P}(\rho,s)}\rho\bigg)+\bigg(\overline{\mathcal{M}_{\delta}(\rho,s)\rho}-\overline{\mathcal{M}_{\delta}(\rho,s)}\rho\bigg)\geqslant 0,\text{ a.e. in }(0,T)\times B
	\end{equation}
	by Lemma \ref{lem:WLP}.
	In view of the a.e. effective viscous flux identity \eqref{effvisfluxae}, the inequality \eqref{inqfrmmono} and the decomposition \eqref{decomposepidelta} we obtain the following from \eqref{inequalityrhologrho}
	\begin{equation}\label{boundoscill1}
		\int_{B} \bigg(\overline{\rho\log(\rho)}-\rho\log(\rho)\bigg)(\cdot,t)\leqslant \frac{1}{(\lambda+2\mu)}\int^{t}_{0}\int_{B}\bigg(\overline{\mathcal{R}(\rho,s)\rho}-\overline{\mathcal{R}(\rho,s)}\rho\bigg)
	\end{equation}
	for $t\in[0,T]$.\\
	The idea next is to majorize the right hand side of \eqref{boundoscill1} by a constant multiple of $\displaystyle\int^{t}_{0}\int_{B}\bigg(\overline{\rho\log(\rho)}-\rho\log(\rho)\bigg)$ and further to use the Gronwall lemma to conclude that the expression on the left hand side of \eqref{boundoscill1} vanishes for $t\in[0,T]$. In that direction we follow the ideas developed in \cite[Sec. 4.3.]{NovoPoko} and present the details for the sake of completeness.\\
	Since for $s\in[\underline{a},\overline{a}],$ $\mathcal{R}(\rho,s)$ is uniformly bounded in $C^{2}([0,\infty))$ and compactly supported, there exists a possibly large constant $\Lambda>0$ such that both the functions $\rho\mapsto\Lambda \rho\log(\rho)-\rho\mathcal{R}(\rho,s)$ and $\rho\mapsto \Lambda\rho\log(\rho)+\mathcal{R}(\rho,s)$ are convex on $[0,\infty)$ for any $s\in[\underline{a},\overline{a}].$ First, the convexity of $\rho\mapsto \Lambda\rho\log(\rho)-\rho\mathcal{R}(\rho,s)$ furnishes the following by using \cite[Cor. 3.33, item (iii)]{NovStr04}) 
	\begin{equation}\label{diffrholnrholb}
		\overline{\rho\mathcal{R}(\rho,s)}-\rho\mathcal{R}(\rho,s)\leqslant \Lambda\bigg(\overline{\rho\log(\rho)}-\rho\log(\rho)\bigg)
	\end{equation} 
	and hence 
	\begin{equation}\label{1ststepbndrhsoscill}
		\int^{t}_{0}\int_{B}\bigg(\overline{\mathcal{R}(\rho,s)\rho}-\overline{\mathcal{R}(\rho,s)}\rho\bigg)\\
		\leqslant\Lambda\int^{t}_{0}\int_{B}\bigg(\overline{\rho\log(\rho)}-\rho\log(\rho)\bigg)+\int^{t}_{0}\int_{B}\bigg(\mathcal{R}(\rho,s)-\overline{\mathcal{R}(\rho,s)}\bigg)\rho.
	\end{equation}
	Further the convexity of $\rho\mapsto\Lambda\rho\log(\rho)+\mathcal{R}(\rho,s)$ renders
	\begin{equation}\label{convexity2nd}
		\mathcal{R}(\rho,s)-\overline{\mathcal{R}(\rho,s)}\leqslant \Lambda\bigg(\overline{\rho\log(\rho)}-\rho\log(\rho)\bigg).
	\end{equation}
	Since for each $s\in[\underline{a},\overline{a}],$ $\mathcal{R}(\rho,s)$ is supported in $[0,\overline{R}],$ so is $\overline{\mathcal{R}(\rho,s)}$ and hence as a consequence of \eqref{convexity2nd} one computes the following
	$$\int^{t}_{0}\int_{B}\bigg(\mathcal{R}(\rho,s)-\overline{\mathcal{R}(\rho,s)}\bigg)\rho\leqslant \Lambda\overline{R}\int^{t}_{0}\int_{B}\bigg(\overline{\rho\log(\rho)}-\rho\log(\rho)\bigg)$$
	which together with \eqref{1ststepbndrhsoscill} and \eqref{boundoscill1} furnishes
	\begin{equation}\label{penultimateGron}
		\int_{B} \bigg(\overline{\rho\log(\rho)}-\rho\log(\rho)\bigg)(\cdot,t)\leqslant \frac{\Lambda}{(\lambda+2\mu)}(1+\overline{R})\int^{t}_{0}\int_{B}\bigg(\overline{\rho\log(\rho)}-\rho\log(\rho)\bigg),
	\end{equation}
	for $t\in[0,T]$. As \eqref{transaftau0} and \eqref{RenEqRLR} with $(r,v)=(\rho,u)$ imply $\overline{\rho\log(\rho)}(0)=\rho_{0,\delta}\log(\rho_{0,\delta})=\rho\log(\rho)(0)$ a.e. in $B$, we infer from \eqref{penultimateGron} by the Gronwall inequality and by the strict convexity of $\rho\mapsto\rho\log(\rho)$ that $\overline{\rho\log(\rho)}=\rho\log\rho$ a.e. in $(0,T)\times B$. Hence we conclude \eqref{almstconvrhotau}.
	
	As the sequence $\{P_\delta(\rho^\tau,Z^\tau)\}$ is equiintegrable and $P_\delta(\rho^\tau, Z^\tau)\to P_\delta(\rho,Z)$ a.e. in $(0,T)\times B$, we conclude by the Vitali convergence theorem
	\begin{equation}\label{idnfcomplete}
		\overline{P_\delta(\rho,Z)}=\overline{p_\delta(\rho)}=p_\delta(\rho)=P_\delta(\rho,Z)
	\end{equation}
	and \eqref{momentumex} is verified by using \eqref{weakformaftertau0}.
	
	Let us note that we have for a nonrelabeled subsequence of $\{Z^\tau\}$ that
	\begin{equation}\label{ZTauConv}
		Z^\tau\to Z\text{ a.e. in }(0,T)\times B
	\end{equation}
	as an immediate consequence of convergence \eqref{almstconvrhotau} and Lemma~\ref{Lem:AlmComp}, identity \eqref{ZeroLimFr} respectively.
	
	$(7)$ \textit{The proof of \eqref{energybalanceex}}\\
	First, we note that by obvious algebraic manipulations and the definition of $\mathbb S^\eta_\omega(\mathbb D u)$ in \eqref{SEODef} we get
	\begin{equation}\label{SEONabU}
		\mathbb S^\eta_\omega(\mathbb D u)\cdot\nabla u=2\mu^\eta_\omega\left|\mathbb D u-\frac{1}{3}\dvr u\mathbb I\right|^2+\lambda^\eta_\omega|\dvr u|^2.
	\end{equation}
	Next, we fix an arbitrary nonnegative function  $\psi\in C^\infty_c((0,T))$, multiply \eqref{UBTau} for $(\rho,Z,u,\eta,v)=(\rho^\tau,Z^\tau,u^\tau,\eta^\tau,v^\tau)$ by $\psi(t)$ $t\in(0,T)$, integrate over $(0,T)$, neglect the terms containing $v$ and conclude \eqref{energybalanceex} by employing convergences \eqref{rhou2conv}, \eqref{CnvTPas}, \eqref{TauConv}$_{1,2,5}$,
	\begin{equation}\label{AuxConv}
		\begin{alignedat}{2}
			\mathcal H_{P,\delta}(\rho^\tau,Z^\tau)&\to     \mathcal H_{P,\delta}(\rho,Z)&&\text{ in }L^1((0,T)\times B),\\
			\sqrt{\mu^{\eta^\tau}_\omega}\left(\mathbb D u^\tau-\frac{1}{3}\dvr u^\tau\mathbb I\right)&\rightharpoonup\sqrt{\mu^{\eta}_\omega}\left(\mathbb D u-\frac{1}{3}\dvr u\mathbb I\right)&&\text{ in }L^2((0,T)\times B),\\
			\sqrt{\lambda^{\eta^\tau}_\omega}\dvr u^\tau&\to\sqrt{\lambda^{\eta}_\omega}\dvr u&&\text{ in }L^2((0,T)\times B)
		\end{alignedat}
	\end{equation}
	and the weak lower semicontinuity of $L^2$--norm taking into account \eqref{SEONabU}. We observe that convergence \eqref{AuxConv}$_1$ follows by the Vitali convergence theorem. Pointwise convergences \eqref{rhou2conv} and \eqref{ZTauConv} imply $\mathcal H_{P,\delta}(\rho^\tau,Z^\tau)\to\mathcal H_{P,\delta}(\rho,Z)$. The continuity of $(\rho,Z)\mapsto \mathcal H_{P,\delta}(\rho,Z)$, defined in \eqref{HPdelta}, on $[0,\infty)^2$ follows from \eqref{HelmFDef} taking into consideration \eqref{PowerEst}. The equintegrability of $\{\mathcal H_{P,\delta}(\rho^\tau,Z^\tau)\}$ is a consequence of the definition of $\mathcal H_{P,\delta}$, the growth of $P$ in \eqref{gambetabnd} and the estimate in \eqref{improvedpressurestawayint}. Finally, convergences \eqref{AuxConv}$_{2,3}$ follow by \eqref{AuxApprCnv}, \eqref{ViscAppr} and \eqref{TauConv}$_5$.
	\subsection{Summary of the proof of Theorem \ref{resultstaulayer}:} In the process of construction of solution we notice that $\rho^{n+1},Z^{n+1}\geq 0$ (cf. item 1. of Section \ref{Sec:FlSub}) and hence by definition of interpolants \eqref{piecewdef} and the weak convergences \eqref{TauConv}$_{3,4}$ one concludes $\rho, Z\geq 0$ in $(0,T)\times B$. Estimate \eqref{TauEst}$_4$ and convergence \eqref{rhou2conv} imply $(\rho+Z)|u|^2\in L^\infty(0,T;L^1(B))$. The other regularities of the unknowns listed in the Definition \ref{WSExtProb} are consequences of \eqref{TauConv}. That $P_{\delta}(\rho,Z)\in L^{1}((0,T)\times B)$ follows from \eqref{afterusingdelavalle} and \eqref{idnfcomplete}.\\
	The momentum balance \eqref{momentumex} is recovered as a consequence of \eqref{weakformaftertau0} and \eqref{idnfcomplete}. In the limiting set-up, we need the continuity of the fluid and the structural velocities at the interface and further the test functions to solve $b\nu=tr_{\Sigma^{\eta}}\phi$ at the interface $(0,T)\times \Sigma_{\eta}.$ We discuss these two properties on item $(4)$ appearing just after \eqref{WLPrId}.\\  Further the recovery of the mass balance \eqref{contrhoex} from \eqref{ConEqAppr} is relatively simple and is discussed after \eqref{WLPrId}.\\
	For the proof of the energy inequality \eqref{energybalanceex}-\eqref{hdelta}, we refer to the discussion from \eqref{SEONabU} and afterwards.\\
	The existence of a positive $T$ such that the weak solution to the extended problem exists on $(0,T)$ can be shown by repeating arguments from Section~\ref{minintex} leading to \eqref{PointBounds}. Moreover, one can apply the extension procedure from Section~\ref{Maxtimeex} to obtain the maximal interval of existence for the weak solution.
	
	%%%%%%%%%%%%%%%%%%%%%%%%%%%%%%%%%%%%%%%%%%%%%%%%%%%%%%%%%%%%%%%%%%%%%%%%%%%%%%%%%%%%%%%%%%%%%%%%%%%%%%%%%%%%%%%%%%%%%%%%%%%%%%%%%%%%%%%%%%%%%%%%%%%%%%%%%%%%%%%%%%%%%%%%%%%%%%%%%%%%%%%%%%%%%%%%%%%%%%%%%%%%%%%%%%%%%%%%%%%%%%%%%%%%%%%%%%%%%%%%%%%%%%%%%%%%%%%%%%%%%%%%%%%%%%%%%%%%%%%%%%%%%%%%%%%%%%%%%%%%%%%%%%%%%%%%%%%%%%%%%%%%%%%%%%%%%%%%%%%%%%%%%%%%%%%%%%%%%%%%%%%%%%%%%%%%%%%%%%%%%%%%%%%%%%%%%%%%%%%%%%%%%%%%%%%%%%%%
	
	\section{Limit as \texorpdfstring{$\omega,\zeta,\delta\to 0_+$}{}}\label{reglim0}
	The goal of this section is to prove Theorem \ref{Thm:main}. We recall that Theorem \ref{Thm:main} has two parts, $\textit{Case I}$ dealing with the existence issue when the adiabatic exponents solve $\max\{\gamma,\beta\}>2$ and $\min\{\gamma,\beta\}>0$ whereas the $\textit{Case II}$ is associated with $\max\{\gamma,\beta\}\geq2$ and $\min\{\gamma,\beta\}>0.$ Indeed for the proof of $\textit{Case I}$ we will rely on estimates independent of the dissipation parameter $\zeta$ whereas for $\textit{Case II}$ we will use the structural dissipation ($\zeta>0$). First we present in details the proof of $\textit{Case I}$ and next we will comment on the proof of $\textit{Case II}.$
	\subsection{{{Proof\,\, of\,\, Case I:}}}\label{ProofcaseI}
	For the proof of $\textit{Case I}$, we set 
	\begin{equation}\label{delomgep0}
		\omega=\zeta=\delta
	\end{equation}
	in \eqref{momentumex}-\eqref{energybalanceex} and perform the limit passage $\delta\to 0_+$. By the end of this section we will be able to conclude the proof of Theorem \ref{Thm:main}. We begin with a collection of estimates which are uniform in $\delta.$
	\subsubsection{Uniform in \texorpdfstring{$\delta$}{} estimates and weak convergences}
	In this section we collect estimates satisfied by $(\eta^{\delta},\rho^{\delta},Z^{\delta},u^{\delta})$ which are independent of $\delta$ and are obtained as a consequence of \eqref{energybalanceex}. Some of these estimates are listed below:
	
	\begin{equation}\label{UnifEstPhDom}
		\begin{split}
			\|{\eta^{\delta}}\|_{L^\infty(0,T;W^{2,2}(\Gamma))}\leq& c,\\
			\|\tder{\eta^{\delta}}\|_{L^\infty(0,T;L^{2}(\Gamma))}+\delta^\frac{1}{2}\|\nabla\tder{\eta^{\delta}}\|_{L^2((0,T)\times\Gamma)}+\delta^\frac{7}{2}\|\nabla^3{\eta^{\delta}}\|_{L^{\infty}(0,T;L^{2}(\Gamma)}\leq& c,\\
			\|{\rho^{\delta}}\|_{L^\infty(0,T;L^\gamma(B))}+\|{Z^{\delta}}\|_{L^\infty(0,T;L^\beta(B))}\leq& c,\\
			\delta^\frac{1}{\kappa}\left(\|{\rho^{\delta}}\|_{L^\infty(0,T;L^\kappa(B))}+\|{Z^{\delta}}\|_{L^\infty(0,T;L^\kappa(B))}\right)\leq& c,\\
			\|\sqrt{{\rho^{\delta}}+{Z^{\delta}}}u^{\delta}\|_{L^\infty(0,T;L^2(B))}\leq& c,\\
			\int_{(0,T)\times B}\mathbb S^{\eta^{\delta}}_{\delta}(\mathbb D u^{\delta})\cdot\nabla u^{\delta}\leq& c,\\
			(\rho^{\delta}(t,x),Z^{\delta}(t,x))\in \overline{\mathcal{O}_{\underline{a}}},\,\,\mbox{for all}\,\, t\in(0,T)\,\,\mbox{and almost all}\,\,x\in B,
		\end{split}
	\end{equation}
	where the last inclusion follows from \eqref{DensComparT}.\\ Estimate \eqref{UnifEstPhDom}$_6$ does not provide directly a uniform bound with respect to $\delta$ on $\nabla u^{\delta}$. In fact, we deduce from \eqref{UnifEstPhDom}$_6$ and \eqref{ExtendingFProp}$_{2}$ that
	\begin{equation}\label{GradDeltaEst}
		\begin{split}
			&\int_{Q^T_{\eta^{\delta}}}\left(2\mu\bigg|\mathbb D u^{\delta}-\frac{1}{3}\dvr u\mathbb I\bigg|^2+\lambda|\dvr u^{\delta}|^2\right)\\&=\int_{Q^T_{\eta^{\delta}}}\left(2\mu\left(\mathbb D u^{\delta}{2}-\frac{1}{3}\dvr u^{\delta}\mathbb I\right)+\lambda\dvr u^{\delta}\mathbb I\right)\cdot\nabla u^{\delta}\leq \int_{(0,T)\times B}\mathbb S^{\eta^{\delta}}_{\delta}(\mathbb D u^{\delta})\cdot\nabla u^{\delta}\leq c
		\end{split}
	\end{equation}
	implying immediately
	\begin{equation}\label{UnifEstPhDom2}
		\begin{split}
			\|\dvr u^\delta\|_{L^2(Q^T_{\eta^\delta})}&\leq c,\\
			\|\mathbb{D}u^{\delta}\|_{L^{2}(Q^{T}_{\eta^{\delta}})}&\leq c
		\end{split}
	\end{equation}
	due to \eqref{CoefAss}, where the notation $\mathbb{D}$ stands for the symmetric gradient.
	Employing the Korn inequality on H\"older domains, i.e., Lemma \ref{Lem:Korn}, we infer
	\begin{equation}
		\begin{split}
			\|u^{\delta}\|^2_{L^2(0,T;W^{1,q}(\Omega_{\eta^{\delta}}(t)))}\leq& c\left(\|\mathbb D u^{\delta}\|^2_{L^2(0,T;L^2(\Omega_{\eta^{\delta}}(t)))}+\int_0^T\int_{\Omega_{\eta^{\delta}}(t)}({\rho^{\delta}}+{Z^{\delta}})|u^{\delta}|^2 \right)\\ \leq& c_1\int_{(0,T)\times B}\mathbb{S}^{\eta^{\delta}}_\delta(\mathbb D u^{\delta})\cdot\nabla u^{\delta}+c_2\|\sqrt{{\rho^{\delta}}+{Z^{\delta}}}u^{\delta}\|^2_{L^\infty(0,T;L^2(B))}\leq c
		\end{split}
	\end{equation}
	with the constant $c_1$ depending on $q$, initial data and $\mu$ and the constant $c_2$ depending on $q$, initial data and $T$.
	We note that the assumptions of Lemma \ref{Lem:Korn} are satisfied. In particular, we can take 
	$$
	L=2\max\{\|{\rho^{\delta}}\|_{L^\infty(0,T;L^{\gamma}(B))},\|{Z^{\delta}}\|_{L^\infty(0,T;L^{\beta}(B))}\}$$ and
	\begin{equation*}
		M=\inf_{\delta}\int_{\Omega_{{\eta}^\delta_0}}({\rho}_{0,\delta}+{Z}_{0,\delta})>0.
	\end{equation*}
	As ${\eta^{\delta}}$ and ${\rho^{\delta}}, {Z^{\delta}}$ satisfy \eqref{UnifEstPhDom}$_{1,3}$ and the formulation of the continuity equations in \eqref{contrhoex} and \eqref{VanDens} imply the conservation of mass, we get 
	\begin{equation*}
		\int_{\Omega_{\eta^{\delta}}(t)}({\rho^{\delta}}+{Z^{\delta}})(t,\cdot)=\int_{\Omega_{{\eta^{\delta}_{0}}}}({\rho}_{0,\delta}+{Z}_{0,\delta})\geq M\text{ for }t\in[0,T].
	\end{equation*} 
	
	Moreover, applying Lemma \ref{Lem:Extension} we get 
	\begin{equation}\label{UestOnB}
		\|u^{\delta}\|_{L^2(0,T;W^{1,q}(\eR^3))}\leq c\text{ for any }q\in[1,2).    
	\end{equation}
	%%%%%%%%%%%%%%%%%%%%%%%%%%%%%%%%%%%%%%%%%%%%%%%%%%%%%%%%%%%%%%%%%%%%%%%%%%%%%%%%%%%%%%%%%%%%%%%%%%%%%%%%%%%%%%%%%%%%%%%%%%%%%%%%%%%%%%%%%%%%%%%%%%%%%%%%%%%%%%%%%%%%%%%%%%%%%%%%%%%%%%%%%%%%%%%%%%%%%%%%%%%%%%%%%%%%%%%%%%%%%%%%%%%%%%%%%%%%%%%%%%%%%%%%%%%%%%%%%%%%%%%%%%%%%%%%%%%%%%%%%%%%%%%%%%%%%%%%%%%%%%%%%%%%%%%%%%%%%%%%%%%%%%%%%%%%%%%%%%%%%%%%%%%%%%%%%%%%%%%%%%%%%%%%%%%
	As a consequence of uniform estimates \eqref{UnifEstPhDom}$_{1,2,3}$ and \eqref{UestOnB} we get the existence of sequence $\{(\eta^\delta, \rho^\delta, Z^\delta, u^\delta)\}$  of solutions in the sense of Definition \ref{WSExtProb} such that 
	\begin{equation}\label{weaklimiturhoeta}
		\begin{alignedat}{2}
			\eta^\delta&\rightharpoonup^*\eta&&\text{ in }L^\infty(0,T;W^{2,2}(\Gamma)),\\
			\tder\eta^\delta&\rightharpoonup^*\tder\eta&&\text{ in }L^\infty(0,T;L^2(\Gamma)),\\
			\rho^\delta&\rightharpoonup^*\rho&&\text{ in }L^\infty(0,T;L^{\max\{\gamma,\beta\}}(B)),\\
			Z^\delta&\rightharpoonup^*Z &&\text{ in }L^\infty(0,T;L^{\max\{\gamma,\beta\}}(B)),\\
			u^\delta&\rightharpoonup u &&\text{ in }L^2(0,T;W^{1,q}(\eR^3))\text{ for any }q\in[1,2).
		\end{alignedat}
	\end{equation}
	Moreover, we conclude
	\begin{equation}\label{EtaDHCnv}
		\eta^\delta\to\eta\text{ in }C^{\frac{1}{4}}([0,T]\times\Gamma)
	\end{equation}
	and  
	\begin{equation}\label{DensVanish}
		\rho|_{B\setminus\Omega_\eta}(t)=Z|_{B\setminus\Omega_\eta}(t)\equiv 0
	\end{equation}
	since $\rho^\delta|_{B\setminus\Omega_{\eta^\delta}}(t)=Z^\delta|_{B\setminus\Omega_{\eta^\delta}}(t)\equiv 0$ as in subsection \ref{Sec:TauLim}. 
	
	\subsubsection{Continuity of the fluid and structural velocities at the interface}\label{contvelo} Next we verify that the limit pair $\tder\eta$ and $u$ satisfies the coupling condition $\tder\eta\nu=\Tr_{\Sigma_\eta}u$. Since $\{u^\delta\circ\tilde\varphi_{\eta^\delta}\}$ is bounded in $L^2(0,T;W^{1,q}(\eR^3))$ for any $q\in[1,2)$ we conclude the existence of a non-relabeled sub-sequence such that
	\begin{equation}\label{CompCnv}
		u^\delta\circ\tilde\varphi_{\eta^\delta}\rightharpoonup w\text{ in }L^2(0,T;W^{1,q}(\eR^3)).
	\end{equation}
	On the other hand convergence \eqref{weaklimiturhoeta}$_2$ implies that $\tder\eta\nu$ coincides with the trace of $w$ on $\partial\Omega$.
	Our task now is to identify the limit $w$. We fix an arbitrary $\zeta\in C^\infty_c((0,T)\times\eR^3)$ and write
	\begin{equation}\label{InIden}
		\int_{(0,T)\times\eR^3}(u^\delta\circ\tilde\varphi_{\eta^\delta}-u\circ \tilde\varphi_\eta)\cdot\zeta=\int_{(0,T)\times\eR^3}\left[(u^\delta-u)\circ\tilde\varphi_{\eta^\delta}\right]\cdot\zeta+\int_{(0,T)\times\eR^3}\left(u\circ \tilde\varphi_{\eta^\delta}-u\circ\tilde\varphi_\eta\right)\cdot\zeta.
	\end{equation}
	First, using the change of variables we obtain
	\begin{equation}
		\int_{(0,T)\times\eR^3}(u^\delta-u)\circ\tilde\varphi_{\eta^\delta}\cdot\zeta=\int_{(0,T)\times\eR^3}(u^\delta-u)\cdot\zeta\circ(\tilde\varphi_{\eta^\delta})^{-1}|\det\nabla(\tilde\varphi_{\eta^\delta})^{-1}|.
	\end{equation}
	We observe that $(\tilde\varphi_{\eta^\delta})^{-1}$ converges to $(\tilde\varphi_\eta)^{-1}$ locally uniformly in $[0,T]\times\eR^3$ as $(\tilde\varphi_\eta)^{-1}$ is obviously uniformly continuous on any compact in $[0,T]\times\eR^3$ and the homeomorphisms $\tilde\varphi_{\eta^\delta}$ converge to $\tilde\varphi_\eta$ locally uniformly in $[0,T]\times\eR^3$, which follows from \eqref{EtaDHCnv}. Hence $\zeta\circ(\tilde\varphi_{\eta^\delta})^{-1}$ converges to $\zeta\circ(\tilde\varphi_\eta)^{-1}$ uniformly in $[0,T]\times\eR^3$ as all the functions $\zeta\circ(\tilde\varphi_{\eta^\delta})^{-1}$ and $\zeta\circ(\tilde\varphi_\eta)^{-1}$ possess their supports in a compact subset of $[0,T]\times\eR^3$. Knowing also that $\nabla(\tilde \varphi_{\eta^\delta})^{-1}$ converges to $\nabla(\tilde \varphi_{\eta})^{-1}$ in $L^s((0,T)\times\eR^3)$ for any $s\in[1,\infty)$  we conclude using \eqref{weaklimiturhoeta}$_5$
	\begin{equation}\label{FirLimZ}
		\lim_{\delta\to 0}\int_{(0,T)\times\eR^3}(u^\delta-u)\circ\tilde\varphi_{\eta^\delta}\cdot\zeta=0.
	\end{equation}
	Let us focus on the second term on the right hand side of \eqref{InIden}. We get
	\begin{equation*}
		\int_{(0,T)\times\eR^3}(u\circ\tilde\varphi_{\eta^\delta}-u\circ\tilde\varphi_{\eta})\cdot\zeta=\int_{\{\tilde\varphi_{\eta^\delta}\neq\tilde\varphi_\eta\}}(u\circ\tilde\varphi_{\eta^\delta}-u\circ\tilde\varphi_\eta)\cdot\zeta.
	\end{equation*}
	Next, we deduce
	\begin{equation}\label{AIneqComp}
		\left|\int_{\{\tilde\varphi_{\eta^\delta}\neq\tilde\varphi_\eta\}}(u\circ\tilde\varphi_{\eta^\delta}-u\circ\tilde\varphi_\eta)\cdot\zeta\right|\leq |\{\tilde\varphi_{\eta^\delta}\neq\tilde\varphi_\eta\}|^\frac{1}{r'}(\|u\circ\tilde\varphi_{\eta^\delta}\|_{L^r(0,T;L^r(\eR^3))}+\|u\circ\tilde\varphi_{\eta}\|_{L^r(0,T;L^r(\eR^3))})\|\zeta\|_{L^\infty(0,T)\times\eR^3}
	\end{equation}
	for some $r\in(1,2)$.
	Using the change of variables and the fact that $\{\nabla(\tilde\varphi_{\eta^\delta})^{-1}\}$ is bounded in $L^s((0,T)\times\eR^3)$ for any $s\in[1,\infty)$ following from its convergence in the latter space we get
	$$
	\|u\circ\tilde\varphi_{\eta^\delta}\|^r_{L^r(0,T;L^r(\eR^3))}=\int_{(0,T)\times\eR^3}|u|^r|\det\nabla(\tilde\varphi_{\eta^\delta})^{-1}|\leq \|u\|_{L^2(0,T;L^2(\eR^3))}^{r}\|\det\nabla(\tilde\varphi_{\eta^\delta})^{-1}\|_{L^{\frac{2}{2-r}}((0,T)\times\eR^3)}.$$
	Observing that the right hand side of the latter inequality is bounded and  $|\{\tilde\varphi_{\eta^\delta}\neq\tilde\varphi_\eta\}|\to 0_+$ as $\delta\to 0_+$, which follows from convergence \eqref{EtaDHCnv}, we get from \eqref{AIneqComp}
	\begin{equation}\label{SecLimZ}
		\lim_{\delta\to 0_+}\int_{(0,T)\times\eR^3}(u\circ\tilde\varphi_{\eta^\delta}-u\circ\tilde\varphi_{\eta})\cdot\zeta=0.
	\end{equation}
	Having \eqref{FirLimZ} and \eqref{SecLimZ} at hand we infer from \eqref{InIden} that $w=u\circ\tilde\varphi_\eta$ in \eqref{CompCnv}. Accordingly, we have showed the coupling $\tder\eta\nu=\Tr_{\Sigma_\eta}u$.\\[2.mm]
	\subsubsection{Compactness of the shell-energy}\label{compactnessshelenergy}
	The next result is a very important one stating that $\eta^{\delta}$ (uniformly in $\delta$) enjoys better regularity in space. On one hand this extra regularity helps in the limit passage $\delta\rightarrow 0_{+}$ in the term $\int^{t}_{0}\langle K'_{\delta}(\eta^{\delta}),b\rangle$ (concerning the shell energy) and on the other hand it renders the fact that the fluid boundary is not just H\"{o}lder but is Lipschitz in space (to be precise in $L^{2}(C^{0,1})$).\\ 
	In the context of an incompressible fluid-structure interaction problem the ingenious idea of improving the structural regularity first appeared in \cite{MuhaSch} and was later adapted for compressible fluid-structure interaction problems in \cite{Breit2}. One can adapt the proof presented in \cite{Breit2} to the present scenario but only after  guaranteeing the use of some particular test functions using density argument. Hence we will state the lemma and detail the required density argument to adapt the test functions used in \cite{Breit2}.
	\begin{lem}\label{improvebndetadelta}
		Let the quadruple $(\rho^{\delta},Z^{\delta},u^{\delta},\eta^{\delta})$ solves an extended problem in the sense of Definition \ref{WSExtProb} (which is obtained as a consequence of Theorem \ref{resultstaulayer}). Then the following holds uniformly in $\delta:$
		\begin{equation}\label{improvedregeta}
			\begin{array}{l}
				\displaystyle \|\eta^{\delta}\|_{L^{2}(0,T;W^{2+\sigma^{*},2}(\Gamma))}\leqslant c,\,\,\|\partial_{t}\eta^{\delta}\|_{L^{2}(0,T;W^{\sigma^{*},2}(\Gamma))}\leqslant c,
			\end{array}
		\end{equation}
		for some $0<\sigma^{*}<\frac{1}{2},$ where the constant $c$ in the last inequality may depend on $\Gamma,$ the initial data and the $W^{2,2}$ coercivity size of $\overline{\gamma}(\eta)$ ($\overline{\gamma}(\eta)$ has been introduced in \eqref{ovgamma}).
	\end{lem}  
	{\textit{Comments on adapting the arguments from \cite[Section 5.2]{Breit2}} in order to prove Lemma \ref{improvebndetadelta}:}  Roughly the idea of the proof of Lemma \ref{improvebndetadelta}, is to test the structure by a correction of $\Delta^{\sigma}_{-h}\Delta^{\sigma}_{h}\eta^{\delta}$ (where $\Delta^{\sigma}_{h}=h^{-\sigma}\bigg(\eta(y+h e_{\alpha})-\eta(y)\bigg)$ is the fractional difference quotient in the direction of $e_{\alpha},$ $\alpha\in\{1,2\}$ ) and the fluid by a solenoidal lifting of the same. More precisely, we use in \eqref{momentumex} test functions of the form $(\phi^{\delta},b^{\delta})=(\mathcal{F}^{\textit{div}}_{\eta^{\delta}}(\Delta^{\sigma}_{-h}\Delta^{\sigma}_{h}\eta^{\delta}-\mathcal{K}_{\eta^{\delta}}(\Delta^{\sigma}_{-h}\Delta^{\sigma}_{h}\eta^{\delta})),\Delta^{\sigma}_{-h}\Delta^{\sigma}_{h}\eta^{\delta}-\mathcal{K}_{\eta^{\delta}}(\Delta^{\sigma}_{-h}\Delta^{\sigma}_{h}\eta^{\delta}))$, where $\mathcal{F}^{\textit{div}}$ and $\mathcal{K}_{\eta}$ were introduced in Proposition \ref{smestdivfrex}. Since so far we have used test functions which possess the regularity $C^{\infty}([0,T]\times\RR^{3})\times (L^{\infty}(0,T;W^{3,2}(\Gamma))\cap W^{1,\infty}(0,T;L^{2}(\Gamma)),$ we need to justify the use of test functions $(\phi^{\delta},b^{\delta})$ in \eqref{momentumex}, since their regularity is restricted.\\
	We further note that with the use of this test function, one needs to take care of the term 
	\begin{equation}\label{extradissterm}
		\zeta\int_{(0,t)\times\Gamma}\partial_{t}\nabla\eta^{\delta}\cdot\Delta^{\sigma}_{-h}\Delta^{\sigma}_{h}\nabla\eta^{\delta}=-\frac{\zeta}{2}\int_{(0,t)\times\Gamma}\partial_{t}|\Delta^{\sigma}_{h}\nabla\eta^{\delta}|^{2}
	\end{equation}
	originating from the structural dissipation (and this does not appear in \cite{Breit2}). We observe that the right hand side of \eqref{extradissterm} can be bounded by $c\zeta\left(\|\Delta^{\sigma}_{h}\nabla\eta^{\delta}\|_{L^{\infty}(0,T;L^{2}(\Gamma))}+\|\nabla^{2}\eta_{0}\|_{L^{2}(\Gamma)}\right),$ for some $c$ independent of $\delta,$ which in turn is bounded by $c\zeta\left(\|\eta^{\delta}\|_{L^{\infty}(0,T;W^{2,2}(\Gamma))}+\|\nabla^{2}\eta_{0}\|_{L^{2}(\Gamma)}\right)$ and hence by a constant $c\zeta$ by virtue of \eqref{UnifEstPhDom}$_{1}.$ Hence the dissipation term does not make any considerable difference in the calculation.   
	
	{\textit{A density argument justifying the use of test functions of the form $(\phi^{\delta},b^{\delta})$ in \eqref{momentumex}}:} One first observes that $b^{\delta}\in W_{1}=L^{\infty}(0,T;W^{3,2}(\Gamma))\cap W^{1,\infty}(0,T;L^{2}(\Gamma))$ (of course this is not true uniformly in $\delta$). This regularity along with the estimates \eqref{Fdivetaest1} and \eqref{higherderest1} from Proposition \ref{smestdivfrex} guarantees that $\phi^{\delta}\in W_{2}=L^{\infty}(0,T;W^{3,2}(B))\cup W^{1,\infty}(0,T; L^{2}(B))$ (where the notation $B$ to denote a neighborhood of $\Omega$ was introduced in \eqref{defB}) and $\phi^{\delta}$ is solenoidal. Indeed for small $\delta>0,$ $\Omega_{\eta^{\delta}}\Subset \Omega\cup S_{m,M}.$ Now we need a density argument to show that $(\phi^{\delta},b^{\delta})$ is an admissible test function in \eqref{momentumex}. In that direction by a standard argument involving convolution with mollifiers, we construct a sequence $\phi_{M}^{\delta}\in C^{\infty}([0,T]\times \mathbb{R}^{3})$ such that $\phi_{M}^{\delta}\rightarrow \phi^{\delta}$ in the weak$^{*}$ topology of $W_{2}.$ Next let us set $b^{\delta}_{M}\nu=\Tr_{\Sigma_{\eta^\delta}}\phi^{\delta}_{M},$ where the notion of trace $\Tr_{\Sigma_{\eta}}$ was introduced in Lemma \ref{Lem:TrOp}. Now the definition of $b^{\delta}_{M}$ along the uniform in $M$ bound of $\phi^{\delta}_{M}$ in $W_{2}$ and the bound of $\eta^{\delta}$ in $W_{1}$ together furnish that $b^{\delta}_{M}$ converges to $b^{\delta}=\Tr_{\Sigma_{\eta^\delta}}\phi^{\delta}\cdot\nu$ in the weak$^{*}$ topology of $W_{1}.$ Finally the weak$^{*}$ convergence of $(\phi^{\delta}_{M},b^{\delta}_{M})$ to $(\phi^{\delta},b^{\delta})$ in $W_{2}\times W_{1}$ (indeed one observes that the continuity of the trace $tr_{\eta^{\delta}}\phi^{\delta}_{M}=b^{\delta}_{M}\nu$ holds by construction) allows to verify that $(\phi^{\delta},b^{\delta})$ is an admissible pair of test function in \eqref{momentumex}.\\[3.mm]
	As a consequence of \eqref{improvedregeta} and the classical Aubin-Lions lemma we have the following strong convergence of $\eta^{\delta}:$\\
	\begin{equation}\label{strngconvetadel}
		\eta^{\delta}\rightarrow \eta \,\,\mbox{in}\,\, L^{2}(0,T;W^{2+,2}(\Gamma)).
	\end{equation}
	The convergence \eqref{strngconvetadel} will be used in particular for the limit passage in the term related to the Koiter energy in the weak-formulation of the momentum equation. The notation $2+$ signifies a number greater than $2$.
	
	\boxed{\textbf{Limit passage in the non-linearities and recovering a momentum balance:}}
	
	\subsubsection{ Some further convergences and the weak limit of the pressure} Employing \eqref{UestOnB} and the Sobolev embedding theorem we infer
	\begin{equation}\label{bndudelta}
		\|u^\delta\|_{L^2(0,T;L^p(B))}\leq c, \text{ for any }p\in[1,6).
	\end{equation}
	Combining the latter bound with \eqref{UnifEstPhDom} we obtain that
	\begin{equation*}
		\|(\rho^\delta+Z^\delta)u^\delta\otimes u^\delta\|_{L^2(0,T;L^r(B))}\leq c
	\end{equation*}
	for certain $r>1$. Employing the arguments based on the Bogovskii operator, we deduce
	\begin{equation}\label{BetterIntRho}
		\int_{\mathcal Q}\left((\rho^{\delta})^{\gamma+\theta_1}+(Z^{\delta})^{\beta+\theta_2}+\delta((\rho^{\delta})^{\kappa+\theta_1}+(Z^{\delta})^{\beta+\theta_2})\right)\leq c(Q)
	\end{equation}
	for some $\theta_1,\theta_2>0$ and any $\mathcal Q\Subset [0,T]\times (\overline B\setminus \Sigma_{\eta^\delta})$. Repeating arguments from Section \ref{tauto0eq}, we conclude the following convergences 
	\begin{equation}\label{AddConvDel}
		\begin{alignedat}{2}
			\eta^\delta&\to\eta&&\text{ in }C^\frac{1}{4}([0,T]\times\Gamma),\\
			(\rho^\delta, Z^\delta)&\to (\rho, Z)&&\text{ in }C_w([0,T]; L^{\max\{\gamma,\beta\}}(B)),\\
			(\rho^\delta+Z^\delta) u^\delta&\to(\rho+Z)u&&\text{ in }C_w([0,T];L^\frac{2\max\{\gamma,\beta\}}{\max\{\gamma,\beta\}+1}(B)),\\
			(\rho^\delta+Z^\delta) u^\delta\otimes u^\delta&\rightharpoonup(\rho+Z)u\otimes u&&\text{ in }L^1((0,T)\times B),\\
			(\rho^\delta+Z^\delta) |u^\delta|^2&\rightharpoonup(\rho+Z)|u|^2&&\text{ in }L^1((0,T)\times B).    
		\end{alignedat}
	\end{equation} Let us point out that when showing the uniform continuity of $\{(\rho^\delta+Z^\delta)u^{\delta}\}$ in $C([0,T];W^{-3,2}(B))$ one uses the pointwise inequality 
	\begin{equation}\label{BoundModS}
		|\mathbb S^{\eta^\delta}_\delta(\mathbb D u^\delta)|^2\leq \max\{\mu_\delta,\lambda_\delta\}\mathbb S^{\eta^\delta}_\delta(\mathbb D u^\delta)\cdot\nabla u^\delta
	\end{equation}  
	$\mu_\delta\leq\mu$, $\lambda_\delta\leq\lambda$ and \eqref{UnifEstPhDom}$_5$.
	
	%%%%%%%%%%%%%%%%%%%%%%%%%%%%%%%%%%%%
	
	Next we state a result on the strong convergence of $\partial_{t}\eta^{\delta},$ which is going to play a crucial role later in the proof of \eqref{convintetatest}$_{1}.$ At a Galerkin level the proof of such a convergence is detailed in \cite{Breit2} but the authors only use bounds uniform with respect to all the parameters and hence we can adapt arguments from \cite{Breit2} without much difficulties. This is the reason we do not provide a full proof of the following lemma but we sketch it in the appendix.
	\begin{lem}\label{strngconvdteta}
		Let the assertions of Lemma \ref{improvebndetadelta} hold. Then
		\begin{equation}\label{TDerEtaStrong}
			\tder\eta^\delta\to\tder\eta\text{ in }L^2(0,T;L^2(\Gamma)).
		\end{equation}
	\end{lem}
	We will comment on the proof of Lemma \ref{strngconvdteta} in Section \ref{lempfstrnged}.
	
	%%%%%%%%%%%%%%%%%%%%%%%%%%%%%%%%%%%%%%%%%%%%%%%%%%%%%%%%%%%%%%%%%%%%%%%%%%%%%%%%%%%%%%%%%%%%%%%%%%%%%%%%%%%%%%%%%%%%%%%%%%%%%%%%%%%%%%%%%%%%%%%%%%%%%%%%%%%%%%%%%%%%%%%%%%%%%%%%%%%%%%%%%%%%%%%%%%%%%%%%%%%%%%%%%%%%%%%%%%%%%%%%%%%%%%%%%%%%%%%%%%%%%%%%%%%%%%%%%%%%%%%%%%%%%%%%%%%%%%%%%%%%%%%%%%%%%%%%%%%%%%%%%%%%%%%%%%%%%%%%%%%%%%%%%%%%%%%%%%%%%%%%%%%%%
	
	Next the estimate \eqref{BetterIntRho} implies that 
	\begin{equation}\label{L1convpress}
		P_\delta(\rho^\delta,Z^\delta)\rightharpoonup \overline{P(\rho,Z)} \text{ in }L^1(\mathcal Q)
	\end{equation}
	for any $\mathcal Q\Subset (\overline B\setminus\Sigma_\eta(t))\times [0,T]$. We consider a sequence of compact sets $\{K_i\}$ such that $K_i\subset K_{i+1}$, $K_i\cap \left([0,T]\times\Sigma_\eta(t)\right)=\emptyset$ and $|([0,T]\times\overline B)\setminus K_i|\to 0$. In order to exclude the possible concentration of $\{P_\delta(\rho^\delta,Z^\delta)\}$ at the moving boundary, we employ the next lemma which concerns the equi-integrability of $\{P_{\delta}(\rho^{\delta},Z^{\delta})\}$ close to the non-Lipschitz hyper-surface $\Sigma_{\eta}.$ 
	\begin{lem}\label{pressureestnearinterface}
		For any $\varepsilon>0,$ there exists a $\delta_{0}>0$ and $\mathcal{A}_{\varepsilon}\Subset B\times (0,T)$ such that for all $\delta<\delta_{0}$ the following holds
		\begin{equation}\label{equiintnearinterface}
			\mathcal{A}_{\varepsilon}\cap (\Sigma_{\eta^{\delta}}\times[0,T])=\emptyset,\qquad \int_{((0,T)\times B)\setminus \mathcal{A}_{\varepsilon}}P_{\delta}(\rho^{\delta},Z^{\delta})\leqslant \varepsilon.
		\end{equation} 
	\end{lem}
	In the context of non-Lipschitz domains (without the structure) a result of the form \eqref{equiintnearinterface} was first proved in \cite[Lemma 8]{Kuku} and later adapted for fluid-structure interaction problems in \cite[Lemma 6.4]{Breit} and \cite[Lemma 3.4]{MaMuNeRoTr}. The proof of Lemma \ref{pressureestnearinterface} can be done by imitating the arguments of \cite{Breit} since the proof does not depend on the structure of the pressure. \\
	%%%%%%%%%%%%%%%%%%%%%%%%%%%%%%%%%%%%%%%%%%%%%%%%%%%%%%%%%%%%%%%
	Employing Lemma \ref{pressureestnearinterface} and estimate \eqref{BetterIntRho} we conclude the equiintegrability of $\{P_\delta(\rho^\delta, Z^\delta)\}$ and
	\begin{equation}\label{PiDL1W}
		P_\delta(\rho^\delta,Z^\delta)\rightharpoonup\overline{P(\rho, Z)} \text{ in }L^1((0,T)\times B)
	\end{equation}
	accordingly. Moreover, fixing an arbitrary compact subset $\mathcal K$ of  $[0,T]\times (\overline B\setminus \Sigma_{\eta})$ we infer
	\begin{equation}\label{ArtPressVanish}
		\delta(\rho^\kappa+Z^\kappa+\frac{1}{2}\rho^{\kappa-2}Z^2+\frac{1}{2}\rho^2Z^{\kappa-2})\to 0\text{ in }L^{\frac{\kappa+\theta}{\kappa}}(\mathcal K)
	\end{equation}
	by \eqref{BetterIntRho}. Rewriting next $P(\rho^\delta, Z^\delta)$ by \textbf{H4}, namely by \eqref{eq2.4}, we obtain
	\begin{equation}\label{PressDec}
		P(\rho^\delta, Z^\delta)=P(\rho^\delta,\rho^\delta s^\delta)-P(\rho^\delta,\rho^\delta s)+\mathcal P(\rho^\delta,\rho^\delta s)+\mathcal R(\rho^\delta,\rho^\delta s),
	\end{equation}
	where in agreement with \eqref{conv}
	\begin{equation*}
		s^\delta=\frac{Z^\delta}{\rho^\delta},\ s=\frac{Z}{\rho}.
	\end{equation*}
	Now we will apply almost compactness argument ($i.e.$ Lemma \ref{Lem:AlmComp}) in order to freeze one of the densities in the expression of the pressure. One verifies the assertions of Lemma \ref{Lem:AlmComp}. In particular we verify \eqref{BddAss} and \eqref{valqgb2}. Note that for the case $\max\{\gamma,\beta\}>2$ and $\min\{\gamma,\beta\}>0$ since we have $\eta^{\delta}$ bounded in $L^{\infty}(0,T;W^{2,2}(\Gamma))$ (cf. \eqref{UnifEstPhDom}$_{1}$) Lemma \ref{Lem:AlmComp} applies.  Hence repeating arguments leading to \eqref{limpassfixs} we conclude
	\begin{equation}\label{estdiffpressdel}
		\lim_{\delta\to 0_+}\|P(\rho^\delta,\rho^\delta s^\delta)-P(\rho^\delta,\rho^\delta s)\|_{L^1(\mathcal K)}=0.
	\end{equation}
	Hence combining \eqref{PiDL1W}, \eqref{ArtPressVanish} and \eqref{PressDec} we obtain
	\begin{equation}\label{PressLimId}
		\overline{P(\rho, Z)}=\overline {p(\rho)}\text{ a.e. in }(0,T)\times B,
	\end{equation}
	where $p(r)=\mathcal P(r,rs)+\mathcal R(r,rs)$. 
	
	\subsubsection{Construction of test functions for the momentum balance equation}\label{constestfn}
	We now look for test functions which solve the compatibility condition $b\nu=\Tr_{\Sigma_{\eta}}\phi$ at the limiting interface $\Sigma_{\eta}.$ Indeed the test functions used in the approximate layer solve similar compatibility $b^{\delta}\nu=\Tr_{\Sigma_{\eta^{\delta}}}\phi^{\delta}$ on $\Sigma_{\eta^{\delta}}$ but they might not solve the same on $\Sigma_{\eta}.$\\
	We start by fixing a function $\phi\in C^{\infty}([0,T]\times \mathbb{R}^{3})$ and next define $b^{\delta}$ as follows:
	\begin{equation}\label{constestdeltalev}
		b^{\delta}=\Tr_{\Sigma_{\eta^{\delta}}}\phi\cdot\nu,
	\end{equation}
	where the notion of trace $\Tr_{\Sigma_{\eta}}$ was introduced in Lemma \ref{Lem:TrOp}.\\
	Now in view of the uniform in $\delta$ bounds of $\eta^{\delta}$ (we refer to \eqref{UnifEstPhDom}$_{1,2},$ \eqref{improvedregeta} and \eqref{strngconvetadel}) one has the following convergences of $b^{\delta}$
	\begin{equation}\nonumber
		\begin{alignedat}{2}
			b^{\delta}&\rightharpoonup^{*} b&&\mbox{ in } L^{\infty}(0,T;W^{2,2}(\Gamma))\cap W^{1,\infty}(0,T;L^{2}(\Gamma)),\\
			b^{\delta}&\rightharpoonup b&&\mbox{ in }
			L^{2}(0,T;W^{2+\sigma,2}(\Gamma))\mbox{ for some }\sigma>0,\\
			b^{\delta}&\rightarrow b&&\mbox{ in } L^{2}(0,T;W^{2+,2}(\Gamma)),\\
			b&=\Tr_{\Sigma_{\eta}}&&\phi\cdot \nu.
		\end{alignedat}
	\end{equation}
	The listed convergences above along with \eqref{TDerEtaStrong} is enough to conclude the following
	\begin{equation}\label{convintetatest}
		\begin{array}{ll}
			&\displaystyle \int^{t}_{0}\int_{\Gamma}\partial_{t}\eta^{\delta}\partial_{t}b^{\delta}\rightarrow \int^{t}_{0}\int_{\Gamma}\partial_{t}\eta\partial_{t}b,\qquad\displaystyle \int^{t}_{0}\langle K'_{\delta}(\eta^{\delta}),b^{\delta}\rangle \rightarrow \int^{t}_{0}\langle K'(\eta),b\rangle.
		\end{array}
	\end{equation}
	Specifically for the second convergence in \eqref{convintetatest}, one first verifies that $K'_{\delta}(\eta^{\delta})\rightharpoonup K'(\eta)$ in $L^{2}(0,T;L^{1+}(\Gamma)),$ in view of the convergences \eqref{weaklimiturhoeta}$_{1}$, \eqref{strngconvetadel} and the structure of $K'_{\delta}$ (cf. \eqref{elasticityoperator},\eqref{amab},\eqref{Geta},\eqref{Gij'},\eqref{exaR}, \eqref{regshellenergy} and \eqref{UnifEstPhDom}$_{2}$).\\
	Until now we have all the necessary convergences for the limit passage in the momentum equation. We just need the following in order to guarantee the weak formulation of the momentum equation in the physical domain $\Omega_{\eta}:$
	\begin{equation}\label{LimDelS}
		\lim_{\delta\to 0_+}\int_0^t\int_{B\setminus \Omega_{\eta^\delta}(s)}\mathbb S^{\eta^\delta}_\delta(\mathbb D u^\delta)\cdot\nabla\phi=0,
	\end{equation}
	for $\phi\in C^{\infty}((0,T)\times \mathbb{R}^{3}).$\\
	To this end we estimate
	\begin{equation*}
		\begin{split}
			\left| \int_0^t\int_{B\setminus \Omega_{\eta^\delta}(s)}\mathbb S^{\eta^\delta}_\delta(\mathbb D u^\delta)\cdot\nabla\phi\right|&\leq c \|\nabla\phi\|_{L^\infty((0,T)\times B)}\|\mathbb S^{\eta^\delta}_\delta(\mathbb D u^\delta)\|_{L^1(((0,T)\times B)\setminus Q^T_{\eta^\omega})}\\
			&\leq c\left(\int_{(0,T)\times B}\mathbb S^{\eta^\delta}_\delta(\mathbb D u^\delta)\cdot\nabla u^\delta\right)^\frac{1}{2}\left(\max\{\mu,\lambda\}\|f^{\eta^\delta}_\delta\|_{L^1(((0,T)\times B)\setminus Q^T_{\eta^\delta}))}\right)^\frac{1}{2},
		\end{split}
	\end{equation*}
	further using \eqref{BoundModS}, employing \eqref{UnifEstPhDom}$_5$ and \eqref{ExtendingFProp}$_3$ we conclude \eqref{LimDelS}.\\
	Now considering \eqref{momentumex} with $(\rho,Z,u,\eta)=(\rho^\delta,Z^\delta,u^\delta,\eta^\delta)$, a corresponding sequence of admissible test functions $\{\phi, b^{\delta}\}$ (as constructed in \eqref{constestdeltalev}) with a fixed $\phi$ possessing the regularity $C^{\infty}([0,T]\times \mathbb{R}^{3}),$ employing the convergences \eqref{AddConvDel}$_{3,4},$ \eqref{weaklimiturhoeta}$_{5},$ \eqref{convintetatest}, \eqref{PiDL1W}--\eqref{PressLimId} and further using \eqref{DensVanish}, \eqref{LimDelS} and the convergence of $M_{0,\delta}$ from \eqref{extensionindata} we conclude
	\begin{equation}\label{LimDeltaIdent}
		\begin{split}
			&\int_0^t\int_{\Omega_\eta(t)} (\rho+Z)\left(u\cdot\tder \phi + (u\otimes u)\cdot\nabla \phi\right)+\overline{p(\rho)}\dvr\phi -\mathbb S^\eta_\omega(\mathbb D u)\cdot\nabla\phi+\int_0^t\int_\Gamma \tder\eta \tder b -\int_0^t\langle K'(\eta),b\rangle\\
			&-\int_{\Omega_{\eta}(t)}(\rho+Z) u(t,\cdot)\phi(t,\cdot)+\int_{\Omega_{\eta_0}}M_0\cdot\phi(0,\cdot)+\int_\Gamma\eta_1b(0,\cdot)\\&=\int_\Gamma\tder\eta(t,\cdot)b(t,\cdot)
		\end{split}
	\end{equation}
	for almost all $t\in(0,T)$ and all $(\phi,b)\in C^\infty(\overline{Q^T_\eta})\times L^\infty(0,T;W^{2,2}(\Gamma))\cap W^{1,\infty}(0,T;L^2(\Gamma))\cap L^{2}(0,T;W^{2+\sigma,2})$ for some $\sigma>0,$ such that $b\nu=\Tr_{\Sigma_\eta}\phi.$ We observe that the left hand side of the latter identity is defined for any $t\in[0,T]$.
	
	Next we plan to show that the identity \eqref{LimDeltaIdent} holds for any $t\in[0,T].$ In that direction we first prove that
	\begin{equation}\label{TDerEtaContW}
		\tder\eta\in C_w([0,T];L^2(\Gamma)).
	\end{equation}
	To this end we consider $b\in C^\infty(\Gamma)$ and corresponding $\phi\in C^\infty(\overline{Q^T_\eta})$ such that $b\nu=\Tr_{\Sigma_\eta}\phi$. We know that there is a set $N\subset [0,T]$, $|N|=0$ such that \eqref{LimDeltaIdent} holds for any $t\in[0,T]\setminus N$. Fixing $t\in N$ we find a sequence $\{t^n\}\subset [0,T]\setminus N$ such that $t^n\to t$. Evidently, as the left hand side of \eqref{LimDeltaIdent} is defined for $t$ taking into account that $(\rho+Z)u\in C_w([0,T];L^\frac{2\max{\gamma,\beta}}{\max{\gamma,\beta}+1}(\eR^3))$, there exists 
	\begin{equation*}
		\lim_{t^n\to t}\int_{\Gamma}\tder\eta(t^n,\cdot)b\in\eR.
	\end{equation*}
	Furthemore, we conclude that the mapping
	\begin{equation}
		C^\infty(\Gamma)\ni b\mapsto\lim_{t^n\to t}\int_\Gamma \tder \eta(t^n,\cdot)b
	\end{equation}
	is linear and thanks to the regularity $\tder\eta\in L^\infty(0,T;L^2(\Gamma))$ also bounded. Hence there is $g(t)\in L^2(\Gamma)$ such that 
	\begin{equation*}
		\int_\Gamma g(t)b=\lim_{t^n\to t}\int_\Gamma\tder\eta(t,\cdot)b
	\end{equation*}
	by the Riesz representation theorem. Since $g\in C_w([0,T];L^2(\Gamma))$, which follows from \eqref{LimDeltaIdent}, we can define $\tder\eta(t)=g(t)$ for $t\in N$ to conclude \eqref{TDerEtaContW}. Now we have \eqref{LimDeltaIdent} meaningful for any $t\in[0,T]$.\\[2.mm]
	
	\boxed{\textbf{Strong convergence of \texorpdfstring{$\rho^{\delta}$}{}:}}\\ The next task is to prove strong convergence/ a.e. convergence of the density sequence $\rho^{\delta}.$ This is the key step to identify the weak limit $\overline{p(\rho)}$ with $p(\rho).$ For convenience of the reader we divide the proof into two sections, namely \ref{effvsflxid} and \ref{conc} respectively.
	\subsubsection{An effective-viscous flux identity}\label{effvsflxid} Here we state an effective viscous-flux equality which is a little different compared to the one stated in Lemma \ref{effvisflux}. 
	In view of \eqref{UnifEstPhDom}, we have the following estimates at our disposal
	\begin{equation}\label{estTkLk}
		\begin{split}
			\|T_{k}(\rho^{\delta})\|_{L^{r}(Q^{T}_{\eta^{\delta}})}&\leqslant c,\mbox{ for all } 1\leqslant r\leqslant\max\{\gamma,\beta\},\\
			\|L_{k}(\rho^{\delta})\|_{L^{\infty}(0,T;L^{r}(\Omega_{\eta^{\delta}}))}&\leqslant c, \mbox{ for all } 1\leqslant r <\max\{\gamma,\beta\},
		\end{split}
	\end{equation}
	where $L_{k}$ and $T_{k}$ were introduced in \eqref{Lkrho}-\eqref{truction}.\\
	Now let us state the effective-viscous flux identity in form of the following lemma whose proof can be done following exactly the line of arguments used in showing \cite[(7.23)]{Breit}:
	\begin{lem}\label{effviscflux2}
		Up to a non-explicitly relabeled subsequence of $\delta\rightarrow 0_{+}$ the following identity holds
		\begin{equation}\label{effvisfleq2}
			\int_{Q^T_{\eta^\delta}}\bigg(P(\rho^{\delta},\rho^{\delta}s)-(\lambda+2\mu)\dvr\,u^{\delta}\bigg)T_{k}(\rho^{\delta})\rightarrow \int_{Q^T_\eta}\bigg(\overline{p(\rho)}-(\lambda+2\mu)\dvr\,u\bigg)\overline{T_{k}(\rho)}.
		\end{equation}
	\end{lem}

	\subsubsection{The conclusion}\label{conc} 
	We will use our result on the renormalized continuity equations presented in form of Lemma \ref{Lem:Renormalization}. One notices that both the couples $(\rho^{\delta},u^{\delta})$ and $(\rho,u)$ solve the assertions of Lemma \ref{Lem:Renormalization} and hence they solve the renormalized continuity equations up to the interface. We will always consider extensions of $\rho^{\delta},$ $u^{\delta},$ $\rho$ and $u$ in entire $\RR^{3},$ more specifically $\rho^{\delta}$ and $\rho$ are defined zero outside $\Omega_{\eta^{\delta}}$ and $\Omega_{\eta}$ respectively. Hence the functions $L_{k}(\rho^{\delta})$ and $L_{k}(\rho)$ are also zero outside $\Omega_{\eta^{\delta}}$ and $\Omega_{\eta}.$  Applying Lemma \ref{Lem:Renormalization} with test function $\phi=1$ and the relation $r\nabla_{r}L_{k}(r)-L_{k}(r)=T_{k}(r),$ one infers the following
	\begin{equation}\label{diffLkr}
		\int_{\RR^{3}}\bigg(L_{k}(\rho^{\delta})-L_{k}(\rho)\bigg)(\cdot,t)=\int^{t}_{0}\int_{\RR^{3}}(T_{k}(\rho)\dvr u-\overline{T_{k}(\rho)}\dvr u^{\delta})+\int^{t}_{0}\int_{\RR^{3}}\bigg(\overline{T_{k}(\rho)}-T_{k}(\rho^{\delta})\bigg)\dvr u^{\delta}
	\end{equation}
	for $t\in [0,T],$ where $\overline{T_{k}(\rho)}$ is the weak limit of $T_{k}(\rho^{\delta})$ in $L^{r}((0,T)\times\eR^3)$ for any $1<r<\infty.$ Now employing \eqref{effvisfleq2}, the decomposition \eqref{eq2.4} and further passing $\delta\rightarrow 0$ we furnish the following
	\begin{equation}\label{calafterdiffLKr}
		\begin{split}
			\int_{\RR^{3}}\bigg(\overline{L_{k}(\rho)}-L_{k}(\rho)\bigg)(\cdot,t)&=\int^{t}_{0}\int_{\eR^3}\bigg(T_{k}(\rho)-\overline{T_{k}(\rho)}\bigg)\dvr u\\
			&+ \frac{1}{(\lambda+2\mu)}\int^{t}_{0}\int_{\RR^{3}}\bigg(\overline{\mathcal{P}(\rho,s)}\,\overline{T_{k}(\rho)}-\overline{\mathcal{P}(\rho,s)T_{k}(\rho)}\bigg)\\
			&+ \frac{1}{(\lambda+2\mu)}\int^{t}_{0}\int_{\RR^{3}}\bigg(\overline{\mathcal{R}(\rho,s)}\,\overline{T_{k}(\rho)}-\overline{\mathcal{R}(\rho,s)T_{k}(\rho)}\bigg)=\sum^{3}_{i=1}I^{\rho}_{i}
		\end{split}
	\end{equation}
	for all $t\in[0,T],$ where $\overline{L_{k}(\rho)}$ is the $C_{w}([0,T];L^{r}(\eR^3))$ limit of $\{L_{k}(\rho^{\delta})\}.$ Now in the case $\max\{\gamma,\beta\}>2,$ $\min\{\gamma,\beta\}>0,$ one can choose a $q=q(\max\{\gamma,\beta\})<2$ such the conjugate exponent $q^{*}=q^{*}(\max\{\gamma,\beta\})\in (2,\max\{\gamma,\beta\})$ and $I^{\rho}_{1}$ is estimated as follows using \eqref{UestOnB} and \eqref{estTkLk}$_1$
	\begin{equation}\label{I1E}
		\begin{array}{ll}
			|I_{1}^{\rho}| &\displaystyle \leqslant C\|\dvr\,u\|_{L^{q}((0,T)\times\mathbb{R}^{3})}\|T_{k}(\rho)-\overline{T_{k}(\rho)}\|_{L^{q*}((0,T)\times\mathbb{R}^{3})}\\[2.mm]
			&\displaystyle \leqslant C\limsup_{\delta\rightarrow 0}\|T_{k}(\rho^{\delta})-{T_{k}(\rho)}\|_{L^{q*}((0,T)\times\mathbb{R}^{3})}\\
			&\displaystyle \displaystyle \leqslant C\limsup_{\delta\rightarrow 0}\bigg(\|T_{k}(\rho^{\delta})-{T_{k}(\rho)}\|^{\frac{\max\{\gamma,\beta\}-q^{*}}{q*(\max\{\gamma,\beta\}-1)}}_{L^{1}((0,T)\times\mathbb{R}^{3})}\|T_{k}(\rho^{\delta})-{T_{k}(\rho)}\|^{\frac{\max\{\gamma,\beta\}(q*-1)}{q*(\max\{\gamma,\beta\}-1)}}_{L^{\max\{\gamma,\beta\}}((0,T)\times\mathbb{R}^{3})}\bigg)\\
			&\displaystyle \leqslant C\limsup_{\delta\rightarrow 0} \|T_{k}(\rho^{\delta})-{T_{k}(\rho)}\|^{\frac{\max\{\gamma,\beta\}-q^{*}}{q*(\max\{\gamma,\beta\}-1)}}_{L^{1}((0,T)\times\mathbb{R}^{3})}.
		\end{array}
	\end{equation}
	The term $I^\rho_2$ in \eqref{calafterdiffLKr} is non-positive by virtue of Lemma \ref{lem:WLP}. Following the similar line of arguments as used in \eqref{1ststepbndrhsoscill}-\eqref{convexity2nd}, leads to bound 
	\begin{equation}\label{I3RB}
		|I^{\rho}_{3}|\leqslant\Lambda(1+\overline{R})\int^{t}_{0}\int_{\mathbb{R}^{3}}\left(\overline{\rho\log(\rho)}-\rho\log(\rho)\right),
	\end{equation}
	for sufficiently large $\Lambda>1$ and $\overline{R}$ is as it appears in \eqref{bndbR}.\\ 
	Further since 
	$$\overline{L_{k}(\rho)}\rightarrow \overline{\rho\log(\rho)},\,\,L_{k}(\rho)\rightarrow \rho\log(\rho)\mbox{ in } C_{w}([0,T];L^{r}(\mathbb{R}^{3}))\mbox{ for any } 1\leqslant r<\max\{\gamma,\beta\}$$
	and 
	$$\|T_{k}(\rho)-\overline{T_{k}(\rho)}\|_{L^{1}((0,T)\times\mathbb{R}^{3})}\leqslant \|T_{k}(\rho)-\rho\|_{L^{1}((0,T)\times\mathbb{R}^{3})}+\liminf_{\delta\rightarrow 0}\|T_{k}(\rho_{\delta})-\rho_{\delta}\|_{L^{1}((0,T)\times\mathbb{R}^{3})}\rightarrow 0\,\, \mbox{as}\,\, k\rightarrow\infty$$ we obtain from 
	\eqref{calafterdiffLKr} by \eqref{I1E} and \eqref{I3RB} that
	\begin{equation*}
		\begin{split}
			\int_{\RR^{3}}\bigg(\overline{\rho\log(\rho)}-\rho\log(\rho)\bigg)(\cdot,t)&\leqslant\frac{\Lambda(1+\overline{R})}{(\lambda+2\mu)}\int^{t}_{0}\int_{\mathbb{R}^{3}}\left(\overline{\rho\log(\rho)}-\rho\log(\rho)\right)
		\end{split}
	\end{equation*}
	for all $t\in[0,T]$. The latter inequality is used to conclude \begin{equation}\label{RhoDPoint}
		\rho^\delta\to\rho\text{ a.e. in }(0,T)\times\eR^3
	\end{equation}
	by repeating the arguments based on Gronwall lemma from Section \ref{tauto0eq} (more specifically we refer the readers to the discussion leading to \eqref{idnfcomplete} from \eqref{penultimateGron}).
	%%%%%%%%%%%%%%%%%%%%%%%%%%%%%%%%%%%%%%%%%%%%%%%%%%%%%%%%%%%%%%%%%%%%%%%%%%%%%%%%%%%%%%%%%%%%%%%%%%%%%%%%%%%%%%%%%%%%%%%%%%%%%%%%%%%%%%%%%%%%%%%%%%%%%%%%%%%%%%%%%%%%%%%%%%%%%%%%%%%%%%%%%%%%%%%%%%%%%%%%%%%%%%%%%%%%%%%%%%%%%%%%%%%%%%%%%%%%%%%%%%%%%%%%%%%%%%%%%%%%%%%%%%%%%%%%%%%%%%%%%%%%%%%%%%%%%%%%%%%%%%%%%%%%%%%%

	%%%%%%%%%%%%%%%%%%%%%%%%%%%%%%%%%%%%%%%%%%%%%%%%%%%%%%%%%%%%%%%%%%%%%%%%%%%%%%%%%%%%%%%%%%%%%%%%%%%%%%%%%%%%%%%%%%%%%%%%%%%%%%%%%%%%%%%%%%%%%%%%%%%%%%%%%%%%%%%%%%%%%%%%%%%%%%%%%%%%%%%%%%%%%%%%%%%%%%%%%%%%%%%%%%%%%%%%%%%%%%%%%%%%%%%%%%%%%%%%
	
	\subsubsection{Fulfillment of energy inequality}\label{finalenergybalance}
	Obtaining \eqref{energybalance} from \eqref{energybalanceex}
	is very similar to the proof of \eqref{energybalanceex} in Section \ref{Sec:TauLim}. Hence we comment on differences only. 
	First, we observe that convergence \eqref{AddConvDel}$_1$ implies that 
	\begin{equation}\label{ChFunCon}
		\chi_{Q^t_{\eta^\delta}}\to\chi_{Q^t_{\eta}}\text{ in }L^p(\eR^4)\text{ for any }p\in[1,\infty)\text{ and }t\in(0,T).
	\end{equation}
	Due to pointwise convergence of $\{\rho^\delta\}$ in \eqref{RhoDPoint} implying $Z^\delta\to Z$ a.e. in $(0,T)\times B$ and the continuity of $(\rho,Z)\mapsto H_P(\rho,Z)$ we get $\mathcal{H}_{P,\delta}(\rho^\delta,Z^\delta)\to H_P(\rho,Z)$ a.e. in $(0,T)\times B$ with $H_P$ and $\mathcal H_{P,\delta}$ defined in \eqref{HelmFDef}, \eqref{HPdelta} respectively. Moreover, the growth estimate in \eqref{gambetabnd}, \eqref{UnifEstPhDom}$_7$, Lemma \ref{pressureestnearinterface} and estimate \eqref{BetterIntRho} yield the equiintegrability of the sequence $\{\mathcal H_{P,\delta}(\rho^\delta,Z^\delta)\}$. Hence we conclude
	\begin{equation*}
		\mathcal H_{P,\delta}(\rho^\delta,Z^\delta)\to H_{P}(\rho,Z)\text{ in }L^1((0,T)\times B).
	\end{equation*}
	Next, convergence \eqref{ChFunCon} combined with convergence \eqref{weaklimiturhoeta}$_5$ and the estimate \eqref{GradDeltaEst} implies
	\begin{equation*}
		\begin{alignedat}{2}
			\sqrt{\chi_{Q^T_{\eta^\delta}}}\left(\mathbb Du^\delta-\frac{1}{3}\dvr u^\delta\mathbb I\right)&\rightharpoonup \sqrt{\chi_{Q^T_{\eta}}}\left(\mathbb D u-\frac{1}{3}\dvr u\mathbb I\right)&&\text{ in }L^2(\eR^4),\\
			\sqrt{\chi_{Q^T_{\eta^\delta}}}\dvr u^\delta&\rightharpoonup \sqrt{\chi_{Q^T_{\eta}}}\dvr u\text{ in }L^2(\eR^4).
		\end{alignedat}
	\end{equation*}
	Taking into consideration the latter convergences and applying the weak lower semicontinuity of $L^2$--norm we infer
	\begin{equation*}
		\begin{split}
			\liminf_{\delta\to 0_+}
			\int_{Q^t_{\eta^\delta}}\mathbb S^{\eta^\delta}_{\delta}(\mathbb D u^\delta)\cdot\nabla u^\delta=&\liminf_{\delta\to 0_+}\int_{\eR^4}\chi_{Q^t_{\eta^\delta}}\left(2\mu\left|\mathbb Du^\delta-\frac{1}{3}\dvr u^\delta\mathbb I\right|^2+\lambda|\dvr u^\delta|^2\right)\\ \geq& \int_{\eR^4}\chi_{Q^T_{\eta}}\left(2\mu\left|\mathbb Du-\frac{1}{3}\dvr u\mathbb I\right|^2+\lambda|\dvr u|^2\right)=\int_{Q^t_{\eta}}\mathbb S(\mathbb D u)\cdot\nabla u
		\end{split}
	\end{equation*}
	for $t\in(0,T)$. 
	Finally, to obtain the first, second and fourth term on the right hand side of \eqref{energybalance} we use \eqref{extensionindata}$_{2,3}$. 
	\subsubsection{Attainment of initial data}\label{attinitialdata}
	From \eqref{AddConvDel}$_3$ we have 
	\begin{equation}\label{RZUContW}
		(\rho+Z)u\in C_w([0,T];L^\frac{2\max\{\gamma,\beta\}}{\max\{\gamma,\beta\}+1}(\Omega_{\eta}(t))) 
	\end{equation}
	considering \eqref{DensVanish}. Taking an arbitrary $\phi\in C^\infty(\eR^3)$ such that the support of $\phi\circ\tilde\varphi_\eta$ is compact in $[0,T]\times\Omega$, i.e., $\Tr_{\Sigma_\eta}\phi=0$, in \eqref{momentum} we perform the passage $t\to 0_+$ to deduce $(\rho+Z)u(0)=M_0$ a.e. in $\Omega_{\eta_0}$.  Furthermore, employing \eqref{RZUContW} and \eqref{TDerEtaContW} we get from \eqref{momentum} that $\tder \eta(0)=\eta_1$ a.e. in $\Gamma$, which concludes \eqref{InDatAtt}$_2$. The proof of the equalities in \eqref{InDatAtt}$_1$ follows the proof of the first equality in \eqref{InDatAtt}$_2$.
	\subsubsection{A minimal time where the solution avoids degeneracy}\label{minintex}
	In this section we show that it is possible to avoid both kind of degeneracy (\eqref{degenfstkind} and \eqref{degensndkind}) for a positive minimal time.\\ 
	We begin with specifying a minimal time $T$ independent of the regularizing parameter $\delta$ for which the $W^{2,2}$--coercivity of the displacement is available. The reason why this minimal time need not coincide with an initially given time is that the Koiter energy can achieve degeneracy depending on the sign of $\overline{\gamma}(\eta)$ (the quantity introduced in \eqref{ovgamma}). The minimal time is found with the help of the following lemma dealing with the $W^{2,2}$--coercivity of the non-linear Koiter shell.\\
	The following lemma further ensures that the degeracy of the second kind (we refer to \eqref{degensndkind}) can be avoided for some minimal time. 
	
	\begin{lem}\label{Lem:Coer}
		Let the assumptions of Theorem \ref{Thm:main} be satisfied and
		\begin{equation}\label{abndonelenergy}
			K_{\delta}(\eta)(\cdot,t)\leqslant \bigg(\int_{B}\bigg(\frac{|M_{0,\delta}|^{2}}{2(\rho_{0,\delta}+Z_{0,\delta})}+{\mathcal{H}}_{P,\delta}(\rho_{0,\delta},Z_{0,\delta})\bigg)+\bigg(\frac{1}{2}\int_{\Gamma}|\eta_{1}|^{2}+K_{\delta}(\eta_{0})\bigg)\bigg)=C_{0}
		\end{equation}
		hold, where $K_{\delta}$, $M_{0,\delta},$ $\rho_{0,\delta}$, $Z_{0,\delta}$, $\eta^{\delta}_{0}$ and ${\mathcal{H}}_{P,\delta}$ are introduced in Section \ref{Sec:ExtProb}. Then if $\overline{\gamma}(\eta)\neq 0$ we have $\eta(t)\in W^{2,2}(\Gamma)$ and moreover the following holds
		\begin{equation}\label{coercivityrel}
			\sup_{t\in[0,T]}\int_{\Gamma}\overline{\gamma}^{2}(\eta)|\nabla^{2}\eta|^{2}\leqslant cC_{0},
		\end{equation} 
		where $c$ depends only on $\varphi$. Furthermore, let $\overline{\gamma}(\eta_{0})>0$ then there is a minimal time $T_*$ depending on the initial configuration such that $\overline{\gamma}(\eta)>0$ and \eqref{coercivityrel} hold in $(0,T_{*})$.
		\begin{proof}
			Although there is a change in the structure of $K_{\delta}$ in comparison with the energy considered in the proof of the $W^{2,2}$ coercivity inequality for the displacement in \cite[Lemma 4.3]{MuhaSch}, this change does not affect the proof itself and we refer the reader interested in details therein. \\
			In order to prove the existence of a minimal time until which $\overline{\gamma}(\eta)>0$ and \eqref{coercivityrel} hold, we notice that the assumption $\bar\gamma(\eta_0)>0$ with $\bar\gamma$ defined in \eqref{ovgamma} and the continuity of $z\mapsto\bar\gamma(z)$ imply the existence of a constant $C$ such that if $z\in L^\infty(I;L^\infty(\Gamma))$ satisfies $\|z-\eta_0\|_{L^\infty(I;L^\infty(\Gamma))}\leq C$ on a suitable time interval $I$ then $\bar\gamma(z)>0$. 
			For the function $\eta$ we have by the interpolation inequality $\|g\|_{L^\infty(\Gamma)}\leq c\|g\|^\frac{1}{9}_{L^2(\Gamma)}\|g\|^\frac{8}{9}_{W^{1,4}(\Gamma)}$ and the $L^{\infty}(0,T;W^{1,4}(\Gamma))\cap W^{1,\infty}(0,T;L^{2}(\Gamma))$ uniform bound on $\eta$ obtained from \eqref{abndonelenergy} (as explained in the proof of \cite[Lemma 4.2]{MuhaSch}) that
			\begin{equation*}
				\|\eta(t)-\eta^{\delta}_0\|_{L^\infty(\Gamma)}\leq ct^\frac{1}{9}\|\eta\|_{W^{1,\infty}(0,T;L^2(\Gamma))}^\frac{1}{9}\|\eta-\eta_0^{\delta}\|_{L^\infty(0,T;W^{1,4}(\Gamma))}^\frac{8}{9}\leq ct^\frac{1}{9}
			\end{equation*}
			with a constant $c$ uniform with respect to the parameter $\delta$.
			Especially, we can fix $T$ such that
			\begin{equation*}
				\|\eta-\eta_0^{\delta}\|_{L^\infty(0,T_*;L^\infty(\Gamma))}\leq cT_*^\frac{1}{9}\leq \frac{C}{2}.
			\end{equation*}
			Then we estimate
			\begin{equation*}
				\|\eta-\eta_0\|_{L^\infty(0,T_*;L^\infty(\Gamma))}\leq \|\eta_0-\eta_0^{\delta}\|_{L^\infty(\Gamma)}+\|\eta-\eta_0^{\delta}\|_{L^\infty(0,T_*;L^\infty(\Gamma))}\leq \frac{C}{2}+\frac{C}{2}
			\end{equation*}
			provided $\|\eta_0-\eta_0^{\delta}\|_{L^\infty(\Gamma)}\leq \frac{C}{2}$ for $\delta$ fixed sufficiently small which is achievable due to the uniform convergence $\eta_0^{\delta}\to\eta_0$ following from \eqref{regintdens}$_1$. Hence we conclude $\bar\gamma(\eta)>0$ and \eqref{coercivityrel} holds for $T=T_*$.
		\end{proof}
	\end{lem}
	Next based on the uniform in $\delta-$ estimates we show that the degeneracy of he first kind \eqref{degenfstkind} can also be excluded for a positive minimal time.\\
	From \eqref{EtaDHCnv}, we have that for $(t,x)\in (0,T)\times\Gamma$
	\begin{equation*}
		\eta^{\delta}_{0}-ct^\frac{1}{4}\leq \eta^{\delta}_{0}(x)-|\eta^\delta(t,x)-\eta^{\delta}_{0}(x)|\leq \eta^\delta(t,x)\leq \eta^{\delta}_{0}(x)+|\eta^\delta(t,x)-\eta^{\delta}_{0}(x)|\leq \eta^{\delta}_{0}+ct^\frac{1}{4}
	\end{equation*}
	with $c$ independent of all regularizing parameters
	we can always choose a minimal time $T>0,$ such that the following holds
	\begin{equation}\label{PointBounds}
		a_{\partial\Omega}<m\leq \eta^\delta(t,x)\leq M<b_{\partial\Omega}
	\end{equation}
	for $(t,x)\in(0,T)\times\Gamma$ with $m=\min_{x\in\Gamma}\eta_{0}-cT^\frac{1}{4}$ and $M=\max_{x\in\Gamma}\eta_{0}+cT^\frac{1}{4}$ as $\min_{\Gamma}\eta_0\leq\eta^\delta_0\leq\max_{\Gamma}\eta_0$ by the definition of $\eta_0^\delta$ and $\eta_0\in(a_{\partial\Omega},b_{\partial\Omega})$ by assumption. Since \eqref{PointBounds} holds uniformly in $\delta,$ we notice that a minimal time can always be chosen so that the degeneracy of the first kind \eqref{degenfstkind} can be avoided.\\
	We further would like to summarize that as a consequence of Lemma \ref{Lem:Coer} and \eqref{PointBounds}, we can always choose a minimal time where it is possible to avoid both kind of degeneracy \eqref{degenfstkind} and \eqref{degensndkind}.\\
	
	\subsubsection{Maximal interval of existence}\label{Maxtimeex}
	In Section \ref{minintex}, we have shown that for a positive minimal time degeneracy of the solution can be excluded. We can start by solving the problem for this minimal time $T_{\min}$. Next considering  $(\eta,\tder\eta,\rho, Z,(\rho+Z)u)(T_{\min})$ as new initial conditions we repeat the existence proof in the interval $(T_{\min},2T_{\min}).$ We can iterate the procedure unitl a degeneracy occurs. That way we obtain a maximal time $T_F\in(0,\infty]$ such that $(\rho,Z,u,\eta)$ is a weak solution on the interval $(0,T)$ for any $T<T_F$. If $T_F$ is finite than either the $W^{2,2}$--coercivity of the Koiter energy is violated or $\lim_{s\to T_F}\eta(s,y)\in\{a_{\partial\Omega},b_{\partial\Omega}\}$. The extension procedure is standard nowadays and we refer to \cite[Theorem 3.5]{LeRu14} and \cite[Theorem 1.1]{MuhaSch} for details.
	
	\subsection{Proof of Case II:}\label{ProofcaseII}
	Here we will comment on the adaptions of arguments presented in Section \ref{ProofcaseI} and new regularity of the velocity field which will lead to the proof of $\textit{Case II}$ of Theorem \ref{Thm:main}. Indeed for the proof of $\textit{Case II}$, we set 
	\begin{equation*}
		\omega=\delta
	\end{equation*}
	in \eqref{momentumex}-\eqref{energybalanceex}, fix the value of the dissipation parameter $\zeta>0$ and perform the limit passage $\delta\to 0_+$.\\
	Since $\zeta>0$ is fixed, now the approximates $(\eta^{\delta},\rho^{\delta},Z^{\delta},u^{\delta})$ additionally solve
	\begin{equation}\label{regdelteta}
		\|\nabla\partial_{t}\eta^{\delta}\|_{L^{2}((0,T)\times\Gamma)}\leq c
	\end{equation}
	(where $c$ is independent of $\delta$) along with the other estimates listed in \eqref{UnifEstPhDom}. Indeed this additional bound leads to the convergence
	\begin{equation}\label{adconvetadelta}
		\begin{array}{l}
			\displaystyle \eta^{\delta}\rightharpoonup \eta\,\,\mbox{in}\,\, W^{1,2}((0,T)\times\Gamma).
		\end{array}
	\end{equation}
	Now we claim that $u^{\delta}$ (hence $u$) admits of an extension still denoted by $u_{\delta}$ (and consequently $u$ for the limit) outside $\Omega_{\eta^{\delta}}$ (and $\Omega_{\eta}$ for $u$) such that
	\begin{equation}\label{extraregetadel}
		\begin{array}{l}
			\|u^{\delta}\|_{L^{2}(0,T;W^{1,2}(B))}\leq c,\,\,\|u\|_{L^{2}(0,T;W^{1,2}(B))}\leq c,
		\end{array}
	\end{equation}
	where $c$ is independent of $\delta$ and $B$ is a neighborhood of both $\Omega_{\eta^{\delta}}$ (for all $\delta$) and $\Omega_{\eta}$ which was introduced in \eqref{defB}. The proof will be consequence of the extra regularity \eqref{regdelteta} of the structure due to dissipation. Provided the structure undergoes a frictional dissipation, \eqref{extraregetadel} improves the estimates \eqref{UestOnB} and \eqref{weaklimiturhoeta}$_{5}$ for non-dissipative hyperbolic elastic strmaxucture.\\
	For the proof of \eqref{extraregetadel} we introduce an extension operator $\mathcal F_\eta$ of functions defined on $\Gamma$ to $B.$ This extension operator is defined as follows
	\begin{equation}\label{FEDef}
		\mathcal{F}_{\eta}b=\mathcal F_{\Omega}((b\nu\circ\varphi^{-1}))\circ(\tilde\varphi_\eta)^{-1}
	\end{equation}  
	keeping in mind \hyperlink{Assumption}{Assumptions (A)}.  We note that $\mathcal F_\Omega$ stands for a composition of the standard extension of Sobolev functions on $\Omega$ to the whole space and the right inverse of the trace operator.\\ Such an extension was used in \cite{LeRu14} and \cite{Breit2} for different purpose. For the properties of $\mathcal{F}_{\eta}$ we refer the readers to \cite[Section 2.3.]{Breit2}.\\
	Indeed by construction $\mathcal{F}_{\eta}b\mid_{\partial\Omega_{\eta}}=\mathcal{F}_{\eta}b\circ\tilde{\varphi}_{\eta}=b\nu.$\\
	Now as a first step of the proof of \eqref{extraregetadel} let us show the following with the aid of a lifting argument involving $\mathcal{F}_{\eta}:$
	\begin{equation}\label{udeltaOmegaeta}
		\begin{array}{l}
			\displaystyle \|u^{\delta}\|_{L^{2}(0,T;W^{1,2}(\Omega_{\eta^{\delta}}))}\leq c,
		\end{array}
	\end{equation}
	for some $c$ independent of $\delta.$ Let us define
	\begin{equation}\label{wdelta}
		\begin{array}{l}
			\displaystyle  w^{\delta}=u^{\delta}(t,x)-\mathcal{F}_{\eta^{\delta}}\partial_{t}\eta^{\delta}(t,x)\,\,\mbox{for all}\,\,(t,x)\in (0,T)\times \Omega_{\eta^{\delta}},
		\end{array}
	\end{equation}
	Since $u^{\delta}\mid_{\partial\Omega_{\eta^{\delta}}}=(\partial_{t}\eta^{\delta}\nu)\circ\varphi^{-1}_{\eta^{\delta}}$ (indeed we are identifying $\partial\Omega$ and $\Gamma$ as stated in Remark~\ref{Rem:SimplNot} and skip writing $\varphi^{-1}$ throughout), one at once obtains that $w^{\delta}$ vanishes on $\partial\Omega_{\eta^{\delta}}.$ Since $\partial_{t}\eta^{\delta}\nu$ is bounded uniformly in $L^{2}(0,T;W^{1,2}(\Gamma))$ we have the following
	\begin{equation}\label{estliftparteta}
		\begin{array}{ll}
			&\displaystyle \|\mathcal{F}_{\eta^{\delta}}\partial_{t}\eta^{\delta}\|_{L^{2}(0,T;W^{1,2}(B))}\leq c\|\partial_{t}\eta^{\delta}\|_{L^{2}(0,T;W^{1,2}(\Gamma))}
		\end{array}
	\end{equation}
	for some constant $c$ independent of $\delta$ (we refer to \cite[Lemma 2.7, item (a)]{Breit2} for this estimate). Notice that \eqref{estliftparteta} is a time integrated version of the inequality stated in \cite[Lemma 2.7, item (a)]{Breit2}. This is possible to achieve since the constant appearing in \cite[Lemma 2.7, item (a)]{Breit2} depends only on $\Omega,$ $\|\eta^{\delta}\|_{W^{2,2}(\Gamma)},$ $b_{\partial\Omega}$ and $M$ and in the present scenario $\eta^{\delta}\in W^{2,2}(\Gamma)$ both uniformly in time and the parameter $\delta.$\\
	In view of \eqref{UnifEstPhDom2} 
	and \eqref{wdelta}, one at once concludes that 
	\begin{equation}\label{EwdelL2}
		\|\mathbb{D}w^{\delta}\|_{L^{2}(Q^{T}_{\eta^{\delta}})}\leq c, 
	\end{equation}
	where $c$ is independent of $\delta.$\\
	Let us now compute the following:\\
	\begin{equation}\label{Sdelw}
		\begin{split}
			\|\mathbb{D} w^{\delta}\|^{2}_{L^{2}(Q^{T}_{\eta^{\delta}})}&=\frac{1}{4}\int_{Q^{T}_{\eta^{\delta}}}(\nabla w^{\delta}+\nabla^\top w^{\delta})\cdot(\nabla w^{\delta}+\nabla^\top w^{\delta})\\
			&= \frac{1}{4}\int_{Q^{T}_{\eta^{\delta}}}\left(|\nabla w^{\delta}|^{2}+|\nabla^\top w^{\delta}|^{2}+2|\dvr\,w^{\delta}|^{2}\right),
		\end{split}
	\end{equation}
	where we have used integration by parts in space variables (which is justified since the boundary of $Q^{T}_{\eta^{\delta}}$ is Lipschitz in space for a.e. $t\in[0,T],$ we refer to Section \ref{compactnessshelenergy} for details), by a density argument and the fact that $w^{\delta}$ vanishes on $\partial\Omega_{\eta^{\delta}}.$\\
	Now \eqref{EwdelL2} and \eqref{Sdelw} together furnish that
	\begin{equation}\label{bndgradw}
		\begin{array}{l}
			\displaystyle \|w^{\delta}\|_{L^{2}(0,T;W^{1,2}(\Omega_{\eta^{\delta}}))}\leq c.
		\end{array}
	\end{equation}
	Finally in view of \eqref{wdelta}, \eqref{estliftparteta} and \eqref{bndgradw} we conclude the proof of \eqref{udeltaOmegaeta}.\\
	Next we plan to show that there exists an extension of $u^{\delta},$ still denoted by the same such that
	\begin{equation}\label{extenu}
		\begin{array}{l}
			\displaystyle \|u^{\delta}\|_{L^{2}(0,T;W^{1,2}(B))}\leq c,
		\end{array}
	\end{equation}
	for some $c$ independent of $\delta.$ Note that such an extension is necessary for the application of the almost compactness Lemma \eqref{Lem:AlmComp}. 
	\begin{remark}
		Note that for the proof of \eqref{extenu} we can not directly use Lemma \ref{Lem:Extension}, since that would render a loss of regularity and the extended $u^{\delta}$ will only belong to $L^{2}(0,T;W^{1,r}(B))$ for $r<2.$ We further specify that the proof of \eqref{extenu} will rely strongly on the dissipation of the structure. 
	\end{remark}
	For the proof of \eqref{extenu} we first extend the function $w^{\delta}$ by zero in $B\setminus \Omega_{\eta^{\delta}}$ (equivalently defining $u^{\delta}=\mathcal{F}_{\eta^{\delta}}\partial_{t}\eta^{\delta}$ in $B\setminus\Omega_{\eta^{\delta}}$) for $a.e.$ $t\in[0,T]$. Since for $a.e.$ $t,$ $\Omega_{\eta^{\delta}}$ has a Lipschitz boundary (we recall the improved $L^{2}(0,T;W^{2+,2}(\Gamma))$ regularity of $\eta^{\delta}$ proved in Section \ref{compactnessshelenergy}) and $w^{\delta}\in L^{2}(0,T;W^{1,2}_{0}(\Omega_{\eta^{\delta}})),$ the zero extension of $w^{\delta}$ (still denoted by the same) belongs to $W^{1,2}(B)$ for $a.e$ $t$. Since $\|w^{\delta}\|_{W^{1,2}(B)}=\|w^{\delta}\|_{W^{1,2}(\Omega_{\eta^{\delta}})}$ for $a.e.$ $t,$ one further has 
	\begin{equation}\label{normextwdel}
		\|w^{\delta}\|_{L^{2}(0,T;W^{1,2}(B))}=\|w^{\delta}\|_{L^{2}(0,T;W^{1,2}(\Omega_{\eta^{\delta}}))}.
	\end{equation}
	In view of \eqref{normextwdel}, \eqref{estliftparteta} and \eqref{regdelteta} we conclude the proof of \eqref{extenu}.\\
	Defining $u$ as the weak limit of $u^{\delta}$ in $(0,T)\times B$ one finishes the proof of \eqref{extraregetadel}.\\
	The part of the proof which differs 'Case I' (presented in Section \ref{ProofcaseI}) with 'Case II' is the compactness of the pressure. Let us the recall the arguments used to show \eqref{PressLimId} from \eqref{PiDL1W}, where the almost compactness compactness argument, more precisely Lemma \ref{Lem:AlmComp} was used. Here we can still use Lemma \ref{Lem:AlmComp} associated with the adiabatic exponents $\max\{\gamma,\beta\}=2$ (we refer to the second case of \eqref{valqgb2}). This is doable in view of the improved regularities presented in \eqref{extraregetadel}. Hence similar to the proof of 'Case I,' we can still freeze one of the densities to infer \eqref{PressLimId}.\\
	The final task is to identify the limit of the pressure. This will be done by modifying some of the arguments presented in Section \ref{conc}. More specifically to handle the critical case $\max\{\gamma,\beta\}=2,$ we need an estimate of the density oscillation presented in the following section.
	\subsubsection{Controlling the amplitude of density oscillations}\label{ampdenosc} First we extend $\rho^{\delta}$ and $\rho$ by zero outside $\Omega_{\eta^{\delta}}$ and $\Omega_{\eta}$ respectively. Further in view of \eqref{extraregetadel}, $u^{\delta}$ and $u$ can always be considered as functions defined on $\RR^{3}$ with uniform in $\delta$ estimates in the space $L^{2}(0,T;W^{1,2}(\mathbb{R}^{3})).$ Indeed this follows by a simple cut-off argument by using a cut-off function with value one in a  neighborhood of $\Omega_{\eta^{\delta}}$ (such a cut-off is possible since $Q^{T}_{\eta^{\delta}}$ is uniformly H\"{o}lder w.r.t $\delta$, $cf.$ \eqref{EtaDHCnv}) contained in $B$ and zero outside $B.$\\ 
	Next we need the following inequality which estimates the amplitude of density oscillations for the case $\max\{\gamma,\beta\}=2$ and $\min\{\gamma,\beta\}>0$:
	\begin{equation}\label{osdensity}
		\sup_{k>0}\limsup_{\delta\rightarrow 0}\int_{(0,T)\times\mathbb{R}^{3}}|T_{k}(\rho^{\delta})-T_{k}(\rho)|^{\max\{\gamma,\beta\}+1}=\sup_{k>0}\limsup_{\delta\rightarrow 0}\int_{(0,T)\times\mathbb{R}^{3}}|T_{k}(\rho^{\delta})-T_{k}(\rho)|^{3}\leqslant C.
	\end{equation}
	The proof of \eqref{osdensity} can be done by following the line of arguments used to show \cite[Proposition 14]{NovoPoko}. The only difference is the critical adiabatic exponent in \cite{NovoPoko} is $\frac{9}{5}$ whereas for our case it is $\max\{\gamma,\beta\}=2.$ Note that while adapting the arguments of \cite{NovoPoko} to show \eqref{osdensity}, one needs to apply the assumption \eqref{?!} concerning the structure of $\mathcal{P}$ when $\max\{\gamma,\beta\}=2.$
	\begin{remark}
		One recalls from the theory of existence of weak solutions for compressible viscous fluids in a Lipschitz domain, that when the adiabatic exponent $\gamma>\frac{9}{5},$ (consequently $\gamma+\gamma_{BOG}>2$ (we refer to \cite{NovStr04} for details)) and the velocity field belongs to $L^{2}(W^{1,2})$ one can prove the existence of renormalized weak solutions to the continuity equations, which is a crucial tool to prove the strong convergence of $\rho^{\delta}.$ A very important observation of \cite{FeireislPet} is to prove an inequality of the form \eqref{osdensity}, in order to show the validity of the renormalized weak solutions for the continuity equation even when the adiabatic exponent $\gamma$ lies in $(\frac{3}{2},\frac{9}{5}].$ Notice that we are working with a moving boundary which is non-Lipschitz and hence we have no improved integrability of the pressure up to the interface by using Bogovskii type argument. That is the reason we separately prove the existence of renormalized weak solutions to the continuity equations when the adiabztic exponent takes value $\geq 2$ up to the to interface in form of Lemma \ref{Lem:Renormalization} and the arguments used in proving the lemma is independent of the inequality \eqref{osdensity}. Despite of this fact we still need \eqref{osdensity}, to deal with the borderline case $\max\{\gamma,\beta\}=2$ and this will be apparent from the analysis done next. 
	\end{remark}
	Now we follow the analysis presented in Section \ref{conc}, by replacing the estimate of $I^{\rho}_{1}$ in \eqref{I1E} with the following
	\begin{equation}\label{I1rhoalter}
		\begin{array}{ll}
			|I^{\rho}_{1}|&\displaystyle \leqslant C\limsup_{\delta\rightarrow 0} \|T_{k}(\rho^{\delta})-{T_{k}(\rho)}\|^{\frac{1}{4}}_{L^{1}((0,T)\times\mathbb{R}^{3})}.
		\end{array}
	\end{equation}
	The last estimate is a consequence of \eqref{osdensity} and $u\in L^{2}(0,T;W^{1,2}(\mathbb{R}^{3})).$ The rest of the analysis can be carried out as it is in Section \ref{conc} until the final conclusion \eqref{RhoDPoint} is made.
	\subsection{Summary of the proof of Theorem \ref{Thm:main}:} First let us summarize the proof of 'Case I' of Theorem \ref{Thm:main}. Since $\rho^{\delta}, Z^{\delta}\geq 0$ for each $\delta,$ from Theorem \ref{resultstaulayer} and the positivity is preserved under weak-convergence, we conclude $\rho, Z\geq 0.$ That $\rho$ and $Z$ are weakly continuous in time with values in $L^{\max\{\gamma,\beta\}}(\mathbb{R}^{3})$ can be obtained from \eqref{AddConvDel}. The regularity of the velocity field $u$ follows from \eqref{weaklimiturhoeta}$_{5}.$ The regularities \eqref{listregsol}$_{5,6}$ are consequences of \eqref{AddConvDel}$_{3,5}.$ The regularity \eqref{listregsol}$_{7}$ of $\eta$ follows from \eqref{weaklimiturhoeta}$_{1,2}$ and \eqref{improvedregeta}. That $P(\rho,Z)$ belongs to $L^{1}(Q^{T}_\eta)$ follows from \eqref{L1convpress}.\\
	For the proof of the continuity of the fluid and structural velocities we refer to Section \ref{contvelo}.\\
	The momentum balance \eqref{momentum} is recovered from \eqref{LimDeltaIdent} and \eqref{RhoDPoint}.\\
	In view of \eqref{AddConvDel}$_{3}$ one can easily pass limit in \eqref{contrhoex} solved by $\rho^{\delta},$ $Z^{\delta}$ and $u^{\delta}$ to recover the continuity equations \eqref{contrho}.\\
	We have shown the validity of the energy inequality \eqref{energybalance} in Section \ref{finalenergybalance}.\\
	The attainment of the initial data in a weak sense is explained in Section \ref{attinitialdata}.\\
	We present the proof of 'Case II' of Theorem \ref{Thm:main} in Section \ref{ProofcaseII}.
	\section{Appendix}
	%%%%%%%%%%%%%%%%%%%%%%%%%%%%%%%%%%%%%%%%%%%%%%%%%%%%%%%%%%%%%%%%%%%%%%%%%%%%%%%%%%%%%%%%%%%%%%%%%%%%%%%%%%%%%%%%%%%%%%%%%%%%%%%%%%%%%%%%%%%%%%%%%%%%%%%%%%%%%%%%%%%%%%%%%%%%%%%%%%%%%%%%%%%%%%%%%%%%%%%%%%%%%%%%%%%%%%%%%%%%%%%%%%%%%%%%%%%%%%%%%%%%%%%%%%%%%%%%%%%%
	The ensuing lemma states that a weak limit of sequence of functions vanishing outside a corresponding varying domain vanishes outside of a varying domain that corresponds to an uniform limit of displacements.
	\begin{lem}\label{Lem:VanSeq}
		Let $\{\eta^i\}\subset C([0,T]\times\Gamma)$ be such that $\eta^i\to\eta$ uniformly on $[0,T]\times\Gamma$ and $\{h^i\}\subset L^1((0,T)\times B)$ be such that $h^i\rightharpoonup h$ in $L^1((0,T)\times B)$ and $h^i\equiv 0$ a.e. in $\left((0,T)\times B\right)\setminus Q^T_{\eta^i}$. Then $h\equiv 0$ a.e. in $\left((0,T)\times B\right)\setminus Q^T_\eta$.
		\begin{proof}
			Let $K\subset\left((0,T)\times B\right)\setminus Q^T_\eta$ be a compact set. As the uniform convergence $\eta^i\to\eta$ in $C([0,T]\times\Gamma)$ is assumed, there is an index $i_0$ such that $K\subset \left((0,T)\times B\right)\setminus Q^T_{\eta^i}$ for each $i>i_0$. Then for an arbitrary $\vartheta\in C(K)$ it follows that 
			\begin{equation}
				\int_{K}h^i\vartheta=0.
			\end{equation}
			Hence we conclude 
			\begin{equation}
				\int_{(0,T)\times B}h\vartheta=\lim_{i\to \infty}\int_{(0,T)\times B}h^i\vartheta=0.
			\end{equation}
			for any $\vartheta\in C_c\left(\left(0,T)\times B\right)\setminus Q^T_\eta\right)$ implying $h\equiv 0$ a.e. in $\left((0,T)\times B\right)\setminus Q^T_\eta$.
		\end{proof}
	\end{lem}
	
	\begin{lem}\label{lem:WLP}
		Let $\mathcal O\subset\eR^d$ be a domain, $P,G:\mathcal O\times[0,\infty)\to[0,\infty)$ be a couple of functions such that for almost all $y\in \mathcal O$ the mappings $z\mapsto P(y,z)$ and $z\mapsto G(y,z)$ are continuous and nondecreasing on $[0,\infty)$. Suppose that $\{z^n\}\subset L^1(\mathcal O;[0,\infty))$ is a sequence such that
		\begin{equation*}
			\begin{alignedat}{2}
				P(\cdot, z^n)\rightharpoonup& \overline{P(\cdot,z)},\\
				G(\cdot, z^n)\rightharpoonup& \overline{G(\cdot,z)},\\
				P(\cdot, z^n)G(\cdot, z^n)\rightharpoonup& \overline{P(\cdot,z)G(\cdot, z)},\\
			\end{alignedat}
		\end{equation*}
		in $L^1(\mathcal O)$.
		Then
		\begin{equation}
			\overline{P(\cdot,z)}\ \overline{G(\cdot,z)}\leq \overline{P(\cdot,z)G(\cdot, z)}\text{ a.e. in }\mathcal O.
		\end{equation}
	\end{lem}
	%%%%%%%%%%%%%%%%%%%%%%%%%%%%%%%%%%%%%%%%%%%%%%%%%%%%%%%%%%%%%%%%%%%%%%%%%%%%%%%%%%%%%%%%%%%%%%%%%%%%%%%%%%%%%%%%%%%%%%%%%%%%%%%%%%%%%%%%%%%%%%%%%%%%%%%%%%%%%%%%%%%%%%%%%%%%%%%%%%%%%%%%%%%%%%%%%%%%%%%%%%%%%%%%%%%%%%%%%%%%%%%%%%%%%%%%%%%%%%%%%%%%%%%%%%%%%%%%%%%%%%%%%%%%%%%%%%%%%%%%%%%%%%%%%%%%%%%%%%%%%%%%%%%%%%%%%%%%%%%%%%%%%%%%%%%%%%%%%%%%%%%%%%%%%%%%%%%%%%%%%%%%
	We will use several times the compactness result in the ensuing lemma, for its proof see \cite[Theorem 5]{Sim87}.
	\begin{lem}\label{Lem:RelComp}
		Let $T>0$, $p\in[1,\infty]$ and Banach spaces $X_1,X_2,X_3$ satisfy $X_1\stackrel{C}{\hookrightarrow}X_2\hookrightarrow X_3$. Assume that $F\subset L^p(0,T;X_1)$ fulfills
		\begin{enumerate}
			\item $\sup_{f\in F}\|f\|_{L^p(0,T X_1)}<\infty$,
			\item $\sup_{f\in F}\|\tau_s f-f\|_{L^p(0,T-s;X_3)}\to 0$ as $s\to 0$.
		\end{enumerate}
		Then $F$ is relatively compact in $L^p(0,T;X_2)$ and $C([0,T];X_2)$ if $p=\infty$.
	\end{lem}
	The following lemma is a particular case of a more abstract result, cf.\ \cite[Lemma 9.1]{Alt12}.
	\begin{lem}\label{Lem:RefIneq}
		Let $M\in\eN$, a Hilbert space $H$ and $\{f^{m}\}\subset H$ be given. Moreover, assume that the function $f^{M}$ being defined via $f^{M}(t)=f^m$ for $t\in [(m-1)\Delta t,m\Delta t)$, $m\in\eN,$ satisfies
		\begin{equation}\label{BoundIntDiff}
			\int_0^{kh-s}\|f^M(t+s)-f^M(t)\|^2_H\dt\leq c s^{q}
		\end{equation}
		where $s=lh$, $l\in\eN$, $l\leq k$ and $q\in(0,1]$. Then \eqref{BoundIntDiff} holds with any $0<s<kh.$
	\end{lem}
	
	%%%%%%%%%%%%%%%%%%%%%%%%%%%%%%%%%%%%%%%%%%%%%%%%%%%%%%%%%%%%%%%%%%%%%%%%%%%%%%%%%%%%%%%%%%%%%%%%%%%%%%%%%%%%%%%%%%%%%%%%%%%%%%%%%%%%%%%%%%%%%%%%%%%%%%%%%%%%%%%%%%%%%%%%%%%%%%%%%%%%%%%%%%%%%%%%%%%%%%%%%%%%%%%%%%%%%%%%%%%%%%%%%%%%%%%%%%%%%%%%%%%%%%%%%%%%%%%%%%%%%%%%%%%%%%%%%%%%%%%%%%%%%%%%%%%%%%%%%%%%%%%%%%%%%%%%%%%%%%%%%%%%%%%%%%%%%%%%%%%%%%%%%%%%%%%%%%%%%%%%%%%%%%%%%%%%%%%%%%%%%%%%%%%%%%%%%%%%%%%%%%%%%%%%
	\subsection{Proof of Theorem \ref{labelstructuralsp}}\label{proofdecstr}
	\begin{proof}
		We divide the proof into three parts and they are presented in three sections.
		\subsubsection{Existence of a time discrete problem, $\Delta t<< \tau$ layer}
		In this section we will further divide the time interval $(0,\tau)$ into subintervals of length $\Delta t<<\tau$ and introduce a further discretization of the structural subproblem. Compared to the $\tau-$ layer here we discretize the structural velocity $\partial_{t}\eta.$\\ 
		{\textit{\bf Introducing further time discretization and a fixed point map:}} We first assume for $m\in\eN$ $(\eta^{m},w^{m})\in W^{3,2}(\Gamma)\times L^{2}(\Gamma)$ and solve for $(\eta^{m+1},w^{m+1})$ in the following discrete problems (in the following the time discretization $\Delta t \ll \tau$)
		\begin{equation}\label{subsubstructure}
			\begin{split}
				&\int_{\Gamma}\frac{\eta-\eta^{m}}{\Delta t}b_{1}=\int_{\Gamma}wb_{1},\\
				& (1-\delta)\int_{\Gamma}\frac{w-w^{m}}{\Delta t}b+\delta\int_{\Gamma}\frac{w-v^{n}\cdot \nu}{\tau}b+\zeta\int_{\Gamma}\nabla w\cdot\nabla b+\langle K'_{\delta}(\eta,\eta^{m}),b\rangle=0
			\end{split}
		\end{equation}    
		where $K'_{\delta}(\eta^{m+1},\eta^{m})$ approximates $K'_{\delta}(\eta)$ and $\langle K'_{\delta}(\eta^{m+1},\eta^{m}),b\rangle$ is given as follows
		\begin{equation}\label{Kepsilon0b}
			\langle K'_{\delta}(\eta,\eta^{m}),b\rangle=\frac{h}{2}\int_{\Gamma}\mathcal{A}\mathbb{G}(\eta):\mathbb{G}'(\eta,\eta^{m})b+\frac{h^{3}}{24}\int_{\Gamma}\mathcal{A}\mathbb{R}(\eta):\mathbb{R}'(\eta,\eta^{m})b+\delta^7\int_{\Gamma}\nabla^{3}\eta\cdot \nabla^{3}b,
		\end{equation}
		with $(b_{1},b)\in L^{2}(\Gamma)\times W^{3,2}(\Gamma)$.
		
		Following \cite{MuhaSch}, above we have approximated ${\GG}'(\eta)$ and ${\RR}'(\eta)$ as follows
		\begin{equation}\label{discreteGR'}
			\begin{split}
				& {\GG}'(\eta,\eta^{m})b=\frac{1}{6}\bigg(\GG'(\eta^{m})+4\GG'(\overline{\eta})+\GG'(\eta)\bigg)b,\\
				& \RR'(\eta,\eta^{m})b=\frac{1}{6}\bigg(\RR'(\eta^{m})+4\RR'(\overline{\eta})+\RR'(\eta)\bigg)b,
			\end{split}
		\end{equation}
		where $\overline\eta=\frac{\eta+\eta^m} {2}$ and tensors $\GG'$, $\RR'$ are defined in Section \ref{Sec:KE}.
		One notices that an approximation of the form \eqref{discreteGR'} is useful, since while testing the approximated equations by $\displaystyle\frac{\eta-\eta^{m}}{\Delta t}$ (a discrete version of $\partial_{t}\eta$), one has
		\begin{equation}\label{discreteaftrtesttimeder}
			\begin{split}
				& \GG'(\eta,\eta^{m})\frac{\eta-\eta^{m}}{\Delta t}=\frac{1}{\Delta t}\left(\GG(\eta)-\GG(\eta^{m})\right),\\
				& \RR'(\eta,\eta^{m})\frac{\eta-\eta^{m}}{\Delta t}=\frac{1}{\Delta t}\left(\RR(\eta)-\RR(\eta^{m})\right).
			\end{split}
		\end{equation}
		
		Next eliminating $w$ from \eqref{subsubstructure}$_{2}$ we obtain
		\begin{equation}\label{afterelew}
			\begin{split}
				&(1-\delta+\frac{\delta\Delta t}{\tau})\int_{\Gamma}\eta\cdot b+\delta(\Delta t)^{2}\int_{\Gamma}\nabla^{3}\eta\cdot\nabla^{3}b+\zeta\Delta t\int_{\Gamma}\nabla\eta\cdot\nabla b\\
				&=-\frac{h}{2}(\Delta t)^{2}\int_{\Gamma}\mathcal{A}\mathbb{G}(\eta):\mathbb{G}'(\eta,\eta^{m})b-\frac{h^{3}}{24}(\Delta t)^{2}\int_{\Gamma}\mathcal{A}\mathbb{R}(\eta):\mathbb{R}'(\eta,\eta^{m})b+(1-\delta+\frac{\delta\Delta t}{\tau})\int_{\Gamma}\eta^{m}b\\
				&\quad+\zeta\Delta t\int_{\Gamma}\nabla\eta^{m}\cdot\nabla b+(1-\delta)\Delta t\int_{\Gamma}w^{m}\cdot b+\delta(\Delta t)^{2}\int_{\Gamma}(v^{n}\cdot\nu)b.
			\end{split}
		\end{equation}
		Notice that since at this level we are interested in finding a solution of \eqref{afterelew}, for fixed $\Delta t,$ $\tau,$ $\delta,$ $\eta^{m},$ $w^{m}$ and $v^{n}.$ Now to solve \eqref{afterelew}, we introduce a map
		\begin{equation}\label{fpm}
			\mathcal{F}:  W^{2,4}(\Gamma)\displaystyle\rightarrow W^{2,4}(\Gamma)\quad\mbox{such that}\,\, \widetilde{\eta}\longmapsto\mathcal{F}(\widetilde{\eta})
		\end{equation}
		where for $b\in W^{3,2},$ $\mathcal{F}(\widetilde{\eta})$ solves
		\begin{equation}\label{eqFteta}
			\begin{split}
				& (1-\delta+\frac{\delta\Delta t}{\tau})\int_{\Gamma}\mF(\weta)\cdot b+\delta(\Delta t)^{2}\int_{\Gamma}\nabla^{3}\mF(\weta)\cdot\nabla^{3}b+\zeta\Delta t\int_{\Gamma}\mF(\weta)\cdot\nabla b\\
				& =-\frac{h}{2}(\Delta t)^{2}\int_{\Gamma}\mathcal{A}\mathbb{G}(\weta):\mathbb{G}'(\weta,\eta^{m})b-\frac{h^{3}}{24}(\Delta t)^{2}\int_{\Gamma}\mathcal{A}\mathbb{R}(\weta):\mathbb{R}'(\weta,\eta^{m})b+(1-\delta+\frac{\delta\Delta t}{\tau})\int_{\Gamma}\eta^{m}b\\
				&\quad+\zeta\Delta t\int_{\Gamma}\nabla\eta^{m}\cdot\nabla b+(1-\delta)\Delta t\int_{\Gamma}w^{m}\cdot b+\delta(\Delta t)^{2}\int_{\Gamma}(v^{n}\cdot\nu)b=\sum_{i=1}^{6}\mathcal{L}_{i}b=\langle \mathcal{L},b \rangle.
			\end{split}
		\end{equation}
		Since $\widetilde{\eta}\in W^{2,4}(\Gamma),$ in view of the structure of $\langle\mathcal{L}_{1},b\rangle$ and $\langle\mathcal{L}_{2},b\rangle$ (we refer to \eqref{amab}, \eqref{Gij'} and \eqref{exaR}), $\mathcal{L}\in (W^{3,2}(\Gamma))'.$\\
		Now using Lax-Milgram theorem one at once proves the existence of a unique $\mathcal{F}(\weta)\in W^{3,2}(\Gamma)$ which solves \eqref{eqFteta}. The continuous embedding $W^{3,2}(\Gamma) \hookrightarrow W^{2,4}(\Gamma)$ renders the well-defineness of $\mathcal{F}.$\\
		Next to prove the existence of a fixed point of the map $\mathcal{F},$ we will use (as in \cite{CanicMuha}) Schaefer's fixed point theorem (the statement can be found in \cite[Theorem 4]{CanicMuha}). For that we first observe that $\mathcal{F}$ is compact since $W^{3,2}(\Gamma)$ is compactly embedded into $W^{2,4}(\Gamma).$\\
		{\textit{\bf Existence of fixed point of the map $\mF,$ Step-1 (boundedness independent of a parameter $\lambda$)}:} 
		Next we consider the operator equation $\weta=\lambda\mathcal{F}(\weta)$ for $\lambda\in[0,1],$ or in other words
		\begin{equation}\label{etlambda}
			\begin{split}
				&(1-\delta+\frac{\delta\Delta t}{\tau})\int_{\Gamma}\weta\cdot b+\delta(\Delta t)^{2}\int_{\Gamma}\nabla^{3}\weta\cdot\nabla^{3}b+\zeta\Delta t\int_{\Gamma}\nabla\weta\cdot\nabla b\\
				& =-\lambda\frac{h}{2}(\Delta t)^{2}\int_{\Gamma}\mathcal{A}\mathbb{G}(\weta):\mathbb{G}'(\weta,\weta^{m})b-\lambda\frac{h^{3}}{24}(\Delta t)^{2}\int_{\Gamma}\mathcal{A}\mathbb{R}(\weta):\mathbb{R}'(\weta,\eta^{m})b+\lambda(1-\delta+\frac{\delta\Delta t}{\tau})\int_{\Gamma}\eta^{m}b\\
				&\quad+\lambda\zeta\Delta t\int_{\Gamma}\nabla\eta^{m}\cdot\nabla b+\lambda(1-\delta)\Delta t\int_{\Gamma}w^{m}\cdot b+\lambda\delta(\Delta t)^{2}\int_{\Gamma}(v^{n}\cdot\nu)b,
			\end{split}
		\end{equation}   
		where $\lambda\in[0,1].$ In order to apply Schaefer's fixed point theorem, we first need to establish $W^{2,4}(\Gamma)$ estimate of $\widetilde{\eta}$ in $W^{2,4}(\Gamma)$ independent of $\lambda,$ by using \eqref{etlambda}.\\
		We introduce $\displaystyle\ww=\frac{\weta-\eta^{m}}{\Delta t}$ and rewrite \eqref{etlambda} as
		\begin{equation}\label{aftrreptw}
			\begin{split}
				& (1-\delta+\frac{\delta\Delta t}{\tau})\Delta t\int_{\Gamma}(\ww-w^{m})\cdot b+\delta(\Delta t)^{2}\int_{\Gamma}\nabla^{3}\weta\cdot\nabla^{3}b+\zeta\Delta t\int_{\Gamma}\nabla\weta\cdot\nabla b\\
				&\quad+\lambda\frac{h}{2}(\Delta t)^{2}\int_{\Gamma}\mathcal{A}\mathbb{G}(\weta):\mathbb{G}'(\weta,\weta^{m})b+\lambda\frac{h^{3}}{24}(\Delta t)^{2}\int_{\Gamma}\mathcal{A}\mathbb{R}(\weta):\mathbb{R}'(\weta,\eta^{m})b\\
				&=(\lambda-1)(1-\delta+\frac{\delta\Delta t}{\tau})\int_{\Gamma}\eta^{m}b
				\displaystyle\quad+\lambda\zeta\Delta t\int_{\Gamma}\nabla\eta^{m}\cdot\nabla b+\lambda\delta(\Delta t)^{2}\int_{\Gamma}(v^{n}\cdot\nu)b\\
				&\quad+\bigg(\lambda(1-\delta)-(1-\delta+\frac{\delta\Delta t}{\tau}\bigg)\Delta t\int_{\Gamma}w^{m}\cdot b.
			\end{split}
		\end{equation}
		Next since $\Delta t<< \tau,$ we can use $b=\widetilde{w}$ as a test function in \eqref{aftrreptw} to furnish
		\begin{equation}\label{usingbeqww}
			\begin{split}
				&\|\ww\|^{2}_{L^{2}(\Gamma)}+\|\ww-w^{m}\|^{2}_{L^{2}(\Gamma)}+\|\nabla^{3}\weta\|^{2}_{L^{2}(\Gamma)}+\|\nabla^{3}(\weta-\eta^{m})^{2}_{L^{2}(\Gamma)}\|\\
				&+\zeta\bigg(\|\nabla\weta\|^{2}_{L^{2}(\Gamma)}+\|\nabla(\weta-\eta^{m})\|^{2}_{L^{2}(\Gamma)}\bigg)\\
				&+\lambda\bigg(\int_{\Gamma}\mA\GG(\weta):\GG(\weta)
				+\int_{\Gamma}\mA(\GG(\weta)-\GG(\eta^{m})):(\GG(\weta)-\GG(\eta^{m})\bigg)\\
				&+\lambda\bigg(\int_{\Gamma}\mA\RR(\weta):\RR(\weta)
				+\int_{\Gamma}\mA(\RR(\weta)-\RR(\eta^{m})):(\RR(\weta)-\RR(\eta^{m})\bigg)\\
				&\leqslant C(\|w^{m}\|_{L^{2}(\Gamma)},\|\eta^{m}\|_{W^{3,2}(\Gamma)},\|v^{n}\|_{L^{2}(\Gamma)},h,\tau,\Delta t),
			\end{split}
		\end{equation}
		where during the calculation we have used H\"{o}lder and Young's inequality, to estimate the terms appearing in the right hand side of \eqref{aftrreptw} and to absorb $\varepsilon(\|\ww\|^{2}_{L^{2}(\Gamma)}+\|\weta\|^{2}_{L^{2}(\Gamma)})$ with suitable terms in the left hand side of \eqref{aftrreptw} for sufficiently small value of $\varepsilon.$
		From \eqref{usingbeqww}, one renders
		$$\|\weta\|_{W^{2,4}(\Gamma)}\leqslant C\|\weta\|_{W^{3,2}(\Gamma)}\leqslant C(\|w^{m}\|_{L^{2}(\Gamma)},\|\eta^{m}\|_{W^{3,2}(\Gamma)},\|v^{n}\|_{L^{2}(\Gamma)},h,\tau,\Delta t).$$
		{\textit{\bf Step-2 (Continuity of the map $\mF$)}:} One finally needs to verify the continuity of $\mathcal{F}$ in order to apply Schaefer's fixed point theorem to show the existence os a fixed point of the map $\mathcal{F}$ (introduced in \eqref{fpm}).\\
		In that direction let us assume that $\weta_{k}\rightarrow \weta$ in $W^{2,4}(\Gamma).$ We claim that 
		$$\mathcal{F}(\weta_{k})=\eta_{k}\rightarrow \mathcal{F}(\weta)=\eta\quad\mbox{in}\quad W^{2,4}(\Gamma).$$
		One observes that $r_{k}=\eta-\eta_{k}$ solves
		\begin{equation}\label{eqsolverk}
			\begin{split}
				& (1-\delta+\frac{\delta\Delta t}{\tau})\int_{\Gamma}r_{k}\cdot b+\delta(\Delta t)^{2}\int_{\Gamma}\nabla^{3}r_{k}\cdot\nabla^{3}b+\zeta\Delta t\int_{\Gamma}\nabla r_{k}\cdot\nabla b\\
				& =\frac{h}{2}(\Delta t)^{2}\int_{\Gamma}\mathcal{A}\mathbb{G}(\weta_{k}):\mathbb{G}'(\weta_{k},\eta^{m})b-\frac{h}{2}(\Delta t)^{2}\int_{\Gamma}\mathcal{A}\mathbb{G}(\weta):\mathbb{G}'(\weta,\eta^{m})b\\
				&\quad +\frac{h^{3}}{24}(\Delta t)^{2}\int_{\Gamma}\mathcal{A}\mathbb{R}(\weta_{k}):\mathbb{R}'(\weta_{k},\eta^{m})b-\frac{h^{3}}{24}(\Delta t)^{2}\int_{\Gamma}\mathcal{A}\mathbb{R}(\weta):\mathbb{R}'(\weta,\eta^{m})b,
			\end{split}
		\end{equation}
		where $b\in W^{3,2}(\Gamma).$\\
		Taking $r_{k}$ as a test function in \eqref{eqsolverk}, we render
		\begin{equation}\label{rkastfn}
			\begin{split}
				&(1-\delta+\frac{\delta\Delta t}{\tau})\|r_{k}\|^{2}_{L^{2}(\Gamma)}+\delta(\Delta t)^{2}\|\nabla^{3}r_{k}\|^{2}_{L^{2}(\Gamma)}+\zeta\Delta t\|\nabla r_{k}\|^{2}_{L^{2}(\Gamma)}\\
				& =\frac{h}{2}(\Delta t)^{2}\int_{\Gamma}\mathcal{A}\mathbb{G}(\weta_{k}):\mathbb{G}'(\weta,\eta^{m})r_{k}-\frac{h}{2}(\Delta t)^{2}\int_{\Gamma}\mathcal{A}\mathbb{G}(\weta):\mathbb{G}'(\weta,\eta^{m})r_{k}\\
				&\quad +\frac{h^{3}}{24}(\Delta t)^{2}\int_{\Gamma}\mathcal{A}\mathbb{R}(\weta_{k}):\mathbb{R}'(\weta,\eta^{m})r_{k}-\frac{h^{3}}{24}(\Delta t)^{2}\int_{\Gamma}\mathcal{A}\mathbb{R}(\weta):\mathbb{R}'(\weta,\eta^{m})r_{k}\\
				&=-\frac{h}{2}(\Delta t)^{2}\int_{\Gamma}\bigg(\mathcal{A}\mathbb{G}(\weta):(\mathbb{G}'(\weta,\eta^{m})-\GG'(\weta_{k},\eta^{m}))+\mA(\GG(\weta)-\GG(\weta_{k})):\GG'(\weta_{k},\eta^{m})\bigg)r_{k}\\
				&\quad-\frac{h}{2}(\Delta t)^{2}\int_{\Gamma}\bigg(\mathcal{A}\mathbb{R}(\weta):(\mathbb{R}'(\weta,\eta^{m})-\RR'(\weta_{k},\eta^{m}))+\mA(\RR(\weta)-\RR(\weta_{k})):\RR'(\weta_{k},\eta^{m})\bigg)r_{k}.
			\end{split}
		\end{equation}
		Since $\GG'(\eta)$ is linear in $\nabla\eta$ (cf. \eqref{Gij'}), $r_{k}$ is bounded in $W^{3,2}(\Gamma)$ uniformly in $k$ (since $\|\weta_{k}\|_{W^{2,4}(\Gamma)}$ can be bounded by $\|\weta\|_{W^{2,4}(\Gamma)}+1$ independently of $k$), the following convergence holds in view of the definition \eqref{discreteaftrtesttimeder}$_{1}$  
		$$
		\displaystyle(\mathbb{G}'(\weta,\eta^{m})-\GG'(\weta_{k},\eta^{m}))r_{k}\rightarrow 0\quad\mbox{in}\quad L^{2}(\Gamma).
		$$
		Further the above convergence combined with the boundedness of $\mA\GG(\eta)$ in $L^{2}(\Gamma)$ furnishes that
		\begin{equation}\label{G1con}
			\begin{array}{l}
				\displaystyle\int_{\Gamma}\mathcal{A}\mathbb{G}(\weta):(\mathbb{G}'(\weta,\eta^{m})-\GG'(\weta_{k},\eta^{m}))r_{k}\rightarrow 0.
			\end{array}
		\end{equation}
		The boundedness of $\GG'(\weta_{k},\eta^{m})r_{k}$ in $L^{2}(\Gamma)$ and the convergence $(\GG(\weta)-\GG(\weta_{k}))\rightarrow 0$ in $L^{2}(\Gamma)$ (which follows from the fact that $\weta_{k}\rightarrow \weta$ in $W^{2,4}(\Gamma)$ and \eqref{Geta}) readily implies
		\begin{equation}\label{G2con}
			\begin{array}{l}
				\displaystyle\int_{\Gamma}\mA(\GG(\weta)-\GG(\weta_{k})):\GG'(\weta_{k},\eta^{m})r_{k}\rightarrow 0.
			\end{array}
		\end{equation}
		Next one uses \eqref{Rijcomp} to observe that $\RR(\widetilde{\eta})$ is bounded in $L^{\infty}(\Gamma)$ (independent of $k$), when $\weta\in W^{2,4}(\Gamma).$ Further one uses \eqref{rewriteRij} and \eqref{discreteaftrtesttimeder}$_{2}$ to verify that  We now claim that
		\begin{equation}\label{apR1con}
			\begin{array}{l}
				(\RR'(\weta,\eta^{m})-\RR'(\weta_{k},\eta^{m}))r_{k}\rightarrow 0\,\,\mbox{in}\,\,L^{2}(\Gamma).
			\end{array}
		\end{equation}
		Since $\nabla^{2}\weta_{k}$ converges to $\nabla^{2}\weta$ in $L^{4}(\Gamma)$ and $\overline{\gamma}'(\weta_{k})$ converges to $\overline{\gamma}'(\weta)$ in $L^{p}(\Gamma)$ for any $p<\infty$ (follows since $\overline{\gamma}'(\eta)$ is linear in $\eta$), one has that $(\overline{\gamma}'(\weta_{k})b)\partial^{2}_{ij}\weta_{k}$ converges to $(\overline{\gamma}'(\weta)b)\partial^{2}_{ij}\weta$ in particular in $L^{2}(\Gamma)$ for any $b\in W^{3,2}(\Gamma).$ Further $P'_{0}(\eta,\nabla\eta)$ is a polynomial in $\eta$ and $\nabla \eta$ of order two one verifies that $P'_{0}(\weta_{k},\nabla\weta_{k})b$ converges to $P'_{0}(\weta,\nabla\weta)b$ in $L^{2}(\Gamma).$ In view of the aforementioned arguments we conclude \eqref{apR1con}. Next using \eqref{apR1con} and the boundedness of $\RR(\weta)$ in $L^{2}(\Gamma),$ we furnish
		\begin{equation}\label{R1con}
			\begin{array}{l}
				\displaystyle \int_{\Gamma}\mA\RR(\weta):(\mathbb{R}'(\weta,\eta^{m})-\RR'(\weta_{k},\eta^{m}))r_{k}\rightarrow 0.
			\end{array}
		\end{equation}
		Similar arguments with minor adaptations lead to the convergence of $\RR(\weta_{k})$ to $\RR(\weta)$ in $L^{2}(\Gamma)$ which combined with the boundedness of $\RR'(\weta_{k},\eta^{m})r_{k}$ in $L^{2}(\Gamma)$ yields
		\begin{equation}\label{R2con}
			\begin{array}{l}
				\displaystyle \int_{\Gamma} \mA(\RR(\weta)-\RR(\weta_{k})):\RR'(\weta_{k},\eta^{m})r_{k}\rightarrow 0.
			\end{array}
		\end{equation}
		The convergences \eqref{G1con}, \eqref{G2con}, \eqref{R1con} and \eqref{R2con} together with \eqref{rkastfn} implies that $\|r_{k}\|_{W^{3,2}(\Gamma)}$ or more particularly $\|r_{k}\|_{W^{2,4}(\Gamma)}$ converges to zero. This renders the continuity of the map $\mF.$\\
		{\bf{Conclusion about the existence:}} Finally we can apply Schaefer's fixed point theorem to show the existence of a fixed point of the map $\mathcal{F}$ and thereby proving the existence of a solution $\eta=\eta^{m+1}\in W^{2,4}(\Gamma)$ of \eqref{afterelew}. By a boot strapping argument one shows that $\eta=\eta^{m+1}\in W^{3,2}(\Gamma).$ Further $w=w^{m+1}\in W^{3,2}(\Gamma)$ is uniquely determined by the relation \eqref{subsubstructure}$_{1}.$ Indeed this regularity of $w^{m+1}$ is true for a fixed $\Delta t>0$ and not uniformly in $\Delta t.$ Hence there exists a couple $(\eta^{m+1},w^{m+1})\in W^{3,2}(\Gamma)\times W^{3,2}(\Gamma)$ solving \eqref{subsubstructure}.
		\subsubsection{Energy analogue at $\Delta t$ layer and convergence of interpolants:}
		Since $(\eta^{m+1}-\eta^{m})=\Delta t w^{m+1}\in W^{3,2}(\Gamma),$ we can use $\Delta t w^{m+1}$ as a test function in \eqref{subsubstructure}$_{2}$ in order to furnish the following estimate independent of $\Delta t:$
		\begin{equation}\label{testingsss}
			\begin{split}
				& \frac{(1-\delta)}{2}\|w^{m+1}\|^{2}_{L^{2}(\Gamma)}+\zeta\Delta t\|\nabla w^{m+1}\|^{2}_{L^{2}(\Gamma)}+\delta\|\nabla^{3}\eta^{m+1}\|^{2}_{L^{2}(\Gamma)}\\
				&\quad+\frac{h}{2}\int_{\Gamma}\mA\GG(\eta^{m+1}):\GG(\eta^{m+1})
				+\frac{h^{3}}{24}\int_{\Gamma}\mA\RR(\eta^{m+1}):\RR(\eta^{m+1})\\
				&\quad +\frac{\delta\Delta t}{2\tau}\|w^{m+1}-v^{n}\cdot\nu\|^{2}_{L^{2}(\Gamma)}+\frac{\delta\Delta t}{2\tau}\|w^{m+1}\|^{2}_{L^{2}(\Gamma)}
				\\
				& \leqslant \frac{(1-\delta)}{2}\|w^{m}\|^{2}_{L^{2}(\Gamma)}+\frac{\delta}{\tau}\Delta t\|v^{n}\|^{2}_{L^{2}(\Gamma)}+\frac{h}{2}\int_{\Gamma}\mA\GG(\eta^{m}):\GG(\eta^{m})\\
				&\quad+\frac{h^{3}}{24}\int_{\Gamma}\mA\RR(\eta^{m}):\RR(\eta^{m})+\delta\|\nabla^{3}\eta^{m}\|^{2}_{L^{2}(\Gamma)}+\frac{\delta\Delta t}{\tau}\|v^{n}\|^{2}_{L^{2}(\Gamma)}.
			\end{split}
		\end{equation}
		Now one can define piecewise constant interpolants in a standard manner. We recall that our goal is to construct a solution for \eqref{structuralpen}, we discretize $[n\tau,(n+1)\tau)$ in subintervals of length $\Delta t$ and define the interpolants as piecewise constant functions in $[n\tau+m\Delta t,n\tau+(m+1)\Delta t+1):$\\
		\begin{equation}\label{interpole}
			\begin{alignedat}{2}
				\eta^{M}(t)=&\eta(n\tau)=\eta^{n\tau}&&\text{ for } t\in [n\tau-\Delta t,n\tau),\\
				\eta^{M}(t)=&\eta^{m}&&\text{ for }t\in [n\tau+(m-1)\Delta t,n\tau+m\Delta t), m\in\mathbb{N}.
			\end{alignedat}
		\end{equation}
		The interpolant $w^{M}$ is defined as 
		\begin{equation}\label{wM}
			\begin{alignedat}{2}
				w^{M}(t)&=w(n\tau)=\eta_{1}^{n\tau}&&\text{ for }t\in [n\tau-\Delta t, n\tau],\\
				w^{M}(t)&=\frac{\eta^{M}(t)-\eta^{M}(t-\Delta t)}{\Delta t}&&\text{ for } t\in [n\tau+(m-1)\Delta t,n\tau+m\Delta t), m\in\mathbb{N},
			\end{alignedat}
		\end{equation}
		where we recall that the notations $\eta^{n\tau}$ and $\eta^{n\tau}_{1}$ was first introduced in the statement of Theorem \ref{labelstructuralsp}.\\
		In view of \eqref{testingsss} and a telescoping argument, the interpolants $\eta^{M}$ and $w^{M}$ solve
		\begin{equation}\label{energytypeinter}
			\begin{split}
				&\frac{(1-\delta)}{2}\|w^{M}(t)\|^{2}_{L^{2}(\Gamma)}+\zeta\int^{t}_{n\tau}\|\nabla w^{M}\|^{2}_{L^{2}(\Gamma)}+\delta\|\nabla^{3}\eta^{M}\|^{2}_{L^{2}(\Gamma)}\\
				&\quad+\frac{h}{2}\int_{\Gamma}\mA\GG(\eta^{M}):\GG(\eta^{M})
				+\frac{h^{3}}{24}\int_{\Gamma}\mA\RR(\eta^{M}):\RR(\eta^{M})
				\\
				&\quad+\frac{\delta}{2\tau}\int^{t}_{n\tau}\bigg(\|w^{M}-v^{n}\cdot\nu\|^{2}_{L^{2}(\Gamma)}+\|w^{M}\|^{2}_{L^{2}(\Gamma)}\bigg)\\
				&\leqslant \frac{(1-\delta)}{2}\|\eta^{n\tau}_{1}\|^{2}_{L^{2}(\Gamma)}+\frac{\delta}{\tau}\int^{t}_{n\tau} \|v^{n}\|^{2}_{L^{2}(\Gamma)}+\frac{h}{2}\int_{\Gamma}\mA\GG(\eta^{n\tau}):\GG(\eta^{n\tau})\\
				&\quad+\frac{h^{3}}{24}\int_{\Gamma}\mA\RR(\eta^{n\tau}):\RR(\eta^{n\tau})+\delta\|\nabla^{3}\eta^{n\tau}\|^{2}_{L^{2}(\Gamma)}+\frac{\delta}{2\tau}\int^{t}_{n\tau}\|v^{n}\|^{2}_{L^{2}(\Gamma)},
			\end{split}
		\end{equation}
		for $t\in[n\tau,(n+1)\tau+1).$\\
		Further the interpolants solve the following (in view of \eqref{subsubstructure}$_{2}$)
		\begin{equation}\label{eqsolveinter}
			\begin{split}
				&(1-\delta)\int^{t}_{n\tau}\int_{\Gamma}\frac{w^{M}(s)-w^{M}(s-\Delta t)}{\Delta t}b+\delta \int^{t}_{n\tau}\int_{\Gamma}\frac{w^{M}-v^{n}\cdot\nu}{\tau}b+\zeta\int^{t}_{n\tau}\int_{\Gamma}\nabla w^{M}\cdot\nabla b\\
				&+\int^{t}_{n\tau}\langle K'_{\delta}(\eta^{M},\eta^{M}(t-\Delta t)),b\rangle=0,
			\end{split}
		\end{equation}
		for $b\in L^{2}((n\tau,(n+1)\tau+1),W^{3,2}(\Gamma))$ and $t\in [n\tau,(n+1)\tau+1].$\\
		The bounds obtained from \eqref{energytypeinter}, infer the following weak type convergences (upto a non-relabeled subsequence)
		\begin{equation}\label{wktypconvweta}
			\begin{alignedat}{2}
				w^{M}&\rightharpoonup^{*} w^{n+1}&&\text{ in } L^{\infty}((n\tau,(n+1)\tau),L^{2}(\Gamma)),\\
				w^{M}&\rightharpoonup w^{n+1}&&\text{ in }  L^{2}((n\tau,(n+1)\tau),\sqrt{\zeta}W^{1,2}(\Gamma)),\\
				\eta^{M}&\rightharpoonup \eta^{n+1}&&\text{ in }  L^{\infty}((n\tau,(n+1)\tau);W^{3,2}(\Gamma)).
			\end{alignedat}
		\end{equation}
		We stress on the fact that we will not use the convergence \eqref{wktypconvweta}$_{2}$ in this proof. Hence this proof remains independent of the viscous nature of the structure and is valid for both hyperbolic Koiter shell and parabolic visco-elastic Koiter shell.\\
		We would now be interested to show the strong convergence of $\eta^{M}.$ For that we will verify the assertions of Aubin-Lions-Simons compactness theorem. In that direction we first define piece-wise affine interpolant $\widetilde{\eta}^{M}$ as
		$$\widetilde{\eta}^{M}=\frac{(n\tau+(m+1)\Delta t)-t}{\Delta t}\eta^{M}(t-\Delta t)+\frac{t-(n\tau +m\Delta t)}{\Delta t}\eta^{M}(t)\text{ for} t\in [n\tau+m\Delta t,n\tau+(m+1)\Delta t),m\in\mathbb{N}_{0}$$
		and observe that 
		$\partial_{t}\widetilde{\eta}^{M}(t)=w^{M}(t).$ Hence in view of \eqref{energytypeinter} we in particular have that 
		\begin{equation}\label{W12timebndetam}
			\widetilde{\eta}^{M}\,\,\mbox{is bounded in}\,\,  W^{1,2}((n\tau,(n+1)\tau+1),L^{2}(\Gamma)).
		\end{equation}
		Next we wish to estimate the difference $\eta^M(t+s)-\eta^M(t)$ for $t\in[n\tau,(n+1)\tau]$ and $s=\tilde k\Delta t$ for some $\tilde k\in\mathbb{N}_{0}.$ Indeed we are interested about small values of $s$ (hence small values of $\tilde k$) such that $(t+s)\in [n\tau, (n+1)\tau+1].$ Obviously, there is $\tilde m\in\eN_0$ such that $t\in [n\tau+\tilde m \Delta t,n\tau+(\tilde m+1)\Delta t)$ and $t+s\in [n\tau+(\tilde k+\tilde m) \Delta t,n\tau+(\tilde k+\tilde m+1)\Delta t)$. Then by the definitions of interpolants we obtain
		\begin{equation}\label{InterpTimeDiffIdent}
			\begin{array}{ll}
				\displaystyle\eta^M(t+s)-\eta^M(t)=\eta^{\tilde k+\tilde m+1}-\eta^{\tilde m+1}&\displaystyle=\widetilde \eta^M\left(n\tau+(\tilde k+\tilde m+1)\Delta t\right) - \widetilde \eta^M\left(n\tau+(\tilde m+1)\Delta t\right)\\
				&\displaystyle=\widetilde \eta^M(\tilde t+s)- \widetilde \eta^M\left(\tilde t\right)
			\end{array}
		\end{equation}
		for $\tilde t=n\tau+(\tilde m+1)h.$\\
		Using \eqref{InterpTimeDiffIdent}, the bound \eqref{W12timebndetam} and the embedding $W^{1,2}(0,T+1;L^{2}(\Gamma))\hookrightarrow C^{0,\frac{1}{2}}([0,T+1];L^{2}(\Gamma))$ we obtain 
		\begin{equation*} 
			\|\eta^M(t+\tilde s)-\eta^M(t)\|_{L^{2}(\Gamma)}=\|\widetilde \eta^M(\tilde t+\tilde s)-\widetilde \eta^M(\tilde t)\|_{L^{2}(\Gamma)} \leq c\tilde s^\frac{1}{2}
		\end{equation*}
		for $t\in [n\tau,(n+1)\tau+1-\tilde{s}]$ with $\tilde s=\tilde k\Delta t$, $\tilde k\in\eN$ and $\tilde{s}<\tau+1.$ Hence we conclude
		\begin{equation*} 
			\int_{n\tau}^{(n+1)\tau+1-\tilde s}\|\eta^M(t+\tilde s)-\eta^M(t)\|^2_{L^{2}(\Gamma)}\leq c(\tau)\tilde s
		\end{equation*}	  
		with $c$ independent of $M$. Then we find $z\in\eN$ such that $(n+1)\tau<zh\leq (n+1)\tau+1$. As a consequence of Lemma~\ref{Lem:RefIneq} we have \begin{equation*} 
			\int_{n\tau}^{(n+1)\tau+1-s}\|\eta^M(t+s)-\eta^M(t)\|^2_{L^{2}(\Gamma)}\leq c(\tau)s
		\end{equation*} 
		for any $0<s<(n+1)\tau+1.$\\
		Taking also into account \eqref{wktypconvweta} and the chain of embeddings $W^{3,2}(\Gamma)\stackrel{C}{\hookrightarrow} W^{2,4}(\Gamma)\hookrightarrow L^{2}(\Gamma)$ Lemma~\ref{Lem:RelComp}  yields the existence of a nonrelabeled subsequence $\{\eta^M\}$ such that 
		\begin{equation}\label{VNStrongly}
			\eta^M\to \eta^{n+1}\text{ in }L^2((n\tau,(n+1)\tau);W^{2,4}(\Gamma))\text{ as }M\to\infty.
		\end{equation}
		Using the boundedness of $\eta^{M}$ in $L^{\infty}((n\tau,(n+1)\tau);W^{3,2}(\Gamma)$ and the strong convergence \eqref{VNStrongly} one in particular renders that
		\begin{equation}\label{VNStrongly2}
			\eta^M\to \eta^{n+1}\text{ in }L^p((n\tau,(n+1)\tau);L^{\infty}(\Gamma))\text{ as }M\to\infty\,\,\mbox{for any}\,\,1<p<\infty.
		\end{equation}
		We next claim that
		\begin{equation}\label{wntmder}
			\displaystyle w^{n+1}=\partial_{t}\eta^{n+1}.
		\end{equation}
		To that end let us observe that
		\begin{equation}\nonumber
			\begin{array}{ll}
				&\displaystyle \widetilde{\eta}^{M}(t)-\eta^{M}(t)\\
				&\displaystyle =({t-(n\tau+(m+1)\Delta t)})\frac{\eta^{M}(t)-\eta^{M}(t-\Delta t)}{\Delta t}\,\,\mbox{when}\,\, t\in [n\tau+m\Delta t,n\tau+(m+1)\Delta t)
			\end{array}
		\end{equation} 
		which leads to 
		$$\|\widetilde{\eta}^{M}(t)-\eta^{M}(t)\|_{L^{2}(\Gamma)}\leqslant \Delta t\|w^{M}(t)\|_{L^{2}(\Gamma)}.$$
		The last estimate along with the bound of $w^{M}$ in $L^{2}(L^{2}(\Gamma))$ furnishes that
		$$\widetilde{\eta}^{M}-\eta^{M}\rightarrow 0\,\,\mbox{as}\,\, M\rightarrow\infty\,\,\mbox{in}\,\, L^{2}((n\tau,(n+1)\tau);L^{2}(\Gamma)).$$
		Since $\partial_{t}\widetilde{\eta}^{M}(t)=w^{M}(t),$ in view of the last convergence, we conclude the proof of \eqref{wntmder}.\\
		Next we use the relation $\partial_{t}\widetilde{\eta}^{M}=w^{M}$ and the boundedness of $w^{M}$ in $L^{2}(L^{2}(\Gamma))$ further to observe that
		\begin{equation}\label{weaketadelt}
			\begin{array}{ll}
				\displaystyle  \eta^M(\cdot-\Delta t)-\eta^{M}(\cdot)\rightarrow 0\,\,\mbox{in}\,\,L^{2}((n\tau,(n+1)\tau);L^{2}(\Gamma)).
			\end{array}
		\end{equation}
		The relation \eqref{weaketadelt}, along with the boundedness of both $\eta^{M}$ and $\eta^{M}(\cdot-\Delta t)$ in $L^{2}((n\tau,(n+1)\tau);W^{3,2}(\Gamma)$ and an application of interpolation argument furnishes the following strong convergence
		\begin{equation}\label{strngconvetaMlag}
			\begin{array}{ll}
				\displaystyle \eta^{M}(\cdot-\Delta t)\rightarrow \eta^{n+1}\,\,\mbox{in}\,\, L^{2}((n\tau,(n+1)\tau);W^{2,4}(\Gamma).
			\end{array}
		\end{equation}
		Indeed 
		\begin{equation}\label{strnfetalagLP}
			\begin{array}{l}
				\displaystyle \eta^{M}(\cdot-\Delta t)\rightarrow \eta^{n+1}\,\,\mbox{in}\,\, L^{p}((n\tau,(n+1)\tau);L^{\infty}(\Gamma),\,\,\mbox{for any}\,\, 1<p<\infty.
			\end{array}
		\end{equation}
		\subsubsection{Limit passage in \eqref{eqsolveinter} and \eqref{energytypeinter}}
		The obtained convergences in the last section, specially \eqref{wktypconvweta}$_{3}$, \eqref{VNStrongly}, \eqref{VNStrongly2}, \eqref{strngconvetaMlag} and \eqref{strnfetalagLP} are enough for the passage $M\rightarrow \infty$ in the approximation of the non-linear Koiter energy $\displaystyle\int^{t}_{n\tau}\langle K'_{\delta}(\eta^{M},\eta^{M}(t-\Delta t)),b\rangle$ (one recalls the definition of $\langle K'_{\delta}(\eta^{M},\eta^{M}(t-\Delta t)),b\rangle$ from \eqref{Kepsilon0b}). The other terms in \eqref{eqsolveinter} are linear in $w^{M}.$ Hence the passage $M\rightarrow \infty$ in the second and third terms of \eqref{eqsolveinter} is trivial. In order to pass  to the limit in the first term one observes that
		\begin{equation}\label{dertimewm}
			\begin{array}{ll}
				&\displaystyle \int^{t}_{n\tau}\int_{\Gamma}\frac{w^{M}(s)-w^{M}(s-\Delta t)}{\Delta t}b\\
				&\displaystyle =\int^{t-\Delta t}_{n\tau}\int_{\Gamma}w^{M}(s)\frac{b(s+\Delta t)-b(s)}{\Delta t}+\int^{t}_{t-\Delta t}\int_{\Gamma}\frac{w^{M}(s)}{\Delta t}b(s)-\frac{1}{\Delta t}\int^{n\tau+\Delta t}_{n\tau}\int_{\Gamma}\eta^{n\tau}_{1}b=\sum^{3}_{i=1}I_{i}.
			\end{array}
		\end{equation}
		As $\Delta t\rightarrow 0,$ (equivalently $M\rightarrow \infty$) one observes that
		$$I_{1}\rightarrow \int^{t}_{n\tau}\int_{\Gamma}w^{n+1}\partial_{t}b= \int^{t}_{n\tau}\int_{\Gamma}\partial_{t}\eta^{n+1}\partial_{t}b$$
		where we have used \eqref{wntmder}.\\
		Next 
		$$I_{2}\rightarrow \int_{\Gamma}w^{n}(t)b(t)$$
		and $$I_{3}\rightarrow -\int_{\Gamma}\eta^{n\tau}_{1}b(n\tau).$$
		Hence one obtains \eqref{structuralpen} by passing $M\rightarrow\infty$ in \eqref{eqsolveinter}.\\ 
		Finally using weak lower semi-continuity convex functionals and \eqref{wntmder} in \eqref{energytypeinter} we furnish \eqref{energytypestr}.
	\end{proof}
	
	\subsection{Proof of Lemma \ref{Lem:Fund}}\label{extension0proof}
	\begin{proof}
		We adapt the level set approach used in the proof of \cite[Lemma 4.1]{FKNNS}.
		Let us choose a function $g_0\in C^\infty(B)$ such that
		\begin{equation*}
			g_0(x)=\begin{cases}
				=0&\text{ if }x\in \partial B\cup \partial\Omega,\\
				>0&\text{ if }x\in B\setminus\Omega,\\
				<0&\text{ else}
			\end{cases}
		\end{equation*}
		and 
		\begin{equation}\label{InG}
			\nabla g_0(x)=h(d(x))\nu(\pi(x))\end{equation}
		in a sufficiently small neighborhood $S$ of $\partial\Omega$,
		where the signed distance function $d$ and the projection of $x\in S$ to a closest point of $\partial\Omega$ to $x$ are defined in section \ref{sec:Geomtry} and $\inf_{\eR}h>0$.
		We consider the function $\tilde\varphi_\eta$ from \eqref{FlowMDef} and define $V(t,x)=\tder\tilde\varphi_\eta(t,(\tilde\varphi_\eta)^{-1}(t,x))$ and $g(t,x)=g_0((\tilde\varphi_\eta)^{-1}(t,x))$. Obviously, it follows that
		\begin{equation*}
			\tder(\tilde\varphi_\eta)^{-1}(t,x)=-(\nabla\tilde\varphi_\eta)^{-1}(t,(\tilde\varphi_\eta)^{-1}(t,x))\tder\tilde\varphi_\eta(t,(\tilde\varphi_\eta)^{-1}(t,x)),\ 
			\nabla (\tilde\varphi_\eta)^{-1}(t,x)=(\nabla\tilde\varphi_\eta)^{-1}(t,(\tilde\varphi_\eta)^{-1}(t,x)).
		\end{equation*}
		Accordingly, we infer that $g$ satisfy the transport equation
		\begin{equation}\label{TrEq}
			\tder g+V\cdot\nabla g=0\text{ in }(0,T)\times\eR^3.
		\end{equation}
		We note that the set $B\setminus\Omega_\eta(t)$ corresponds to $\{g(t,\cdot)>0\}$ and the interface $\Sigma_\eta(t)$ corresponds to the set $\{g(t,\cdot)=0\}$.
		Fixing $\xi>0$ and setting $\psi=\max\{\min\{\frac{1}{\xi}g,1\},0\}$ in \eqref{contrhoex}, which is possible via an approximating procedure, we get 
		\begin{equation}\label{TrG}
			\int_{B\setminus\Omega_{\eta}(t)}\rho\psi=\frac{1}{\xi}\int_0^t\int_{\{0\leq g(\tau,x)<\xi\}}(\rho\tder g+\rho u\cdot\nabla g).
		\end{equation}
		Employing \eqref{TrEq} we infer
		\begin{equation*}
			\rho(\tder g+u\cdot\nabla g)=\rho(u-V)\cdot\nabla g.
		\end{equation*}
		Using the latter identity on the right hand side of \eqref{TrG} we obtain
		\begin{equation}\label{TrGFin}
			\int_{B\setminus\Omega_{\eta}(t)}\rho\psi=\frac{1}{\xi}\int_0^t\int_{\{0\leq g(\tau,x)<\xi\}}\rho(u-V)\cdot\nabla g.
		\end{equation}
		We focus on the regularity of the expression $(u-V)\cdot\nabla g$. By the assumed regularity of $u$ and the definition of $V$ and the assumed regularity of $\eta$
		we deduce 
		\begin{equation}\label{UVDifReg}
			u-V\in L^2(0,T;W^{1,2}(B)).
		\end{equation}
		Employing the assumed regularity of $\eta$, the regularity of the given mapping $\varphi$ and the regularity of the projection $\pi$ accordingly, we conclude from the definition of $g$ and \eqref{FlowMInv}
		\begin{equation}\label{NGReg}
			\nabla g\in L^\infty(0,T;L^\infty(B))\cap L^\infty(0,T;W^{1,2}(B)).
		\end{equation}
		Hence using the Sobolev embedding we infer from \eqref{UVDifReg} and \eqref{NGReg} that for a.a. $t\in(0,T)$
		\begin{equation*}
			(u-V)\cdot\nabla g\in L^2(0,T;W^{1,\frac{3}{2}}(B)).
		\end{equation*}
		Moreover, we know that $\Tr(u-V)=0$ on $\partial B\cup\Sigma_{\eta}(t)$. Hence applying the Hardy inequality we get 
		\begin{equation}\label{HIn}
			\left\|\frac{(u-V)\cdot\nabla g}{\dist(\cdot,\partial B\cup\partial\Sigma_{\eta}(t))}\right\|_{L^\frac{3}{2}(B\setminus\Omega_{\eta}(t))}\leq c\|(u-V)\cdot\nabla g\|_{W^{1,\frac{3}{2}}(B\setminus\Omega_\eta(t))}.
		\end{equation}
		We note that the constant in the Hardy inequality depends also on the Lipschitz constant of $\eta(t)$ that can be estimated uniformly in time due to the assumed regularity of $\eta$. Hence the constant $c$ in \eqref{HIn} can be taken independent of $t$. Using \eqref{HIn}, the assumed regularity of $\rho$  we conclude from \eqref{TrG}
		\begin{equation}\label{AuxI}
			\begin{split}
				\int_{B\setminus\Omega_{\eta}(t)}(\rho\psi)(t)&\leq\xi^{-1}\left|\int_0^t\int_{\{0\leq g(\tau, x)<\xi\}}\rho(u-V)\cdot\nabla g\right|\\
				&\leq T^\frac{1}{2}\sup_{(t,x)\in M,\xi\in(0,\xi_0]}F(t,x,\xi)\|\rho\|_{L^\infty(0,T;L^3(\{0\leq g(\tau,\cdot)<\xi\}))}\|(u-V)\cdot\nabla g\|_{L^2(0,T;W^{1,\frac{3}{2}}(B))},
			\end{split}
		\end{equation}
		where $M=\bigcup_{t\in[0,T]}\{t\}\times \{x\in B\setminus\Omega_\eta(t):0\leq g(t,x)<\xi\}$ and $F(t,x,\xi)=\xi^{-1}\dist(x,\partial B\cup\Sigma_\eta(t))$. The choice of $\xi_0$ is specified in the following way. The number $\xi_0$ is chosen small such that $g(t,x)<\xi_0$ implies one of the following options. The first one is that $x$ belongs to a neighborhood $N$ of $\partial B$ on which $\tilde\varphi_\eta$ is the identity and $\min_{x\in\overline N}|\nabla g_0(x)\cdot\nu(\pi_{\partial B}(x)|>0$, where $\pi_{\partial B}(x)$ is the projection of $x\in \overline N$ on $\partial B$ such that $|x-\pi_{\partial B}(x)|=\dist(x,\partial B)$. The second option is that $(\tilde\varphi_\eta)^{-1}(t,x)\in S$. The next task is to show that
		\begin{equation}\label{SupCond}
			\sup_{(t,x)\in M,\xi\in(0,\xi_0]}F(t,x,\xi)<\infty.
		\end{equation}
		To this end we distinguish the cases $\dist(x,\partial B\cup\Sigma_\eta(t))=\dist(x,\partial B)$ and $\dist(x,\partial B\cup\Sigma_\eta(t))=\dist(x,\Sigma_\eta(t))$.  In the first case we have for fixed $\xi\leq\xi_0$ and any $x\in N$
		\begin{equation}\label{XiFir}
			\xi\geq g(t,x)=g_0(x)-g(\pi_{\partial B}(x))\geq \dist (x,\partial B)\min_{x\in\overline N}|\nabla g_0(x)\cdot\nu(\pi_{\partial B}(x))|
		\end{equation}
		implying \eqref{SupCond} immediately. Concerning the second case we have
		\begin{equation}\label{XISec}
			\xi\geq g(t,x)=g(t,x)-g(t,\tilde\varphi_\eta(t,\pi(x)))\geq \dist (x,\Sigma_\eta(t) )\min_{(t,x)\in O}|\nabla g(t,x)\cdot\nu(\pi(x))|,
		\end{equation}
		where $O=\bigcup_{t\in[0,T]}\{t\}\times\tilde\varphi_\eta(t,\overline{S})$. 
		Taking into account \eqref{InG} and \eqref{FlowMInv} we have in $O$ 
		\begin{equation}\label{NGEx}
			\begin{split}
				\nabla g(t,x)=&\nabla(\tilde\varphi_\eta)^{-1}(t,x)\nabla g_0((\tilde\varphi_\eta)^{-1}(t,x))=\nabla(\tilde\varphi_\eta)^{-1}(t,x)h(d((\tilde\varphi_\eta)^{-1}(t,x)))\nu(\pi((\tilde\varphi_\eta)^{-1}(t,x)))\\
				=&h(d((\tilde\varphi_\eta)^{-1}(t,x)))\partial_{\nu(\pi(x))}(\tilde\varphi_\eta)^{-1}(t,x)=h(d((\tilde\varphi_\eta)^{-1}(t,x)))(1-f'_\Gamma(d(x))\eta(t,\varphi^{-1}(\pi(x))))\nu(\pi(x)),
			\end{split}
		\end{equation}
		denoting by $\partial_{\nu(\pi(x))}$ the derivative in the direction $\nu(\pi(x))$.
		Noticing that $\min\{m,0\}\leq \eta\leq\max\{0,M\}$ in $[0,T]\times\Gamma$ we have due to \eqref{fPrEst}
		\begin{equation*}
			1-f'_\Gamma(d(x))\eta(t,\varphi^{-1}(\pi(x)))\geq 1-\max\left\{\frac{\max\{M,0\}}{M'},\frac{\min\{m,0\}}{m'}\right\}\text{ in }O.
		\end{equation*}
		Hence combining the latter inequality with \eqref{NGEx}  we obtain 
		\begin{equation*}
			\min_{(t,x)\in O}|\nabla g(t,x)\cdot\nu(\pi(x))|>0.
		\end{equation*}
		This along with \eqref{XISec} concludes \eqref{SupCond}. Moreover, it follows from \eqref{SupCond} that for $\xi$ small enough we get
		\begin{equation}\label{MesEst}
			|\{x\in B\setminus\Omega_\eta(t):0\leq g(t,x)<\xi\}|\leq c\xi
		\end{equation}
		with $c$ independent of $t\in[0,T]$. Using \eqref{SupCond}, \eqref{MesEst} and the assumption $\rho\in L^\infty(0,T;L^3(B))$ we pass to the limit $\xi\to 0_+$ in \eqref{AuxI} to conclude
		\begin{equation*}
			\int_{B\setminus\Omega_\eta(t)}\rho(t,\cdot)=0
		\end{equation*}
		implying $\rho(t)|_{B\setminus\Omega_\eta(t)}\equiv 0$ for a.a. $t\in(0,T)$. The conclusion for $Z$ is obtained in the exactly same way.
	\end{proof}
	\subsection{Comments on the proof of Lemma \ref{strngconvdteta}}\label{lempfstrnged}
	The convergence \eqref{DLimIdent}$_{2}$ is a consequence of \eqref{AddConvDel}$_5$ and 
	\begin{equation}\label{DLimIdent}
		\begin{split}
			\lim_{\delta\to 0_+}\left(\int_0^T\int_\Gamma |\tder \eta^\delta|^2+\int_0^T\int_{\Omega_{\eta^\delta}(t)}(\rho^\delta +Z^\delta)u^\delta\cdot\mathcal F_{\eta^\delta}\tder\eta^\delta\right)=&\int_0^T\int_\Gamma |\tder \eta|^2+\int_0^T\int_{\Omega_{\eta}(t)}(\rho^\delta +Z^\delta)u\cdot\mathcal F_{\eta}\tder\eta,\\
			\lim_{\delta\to 0_+}\int_0^T\int_{\Omega_{\eta^\delta}(t)}(\rho^\delta+Z^\delta) u^\delta\cdot(u^\delta-\mathcal F_{\eta^\delta}\tder\eta^\delta)=&\int_0^T\int_{\Omega_{\eta}(t)}(\rho+Z) u\cdot(u-\mathcal F_{\eta}\tder\eta),
		\end{split}
	\end{equation}
	where $\mathcal{F}_{\eta}$ is introduced in \eqref{FEDef}. In order to show \eqref{DLimIdent}$_2$, we note that $\|\mathcal F_{\eta^\delta} \tder\eta^\delta\|_{L^2(0,T;W^{1-\frac{1}{r},p}(B))}\leq c\|\tder\eta^\delta\|_{L^2(0,T;W^{1-\frac{1}{r}},r(\Gamma))}$, $p\in [1,\frac{3r}{2})$ as follows by \cite[Lemma 2.7(a)]{Breit2} and the uniform bound on $\{\tder\eta^\delta\}$ in $L^2(0,T;W^{1-\frac{1}{r},r}(\Gamma))$ for any $r\in[1,2)$ following from the coupling $\tder\eta\nu=\Tr_{\Sigma_{\eta^\delta}}u^\delta$, the bound \eqref{UestOnB} and Lemma \ref{Lem:TrOp}. Hence we get the compactness of $\{\mathcal F_{\eta^\delta}\}$ in the weak topology of $L^2(0,T;W^{\sigma,p}(B))$ for any $\sigma\in[0,\frac{1}{2})$, $p\in[1,3)$. Using the linearity of $\mathcal F_\Omega$,  convergences \eqref{weaklimiturhoeta}$_{1,2}$ and \eqref{EtaDHCnv} we conclude from definition \eqref{FEDef} that up to a nonrelabeled subsequence
	\begin{equation}\label{ExtTDerWCnv}
		\mathcal F_{\eta^\delta}\tder\eta^\delta\rightharpoonup \mathcal F_{\eta}\tder\eta\text{ in }L^2(0,T;W^{\sigma,p}(B))\text{ for any }\sigma\in[0,\frac{1}{2}),\ p\in[1,3).
	\end{equation} 
	Next, $W^{\sigma,p}(B)$ with $\sigma\in[0,\frac{1}{2})$, $p\in[1,3)$ is compactly embedded in $L^s(B)$ with $s<6$ implying $L^{s'}(B)$ with $s'>\frac{6}{5}$ is compactly embedded in $\left(W^{\sigma,p}(B)\right)'$. Therefore we get 
	\begin{equation*}
		(\rho^\delta+Z^\delta)u^\delta\to  (\rho+Z)u\text{ in }L^2(0,T;(W^{\sigma,p}(B))')\text{ for any }\sigma\in [0,\frac{1}{2}),\ p\in[1,3)
	\end{equation*}
	from \eqref{AddConvDel}$_3$ as $\frac{2\max\{\gamma,\beta\}}{\max\{\gamma,\beta\}+1}>\frac{6}{5}$. The latter convergence and \eqref{ExtTDerWCnv} concludes \eqref{DLimIdent}$_2$. 
	Identity \eqref{DLimIdent}$_1$ follows by making use of the general compactness result \cite[Theorem 5.1. and Remark 5.2.]{MuhaSch}. We jus mention that the justification of the assumption of \cite[Theorem 5.1.]{MuhaSch} is performed in \cite[Section 4.3]{Breit2}. In fact, this justification can be easily adapted in our case which is even simpler because the momentum equation is not considered at the Galerkin level and there is no need to project into discrete spaces when justifying the equi--continuity assumption. We notice that the key ingredient used in this justification is the bound on $\{\nabla \eta^\delta\}$ in $L^2(0,T;L^\infty(\Gamma))$ that follows immediately by \eqref{improvedregeta} in Lemma \ref{improvebndetadelta}.
	
	\begin{center}
		\Large\textbf{Acknowledgements} \\[4mm]
	\end{center}
	
	\textit{This work has been supported by the Czech Science Foundation (GA\v CR) through projects 22-08633J (for \v S.N. and M.K.) Moreover, \it \v S. N., M.K. and S. M.  have been supported by  Praemium Academiæ of \v S. Ne\v casov\' a. Finally, the Institute of Mathematics, CAS is supported by RVO:67985840.}
	%%%%%%%%%%%%%%%%%%%%%%%%%%%%%%%%%%%%%%%%%%%%%%%%%%%%%%%%%%%%%%%%%%%%%%%%%%%%%%%%%%%%%%%%%%%%%%%%%%%%%%%%%%%%%%%%%%%%%%%%%%%%%%%%%%%%%%%%%%%%%%%%%%%%%%%%%%%%%%%%%%%%%%%%%%%%%%%%%%%%%%%%%%%%%%%%%%%%%%%%%%%%%%%%%%%%%%%%%%%%%%%%%%%%%%%%%%%%%%%%%%%%%%%%%%%%%%%%%%%%%%%%%%%%%%%%%%%%%%%%%%%%%%%%%%%%%%%%%%%%%%%%%%

\end{document}